\documentclass{acta_2011}
\usepackage{harvard_acta}
\usepackage{amssymb}
\usepackage{amsmath} 
\usepackage{graphicx}
\usepackage{subfigure}
\def \RR   {\mathbb{R}}
\def \RN   {\RR^N}  
\def \del  {\partial}    

\def \Ucal  {\mathcal U}

\def \lam  {\lambda} 
\def \eps  {\epsilon} 
\def \be   {\begin{equation}}
\def \ee   {\end{equation}}
\def \bea  {\begin{eqnarray}} 
\def \eea  {\end{eqnarray}}
\def \Dx {{\Delta x}}
\def \Ord {{\mathcal O}}
\def \bei {\begin{itemize}}
\def \eei {\end{itemize}} 
\newcommand \bel {\be\label} 

\begin{document}

\title[Numerical methods for small-scale dependent shocks]{Numerical methods
\\
with controlled dissipation 
\\
for small-scale dependent shocks}
\author[{Acta Numerica}]{%
Philippe G. LeFloch\\
 Laboratoire Jacques-Louis Lions 
\\
and Centre National de la Recherche Scientifique
\\
Universit\'e Pierre et Marie Curie
\\
4 Place Jussieu, 75258 Paris, France.
\\
{\tt contact@philippelefloch.org}
\\
\and
Siddhartha Mishra\\ 
Seminar for Applied Mathematics, D-Math
\\
Eidgen\"ossische Technische Hochschule
\\
HG Raemistrasse, Z\"urich--8092, Switzerland. 
\\
{\tt siddhartha.mishra@sam.math.ethz.ch}
}

\pagerange{\pageref{firstpage}--\pageref{lastpage}}
\doi{}
%XXXXXX}
\label{firstpage}

\maketitle

 \vskip1.cm 

\begin{abstract} We provide a `user guide' to the literature of the past twenty years concerning the modeling and  approximation of discontinuous solutions to nonlinear hyperbolic systems that admit {\sl small-scale dependent} shock waves. We cover several classes of problems and solutions: nonclassical undercompressive shocks, hyperbolic systems in nonconservative form, boundary layer problems. We review the relevant models arising in continuum physics and describe the numerical methods that have been proposed to capture small-scale dependent solutions. In agreement with the general well-posedness theory, small-scale dependent solutions are characterized by a {\sl kinetic relation,} a {\sl family of paths,} or an {\sl admissible boundary set.}  We provide a review of numerical methods (front tracking schemes, finite difference schemes, finite volume schemes), which, at the discrete level, reproduce the effect of the physically-meaningful dissipation mechanisms of interest in the applications. An essential role is played by the {\sl equivalent equation} associated with discrete schemes, which is found to be relevant even for solutions containing shock waves. 
\end{abstract}

\setcounter{tocdepth}{1}
\tableofcontents

%=================================================================================== 

\section{Introduction}
\label{sec:1}

\subsection{Small-scale dependent shock waves}
\label{sec:11}

Nonlinear hyperbolic systems of partial differential equations, arising in continuum physics and especially in compressible fluid dynamics, 
admit discontinuous solutions containing shock waves that may {\sl depend on underlying small--scale mechanisms} such as the coefficients of viscosity, capillarity, Hall effect, relaxation, heat conduction, etc. 
Such small-scale dependent shocks exist for a variety of problems of physical interest, for instance, those modeled by conservative hyperbolic systems with dispersive phenomena (e.g.~with capillary effects) and nonconservative hyperbolic systems (e.g.~for two-phase fluid flows), as well as boundary layer problems (e.g.~with viscosity terms). In the past twenty years, it is has been successively recognized that a standard entropy inequality (after Lax, Oleinik, Kruzkov, and others) does not suffice for the unique characterization of physically meaningful solutions to such problems, so that {\sl additional criteria} are required in order to characterize these small-scale dependent shock waves uniquely. These criteria are based on the prescription of a \emph{kinetic relation} for conservative hyperbolic systems, a \emph{family of paths} for nonconservative hyperbolic systems, and an \emph{admissible boundary set} for boundary value problems. 

Although standard finite difference, finite volume and finite element methods have been very successful in computing solutions to hyperbolic conservation laws, including those containing shock waves, these well-established methods are found to be inadequate for the approximation of small-scale dependent shocks. Starting with Hou and LeFloch (1994) and Hayes and LeFloch (1996), this lack of convergence was explained in terms of the equivalent equation associated with discrete schemes through a formal Taylor expansion. The leading terms in the equivalent equation represent the numerical viscosity of the scheme and need not match the physically relevant small--scale mechanisms that have been neglected in the hyperbolic modeling.  Consequently, standard shock-capturing schemes, in general, fail to converge to physically meaningful solutions. 

Our purpose here is to review the methods developed in the past twenty years which accurately describe and compute small--scale dependent shocks. The challenge is, both, to develop the proper theoretical tools for the description of such solutions and to ensure the convergence of numerical methods toward physically meaningful solutions. We built here on several earlier reviews by LeFloch (1989, 1999, 2002, 2010) and include the most recent developments and material on numerical methods that were not presented until now ---especially the theory of schemes with well-controlled dissipation (WCD, in short) recently developed by the authors. 

%------------------------------------------------------------------------------------------------------------------ 

\subsection{Physical models and mathematical theories}
\label{sec:12}

In one space dimension, the classes of hyperbolic systems under review either admit the conservative form 
\bel{eq:109}
u_t + f(u)_x = 0
\ee 
or the nonconservative form 
\bel{eq:111}
u_t + A(u) \, u_x = 0, 
\ee 
in which the map $u: \RR^+ \times \RR \to \Ucal$ is the unknown. 
In \eqref{eq:109}, the given flux $f: \Ucal \to \RN$ is defined on a (possibly non-connected) open set $\Ucal\subset \RN$ 
and satisfies the following strict hyperbolicity condition: 
for every $v \in \Ucal$, the matrix $A(v):=Df(v)$ admits real and distinct eigenvalues 
$\lam_1(v) < \ldots < \lam_N(v)$ and a basis of right-eigenvectors $r_1(v), \ldots, r_N(v)$.  
In \eqref{eq:111}, the given matrix-valued field $A=A(u)$ need not be a Jacobian matrix and, also, is required to possess real and distinct eigenvalues and a complete set of eigenvectors.

Small-scale dependent solutions arise with systems of the form \eqref{eq:109} or \eqref{eq:111} and it is our objective to present and investigate specific models of interest in physical applications, especially: 
\begin{itemize}

\item {\sl Strictly hyperbolic systems.} The simplest example of interest is provided by the scalar conservation law with cubic flux and added second- and third-order terms and is presented in Section~\ref{sec:21} below. More challenging models arise in the dynamics of fluids and nonlinear elastic material with viscosity and capillarity effects (cf.~Section~\ref{sec:23}). 

\item {\sl Non-strictly hyperbolic systems.} In presence of certain phase transition phenomena, the models of fluids and elastic materials fail to be globaly strictly hyperbolic. Furthermore, the system of magnetohydrodynamics with viscosity and Hall effects is probably the most challenging model and plays also an essential role in applications, for instance, in the modeling of the solar wind.

\item {\sl Nonconservative hyperbolic systems.} As will be discussed in Section~\ref{sec:5} below, one of the simplest (yet challenging) example in the class of nonconservative hyperbolic systems is obtained by coupling two Burgers equations, while the most challenging model of interest in the application contains five equations (or seven equations, if two thermodynamical variables are introduced) for the evolution of two (fluid and vapor) phases of a fluid mixture.
Other important models are the multi-layer shallow water system and the Lagrangian gas dynamics with internal energy taken as an independent variable.  

\item {\sl Initial and boundary value problems.} Models of interest include, both, the linearized and nonlinear Euler equations with artificial viscosity or physical viscosity.

\end{itemize}
 
A (mathematical) entropy inequality can be naturally associated with all (conservative or nonconservative) systems under consideration, that is, 
\bel{eq:109-ineq}
U(u)_t + F(u)_x \leq 0, 
\ee 
where $U: \Ucal \to \RR$ and $F: \Ucal \to \RR^n$ are refered to as the entropy and entropy-flux, respectively.  
However, in constrast with more classical problems arising in fluid dynamics, \eqref{eq:109-ineq} 
is often insufficiently discriminating in order to characterize physically meaningful solutions to the initial value problem associated with 
\eqref{eq:109} or \eqref{eq:111}. Therefore, additional admissibility criteria are required, as follows: 
\begin{itemize}

\item {\sl Kinetic relations} $\Phi_j =\Phi_j(u_-)$ provide a general tool to define nonclassical entropy solutions to strictly hyperbolic systems of conservation laws (LeFloch 1993, Hayes and LeFloch 1996) and are relevant when the characteristic fields of the systems do not satisfy Lax's genuine nonlinearity condition (Lax 1957, 1973). They were introduced first for an hyperbolic-elliptic model of phase transitions in solids (Abeyaratne and Knowles 1991a, 1991b, Truskinovsky 1987, 1993, 1994). Roughly speaking, a kinetic relation for a nonclassical (undercompressive) shock in the $j$-characteristic family prescribes the right-hand state $\Phi_j(u_-)$ as a function of the left-hand state $u_-$. 

\item {\sl Families of paths} $s \in [0,1] \mapsto \varphi(s; u_-, u_+)$  (Dal Maso, LeFloch, and Murat 1990, 1995) provide  underlying integration paths which are necessary in order to define generalized jump relations and weak solutions to nonconservative hyperbolic systems. The paths $s \in [0,1] \mapsto \varphi(s; u_-, u_+)$ connect left-hand states $u_-$ to right-hand states $u_+$ and are derived by analyzing the trajectories of traveling solutions, once an augmented model taking higher-order effects into account is selected. 

\item {\sl Admissible boundary sets} $\Phi(u_B)$ (Dubois and LeFloch 1988) are required in order to formulate well-posed initial-- boundary value problems associated with nonlinear hyperbolic systems. Waves propagating in weak solutions may collapse near a boundary and generate a boundary layer connecting a given boundary state $u_B$ with an actual state constrained to lie in a prescribed boundary set $\Phi(u_B)$. 
\end{itemize}

%------------------------------------------------------------------------------------------

\subsection{Designing schemes with well-controlled dissipation}  
\label{sec:13}

It is well now recognized (Hou and LeFloch 1994, Hayes and LeFloch 1996, LeFloch and Rohde 2000, LeFloch and Mohamadian 2008, 
Fjordholm and Mishra 2012) that a finite difference or volume scheme (LeVeque 2003) may not converge to physically-relevant weak solutions, unless one can ensure a certain {\sl consistency property with small-scale effects}, that is, a consistency with the prescribed kinetic relation, family of paths, or admissible boundary set associated with any specific problem under consideration. It was suggested in LeFloch (2010) to call the numerical methods satisfying this requirement as the schemes with ``controlled dissipation'' and, recent work by the authors (covering the treatment of shocks with arbitrary strength), it was proposed to refer to them as {\sl schemes with well-controlled dissipation}. 

In particular, the role of the {\sl equivalent equation} associated with a given finite difference scheme has been emphasized and analyzed. The leading terms of the equivalent equations for standard finite difference schemes like the Lax-Friedrichs scheme, for instance, significantly differ from the physically-relevant small-scale mechanisms. For small-scale dependent shocks, the approximate solutions, say $u^\Dx$ converge (when the discretization parameter $\Dx$ approaches zero) toward a limit, say $v$, which is {\sl distinct from the physical solution,} say $u$. In other words, standard finite difference or finite volume techniques, in general, lead to the non-convergence property   
\be
v := \lim_{h \to 0} u^h \neq u
\ee
This holds for a variety of strictly hyperbolic systems, nonconservative hyperbolic systems, and boundary layer problems. 

On the other hand, `schemes with well-controlled dissipation' are precisely designed to overcome this challenge and are built by analyzing discrete dissipation operators arising in equivalent equations, imposing that the latter should match the small-scale mechanisms in the underlying augmented system, at least to leading order. It was emphasized by Hou and LeFloch (1994) and Hayes and LeFloch (1996) that, in fact, the objective need not be to ensure the convergence of the schemes, but rather to control the error term (in a suitable norm)
$$
\| v - u \| 
$$
in terms of the physical parameters arising in the problem: shock strength, order of accuracy of the scheme, ratio of capillarity over viscosity, etc. 

In recent years, extensive studies have demonstrated the relevance of the equivalent equation, as a tool for designing numerical methods for computing small-scale dependent shocks, and have included numerical experiments in physically realistic set-ups. 

This paper is organized as follows.  Section~2 is concerned with strictly hyperbolic systems and the discussion of nonclassical undercompressive shocks to such systems, when the associated augmented models contain diffusive and dispersive terms. Then, Section~3 discusses the numerical methods adapted to these problems. Next, in Section~4 we turn our attention to nonconservative hyperbolic systems. The boundary value problems are discussed in Section~5 and some concluding remarks are made in Section~6. 

%====================================================================================

\section{Nonclassical entropy solutions to nonlinear hyperbolic systems}
\label{sec:2} 

\subsection{The regime of balanced diffusion and dispersion}
\label{sec:21}

We consider first the class of {\sl hyperbolic conservation laws with vanishing diffusion,} i.e. 
\be
\label{eq:356}
u^\eps_t + f(u^\eps)_x = \eps \, \bigl( b(u^\eps) \, u^\eps_x\bigr)_x, 
\ee
where $ u^\eps=u^\eps(t,x) \in \RR,$ is the unknown, the flux $f: \RR \to \RR$ is a given smooth function,
and the diffusion coefficient $b: \RR \to (0, \infty)$ is bounded above and below. 
For any given initial data and, solutions to the initial value problem 
associated with \eqref{eq:356} converge strongly (when $\eps \to 0$) toward a limit $u=u(t,x)$ satisfying the hyperbolic conservation law 
\be
\label{eq:356-0}
u_t + f(u)_x = 0   
\ee
in the weak sense of distributions. Weak solutions to \eqref{eq:356-0} are not uniquely characterized by their initial data, but must also be constrained to satisfy a certain entropy condition, ensuring that they be achieved as limits of \eqref{eq:356} (Oleinik 1963, Kruzkov 1970, Volpert 1967).  

More precisely, solutions $u^\eps$ to \eqref{eq:356} satisfy, for every convex function $U: \RR \to \RR$, 
$$
\aligned
& U(u^\eps)_t  + F(u^\eps)_x = - D^\eps + \eps \, C^\eps_x, 
\\
& D^\eps := \eps \, b(u^\eps) \, U''(u^\eps) \, |u^\eps_x|^2, 
 \qquad C^\eps := b(u^\eps) \, U(u^\eps)_x, 
\endaligned
$$
in which $F(u) := \int^u f'(v) \, U'(v) \, dv$ and $(U,F)$ is refered to as an entropy-entropy flux pair. 
Hence, $u = \lim_{\eps \to 0} u^\eps$  satisfies the so-called {\rm entropy inequalities} 
\be
\label{eq:357}
U(u)_t + F(u)_x \leq 0, \qquad U'' \geq 0. 
\ee 
Weak solutions to \eqref{eq:356-0} that satisfy {\sl all} inequalities \eqref{eq:357} (i.e for every convex $U$) 
are refered to as {\sl classical entropy solutions.} It is also customary to reformulate the entropy condition  in the Kruzkov form: 
$|u-k|_t + \big( \text{sgn}(u-k) (f(u) - f(k)) \big)_x \leq 0$ for all $k \in \RR$.

Nonclassical shock waves arise in weak solutions when, both, diffusion and dispersion are included. The simplest model of interest is provided by the {\sl linear diffusion-dispersion model}   
\be
\label{eq:cons2}
u^\eps_t + f(u^\eps)_x = \eps \, u^\eps_{xx} + \gamma(\eps) \, u^\eps_{xxx}, 
\ee 
which depends upon two parameters $\eps$ and $\gamma$ refered to as the diffusion and the dispersion coefficients. This equation was studied first by Jacobs, McKinney, and Shearer (1993), Hayes and LeFloch (1996, 1997), and Bedjaoui and LeFloch (2002abc, 2004).  Importantly, the relative scaling between $\eps$ and $\gamma=\gamma(\eps)$ determines the limiting behavior of solutions, and we can distinguish between three cases: 
\begin{itemize} 

\item In the {\sl diffusion-dominant regime} $\gamma(\eps) << \eps^2$, the qualitative behavior of solutions to \eqref{eq:cons2} is 
similar to the behavior of solutions to \eqref{eq:356} and, in fact, the limit 
$u:= \lim_{\eps \to 0} u^\eps$ is then {\sl independent} upon $\gamma(\eps)$, and is a {\sl classical} entropy solution characterized by the infinite family of entropy inequalities  \eqref{eq:357}.

\item In the {\sl dispersion-dominant regime} $\gamma(\eps) >> \eps^2$, high oscillations develop (as $\eps$ approaches zero) especially in regions of steep gradients of the solutions and only weak convergence of $u^\eps$ is observed. The vanishing dispersion method developed by Lax and Levermore (1983) for the Korteweg-de Vries equation is the relevant theory in this regime, which can not be covered by the techniques under consideration in the present review.

\item In the regime of {\sl balanced diffusion-dispersion,} corresponding typically to $\gamma(\eps) := \delta \, \eps^2$ for a fixed $\delta$, the limit $u:= \lim_{\eps \to 0} u^\eps$ does exist (in a strong topology) and only {\sl mild oscillations} are observed near shocks, so that the limit is a weak solution to the hyperbolic conservation law \eqref{eq:356-0}. Most importantly, when $\delta >0$, 
the solutions $u$ exhibit nonclassical behavior, as they contain undercompressive shocks (as defined below) and {\sl strongly depend} upon the  coefficient $\delta$.   
\end{itemize}

From now on, we focus our attention on the ``critical regime'' where the small-scale terms are kept in balance and, for instance, we write \eqref{eq:cons2} as
\bel{eq:240}
u^\eps_t + f(u^\eps)_x = \eps \, u^\eps_{xx} + \delta \, \eps^2 \, \, u^\eps_{xxx},
\ee 
in which $\delta$ is a fixed parameter and $\eps \to 0$. For this augmented model, we easily derive the identity 
$$
\aligned
&  (1/2) \, | u_\eps^2 |_t  + F(u^\eps)_x
= - D^\eps + \eps \, C^\eps_x, 
\\
& D^\eps := \eps \, |u^\eps_x|^2 \geq 0, 
\qquad
\\
&
C^\eps := u^\eps u^\eps_x + \delta \, \eps \, \big(u^\eps \, u^\eps_{xx} - (1/2) \, |u^\eps_x|^2 \big). 
\endaligned 
$$
The diffusive contribution decomposes into a non-positive term and a conservative one, while  the dispersive contribution is entirely conservative; hence, formally at least, as $\eps \to 0$ we recover the entropy inequality \eqref{eq:357}, {\sl but only} for the specific entropy function $U(u)= u^2/2$. In other words, we obtain the {\sl (single) quadratic entropy inequality} 
\bel{eq:250}
\bigl( u^2/2 \bigr)_t + F(u)_x \leq 0, \quad \qquad F' := u \, f'. 
\ee
In general, no specific sign is available for arbitrary convex entropies (unlike what we observed in the diffusion-only regime). 

%------------------------------------------------------------------------------------------------------

\subsection{Thin liquid film and  Camassa-Holm models} 
\label{sec:22}

More generally, we can consider the {\sl nonlinear diffusion-dispersion model} 
\bel{eq:129} 
u^\eps_t + f(u^\eps)_x
= \eps \, \bigl(b(u^\eps) \, u^\eps_x \bigr)_x + \delta \, \eps^2 \, \big(
c_1(u^\eps) \, \big( c_2(u^\eps) \, u^\eps_x \big)_x \big)_x,
\ee
where $b, c_1, c_2: \RR \to \RR$ are given smooth and positive functions. 
For this model, the formal limit $u=\lim_{\eps \to 0} u^\eps$ 
satisfies the following {\sl entropy inequality} determined by the function $c_1/c_2$   
\be
\label{eq:Ucc} 
\aligned
& U(u)_t + F(u)_x \leq 0, 
\\
& U'' = {c_2 \over c_1} >0, \qquad F' := f' \, U'. 
\endaligned
\ee
Namely, in the entropy variable $\widehat u = U'(u)$, the dispersive term takes the form 
$
\big(c_1(u) \, \big( c_2(u) \, u_x \big)_x\big)_x
= \big( c_1(u) \, \big( c_1(u) \, \widehat u_x \big)_x \big)_x,  
$
so that any solution to \eqref{eq:129} satisfies  
\be
\aligned
& U(u^\eps)_t  + F(u^\eps)_x = - D^\eps + \eps \, C^\eps_x, 
\quad 
\qquad D^\eps := \eps \, b(u) \, U''(u) \, |u_x|^2, 
\\
& C^\eps :=  b(u) \, U'(u) \, u_x 
           + \delta\eps \, \bigl(c_1(u) \widehat u \, \big( c_1(u) \, \widehat u_x \big)_x  - | c_2(u) \, u_x|^2 /2 \bigr). 
\endaligned 
\ee 

Nonlinear augmented terms also arise in the {\sl model of thin liquid films} 
\be
\label{thin}
u^\eps_t + (u_\eps^2 - u_\eps^3)_x 
= 
\eps \, (u_\eps^3 u_x^\eps)_x - \delta \eps^3 \, ( u_\eps^3 \, u^\eps_{xxx} )_x
\ee
with $\delta>0$ fixed and $\eps \to 0$. The scaling $\delta \eps^3$ is natural since the augmented term 
$( u_\eps^3 \, u^\eps_{xxx} )_x$
now involves four derivatives. 
The right-hand side describes the effects of surface tension in a thin liquid film moving on a surface and $u=u(t,x) \in [0,1]$ denotes the normalized thickness of the thin film layer. The parameters governing the forces and the slope of the surface are represented by the small parameter 
$\eps$. The model \eqref{thin} can be derived from the  so-called lubrication approximation of the Navier-Stokes equation when two counter-acting forces are taken into account: the gravity is responsible for pulling the film down an inclined plane while a thermal gradient (i.e.~the surface tension gradient) pushes the film up the plane. 
This model was studied by Bertozzi, M\" unch, and Shearer (2000), 
Bertozzi and Shearer (2000), as well as by Levy and Shearer (2004, 2005), 
LeFloch and Shearer (2004), 
Otto and Westdickenberg (2005), 
and LeFloch and Mohamadian (2008). 

It was observed by LeFloch and Shearer (2004) that the model \eqref{thin} satisfies the identity 
$$
\aligned
& ( u^\eps \log u^\eps - u^\eps )_t + \big( (u_\eps^2 - u_\eps^3) \log u^\eps - u^\eps + u_\eps^2 \big)_x 
= - D^\eps + \eps \, C^\eps_x,  
\\
&D^\eps := \eps \, u_\eps^3 \, |u^\eps_x|^2 + \gamma(\eps) \, | (u_\eps^2 \, u^\eps_x)_x|^2 \geq 0,
\endaligned
$$
so that, in the limit $\eps \to 0$, the following {\sl log-type entropy inequality} associated with the thin liquid film model holds: 
\be
\label{eq:Ulog}
( u \log u - u )_t  + \big( (u^2 - u^3) \log u - u + u^2 \big)_x \leq 0.
\ee

Consider finally the {\sl generalized  Camassa-Holm model} 
\be
\label{Cama}
u^\eps_t + f(u^\eps)_x = \eps \, u_{xx} + \delta \, \eps^2 \,
 ( u^\eps_{txx} + 2 \, u^\eps_x \, u^\eps_{xx} + u^\eps \, u^\eps_{xxx}), 
\ee 
which arises as a simplified shallow water model when wave breaking takes place. This equation was first investigated by Bressan and Constantin (2007) and Coclite and Karlsen (2006).  It was observed by LeFloch and Mohamadian (2008) that  \eqref{Cama} implies 
$$
\aligned
\bigl( (|u^\eps|^2 + \delta \eps^2 \, |u_x^\eps|^2)/2 \bigr){_t} + {F(u^\eps)_x} 
= - \eps \, |u^\eps_x|^2 + \eps \, C^\eps_x,   
\endaligned
$$
so that the formal limits $u=\lim_{\eps \to 0} u^\eps$ must satisfy the quadratic entropy inequality  \eqref{eq:250}, which was already derived for the linear diffusion-dispersion model \eqref{eq:240}. Although limiting solutions to \eqref{eq:240} and \eqref{Cama} look very similar in numerical tests, careful investigation (LeFloch and Mohamadian 2008) lead to the conclusion that they do not coincide. Consequently, we emphasize that two augmented models obtained by adding vanishing diffusive-dispersive terms to the same hyperbolic conservation laws do not generate the same nonclassical entropy solutions. In particular, if a numerical scheme is consistent with the quadratic entropy inequality, it need not converge to physically-relevant solutions. 
 
%------------------------------------------------------------------------------------------------------- 
 
\subsection{Fluids and elastic materials with phase transitions}
\label{sec:23}

Models for the dynamics of fluids and elastic materials (in one-space variable) are completely similar and, for  definiteness, we present the latter. 
The evolution of elastic materials undergoing phase transition may be described by the {\sl nonlinear elasticity model} 
\bel{5550}
\aligned 
 w_t - v_x &  = 0,
\\
v_t - \sigma(w)_x 
&  = 0.  
\endaligned    
\ee
Here, $w>-1$ denotes the deformation and $v$ the velocity of the material and, for typical materials, the stress-deformation relation $\sigma=\sigma(w)$ satisfies the monotonicity property 
\bel{hypp}
\sigma'(w) >0 \quad \text{ for all } w>-1. 
\ee
Under this condition, \eqref{5550} is strictly hyperbolic and admits the two wave speeds $-\lam_1 = \lam_2 = c(w)$ (the sound speed). 
The two characteristic fields are genuinely nonlinear in the sense of Lax (1974) if and only if $\sigma''$ never vanishes, which, however, fails for most materials encountered in applications as convexity is lost at $w=0$. We thus assume 
\be
\sigma''(w) \gtrless 0 \quad \text{ if } \quad  w \gtrless 0. 
\ee
Furthermore, following Slemrod (1983, 1989), the augmented version of \eqref{5550} reads 
\bel{555}
\aligned 
 w_t - v_x &  = 0,
\\
v_t - \sigma(w)_x 
&  = \eps \, v_{xx} - \delta \, \eps^2 \, w_{xxx},  
\endaligned    
\ee
and is refered to as the {\sl model of viscous-capillary materials} 
where the parameters $\eps$ and $\delta \, \eps^2$ are (rescaled) viscosity and capillarity coefficients.

Material undergoing phase transitions may be described by the model \eqref{555} but with 
a {\sl non-monotone} stress-strain function, satisfying  
\be
\label{888} 
\aligned 
& \sigma'(w) >0, \quad w \in (-1, w^m) \cup (w^M, +\infty), 
\\
& \sigma'(w) <0, \quad w \in (w^m, w^M) 
\endaligned    
\ee
for some constants $w^m < w^M$.
In the so-called {\sl unstable phase} $(w^m, w^M)$, the model admits two complex (conjugate) eigenvalues and is thus elliptic in nature. The system is hyperbolic in the non-connected set 
$\Ucal := \bigl(\RR \times (-1, w^m)\bigr) \cup 
      \bigl(\RR \times (w^M, +\infty)\bigr)$
and all 
solutions of interest for the hyperbolic lie {\sl outside} the unstable region. 

Recall that Slemrod (1983, 1989) first studied self-similar solutions to the Riemann problem 
(i.e.~the initial value problem with piecewise constant data), while Shearer (1986) introduced an explicit construction of the Riemann solutions when $\delta = 0$. 
It is only later that the notion of a kinetic relation for {\sl subsonic phase boundaries} 
was introduced by Truskinovsky (1987, 1993, 1994) 
and Abeyaratne and Knowles (1991a, 1991b). The latter solved the Riemann problem for \eqref{555} and investigated the existence of traveling wave solutions when $\sigma$ is a piecewise linear function. 
Next, LeFloch (1993) introduced a mathematical formulation of the kinetic relation for \eqref{555} and studied the initial value problem 
within a class of nonclassical entropy solutions with bounded variation and established an existence theory based on  Glimm's random choice scheme; therein, the kinetic relation was given an interpretation as the {\sl entropy dissipation measure.}  Further studies on this problem 
were then later by Corli and Sable-Tougeron (1997ab, 2000).

%------------------------------------------------------------------------------------ 

\subsection{Entropy inequality for models of elastodynamics}

 To the models in Section~\ref{sec:23}, we can associate the entropy 
\bel{eee}
\aligned
U(v,w) =& \frac {v^2} 2 + \Sigma(w), 
\qquad F(v,w) = - \sigma(w) \, v, 
\\
\Sigma(w):=& \int_0^w \sigma(s) \, ds,
\endaligned
\ee
which is {\sl strictly convex} under the assumption \eqref{hypp}.  
Namely, for the augmented model \eqref{555}, one has   
$$
\aligned 
& \left(\frac{v^2}{2} +  \Sigma(w) + {\delta \, \eps^2 \over 2} \, w_x^2\right)_t 
-  \bigl(v \, \sigma(w) \bigr)_x
\\
&= 
\eps \, \bigl(v \, v_x \bigr)_x  - \eps \, v_x^2 + \delta \, \eps^2 \, 
\bigl(v_x \, w_x - v \, w_{xx}\bigr)_x,
\endaligned 
$$ 
so that in the limit one formally obtains 
 the following {\rm entropy inequality associated with the phase transition model} 
\be
\label{eq:299}
\Big(\frac{v^2}{2} + \Sigma(w) \Big)_t - \bigl(v \, \sigma(w)\bigr)_x \leq 0.  
\ee 
One important difference between the hyperbolic and the hyperbolic-elliptic cases concerns the entropy (or total mechanical energy) \eqref{eee}, which
is convex in each hyperbolic region, but can not be extended to be globally convex in (the convex closure of) $\Ucal$. This model and its augmented version including viscosity and capillarity terms describe the dynamics of complex fluids with hysteresis and is relevant in many applications to phase transition dynamics: solid-solid interfaces, fluid-gas mixtures, etc.

More generally, we can assume an internal energy function $e = e(w, w_x)$ and then derive the evolution equations from the action
$$
J[v,w]  
: = \int_0^T \int_\Omega \Big(e(w, w_x) 
               - {v^2 \over 2} \Big) \, dx dt. 
$$
Namely, by defining the {\sl total stress} as  
$$
\Sigma(w,w_x, w_{xx}) 
:= { \del e \over \del w }(w,w_x) - \Big( { \del e \over \del w_x }(w,w_x) \Big)_x,   
$$ 
from the least action principle we deduce that a critical point of $J[v,w]$ satisfies 
$$ 
\aligned 
& v_t - \Sigma(w,w_x, w_{xx})_x = 0, \\ 
& w_t  - v_x = 0. 
\endaligned 
$$ 
If we also include the effect of a (nonlinear) viscosity $\mu=\mu(w)$, we arrive at a fully nonlinear phase transition model with viscosity and capillarity:
$$ 
\aligned 
& w_t - v_x = 0,  
\\
& v_t - \Sigma(w,w_x, w_{xx})_x = \bigl(\mu(w) \, v_x \bigr)_x.
\endaligned  
$$  
Again, the total energy 
$E(w,v,w_x) := e(w,w_x) + v^2 /2$ 
plays the role of a mathematical entropy, and we find
$$
\aligned
&
E(w,v,w_x)_t - \bigl(\Sigma(w,w_x, w_{xx}) v\bigr)_x 
\\
&=  \Big( v_x { \del e \over \del w_x }(w,w_x) \Big)_x 
  + \bigl( \mu(w) v v_x \bigr)_x - \mu(w) v_x^2 
\endaligned
$$  
and once more, a single entropy inequality is obtained.  
 
In the case that $e$ is quadratic in $w_x$, for some positive capillarity coefficient $\lam(w)$, we set 
$e(w,w_x) = \eps(w) + \lam(w) \, {w_x^2 \over 2}$, and  
the total stress decomposes as  
$$ 
\Sigma(w,w_x, w_{xx}) =  \sigma(w) + \lam'(w) \, {w_x^2 \over 2} - (\lam(w) \, w_x )_x, 
\quad 
\sigma(w) = \eps' (w). 
$$ 
The evolution equations then take the form  
\be
\label{666}
\aligned 
& w_t  - v_x = 0,
\\
& v_t - \sigma(w)_x  
 = \Big( \lam'(w) \, {w_x^2 \over 2} - \bigl(\lam(w) \, w_x \bigr)_x \Big)_x 
   +  \bigl(\mu(w) \, v_x \bigr)_x. 
\endaligned 
\ee
We then obtain 
$$
\aligned 
& \Big(\eps(w) + {v^2 \over 2} + \lam(w) \, {w_x^2 \over 2}\Big)_t 
  - \bigl(\sigma(w) \, v \bigr)_x 
\\
& =        \bigl( \mu(w) \, v \, v_x \bigr)_x  - \mu(w) \, v_x^2
		 + \Big(v \, {\lam'(w)\over 2} \, w_x^2  
                - v \, \bigl(\lam(w) \, w_x \bigr)_x 
                + v_x \, \lam(w)\, w_x \Big)_x,
\endaligned 
$$  
which, again, leads to the entropy inequality \eqref{eq:299}.  When the viscosity and capillarity are taken to be constants, we recover the previous model and, again, the entropy inequality is identical for both regularizations.   

%------------------------------------------------------------------------------------------------------------------

\subsection{Kinetic relations for nonclassical shocks}
\label{sec:25}

All the models in this section admit shock wave solutions that do satisfy a single entropy inequality in the sense of distributions (that is, \eqref{eq:250}, \eqref{eq:Ucc}, \eqref{eq:Ulog}, or \eqref{eq:299}), but fail to satisfy standard entropy conditions (Oleinik, Kruzkov, Lax, Wendroff, etc.). However, these entropy conditions played an essential role in the design of efficient shock-capturing schemes for standard fluid dynamics problems (LeVeque 2003). 

The new shocks are referred to as as nonclassical shocks and can be checked to be undercompressive, in the sense that they are linearly unstable, since a perturbation passes through them rather than impinging on them (as is the case of compressive shocks). For this reason, an additional admissibility condition is necessary, which is called a kinetic relation. It takes the form of an additional jump relation at shock discontinuities. 

Kinetic relations can be obtained analytically only for the simpler models (scalar conservation laws with linear diffusion and dispersion, nonlinear elasticity model with linear viscosity and capillarity), so numerical approaches are necessary to tackle these problems. In the references already cited, numerical investigations have established that kinetic functions exist (and often satisfy certain monotonicity properties)
 for a large class of physically relevant models including thin liquid films,  generalized Camassa-Holm, nonlinear phase transitions, van der Waals fluids (for small shocks), and magnetohydrodynamics.

%-----------------------------------------------------------------------------------------------------------

\subsection{Other physical models and applications}
\label{sec:71}

The methods and techniques to be presented in this paper are also relevant for other classes of problems. For instance, 
the Buckley-Leverett equation for two-phase flows in porous media provides another model of interest for the applications, and was studied 
Hayes and Shearer (1999) and Van~Duijn, Peletier, and Pop (2007). There are also other models of great 
physical interest which, however, have not yet received as much attention. Since the hyperbolic flux of these models
does admit an inflection point and that dispersive-type effects are important in the modeling of such problems,  
it is expected that undercompressive shocks shoud occur, at least in certain regimes of applications.  This is the case 
of the quantum hydrodynamics models (Marcati, Jerome), phase field models (Caginalp, Ratz), 
Suliciu-type models (Carboux, Frid, Suliciu), non-local models involving fractional integrals (Kissiling et al., Rohde)  
discrete molecular models based on potentials (for instance of Lennard-Jones type) (Bohm, Dreyer, NT, TV, Weinan). 
Cf.~the references at the end of this paper. 

Another rich direction of research where analogue challenges are met 
is provided by the coupling techniques involving two hyperbolic systems with distinct flux-functions, for which we refer to 
LeFloch (1993), 
Seguin and Vovelle (2003), 
Godlevski and Raviart (2004), 
Adimurthi, Mishra, and Gowda (2005), 
Godlevski, Le Thanh, and Raviart (2005), 
Bachman and Vovelle (2006), 
B\"urger and Karlsen (2008), 
Chalons, Raviart, and Seguin (2008),  
Holden, Karlsen, Mitrovic, and Panov (2009),
Boutin, Coquel, and LeFloch (2011, 2012, 2013),
and the references cited therein. 

%==================================================================================

\section{Schemes with controlled dissipation for nonclassical entropy solutions} 
\label{sec:3}

\subsection{Standard finite difference or finite volume schemes}
\label{sec:31}

We consider a nonlinear hyperbolic system in the conservative form \eqref{eq:109} and we now discretize it (in space) on a grid consisting of points $x_i = i \Dx$, with $\Dx$ being a uniform mesh width. (The grid is assumed to be uniform for the sake of simplicity in the exposition.)
 On this grid, a standard (semi-discrete) finite difference scheme provides an approximation of the point values $u_i(t) \approx u(t,x_i)$ of  solutions to \eqref{eq:109}, defined by
\be
\label{eq:fdm1}
\frac{d}{dt} u_i + \frac{1}{\Dx}\left(g_{i+1/2}(t) - g_{i-1/2}(t) \right) = 0.
\ee
Here, $g_{i+1/2} := g(u_i,u_{i+1})$ is a consistent {\sl numerical flux} associated with of the flux $f$, that is, $g(a,a) = f(a)$ for all relevant $a$. Alternatively, by considering the cell averages
$$
u_i(t) \approx \frac{1}{\Dx} \int_{x_{i-1/2}}^{x_{i+1/2}} u(t,x) dx, 
$$
one formulates a finite volume scheme which has the same form \eqref{eq:fdm1}.

A well-known theorem due to Lax and Wendroff (1960) establishes that if the numerical approximation converges (in a suitable sense), it can converge only towards a weak solution of the underlying system \eqref{eq:109}. Furthermore, if some structural conditions on the numerical flux $g$ are assumed as in (Tadmor 1987, 2003), one can show that the scheme \eqref{eq:fdm1} in the limit also satisfies a discrete version of the entropy inequality \eqref{eq:357}, that is, 
\be
\label{eq:efdm1}
\frac{d}{dt} U(u_i) + \frac{1}{\Dx}\left(G_{i+1/2}(t) - G_{i-1/2}(t) \right) \leq 0.
\ee
Here, the \emph{numerical entropy flux} $G_{i+1/2} := G(u_i,u_{i+1})$ is consistent with the entropy flux $F$ in \eqref{eq:357}, in the sense that $G(a,a) = F(a)$ for all relevant $a$. When such a discrete entropy inequality is available, one can readily modify Lax and Wendroff's argument and show that if the approximations generated by the finite difference (or finite volume) scheme converge, then they can only converge toward an entropy solution of \eqref{eq:109}.

A variety of numerical fluxes that satisfy the entropy stability criteria (as stated in Tadmor 1987, 2003) have been designed in the last three decades. Theses classes of schemes include exact Riemann solvers of Godunov-type, approximate Riemann solvers such as Roe and Harten-Lax-van~Leer (HLL) solvers, as well a central difference schemes such as Lax-Friedrichs and Rusanov schemes. A comprehensive description of these schemes and their properties are available in the literature, for instance in the textbook by LeVeque (2003). 

\subsubsection{Failure of standard schemes to approximate nonclassical shocks}

\begin{figure}[htbp]
\centering
\includegraphics[height=4.5cm, width=8cm]{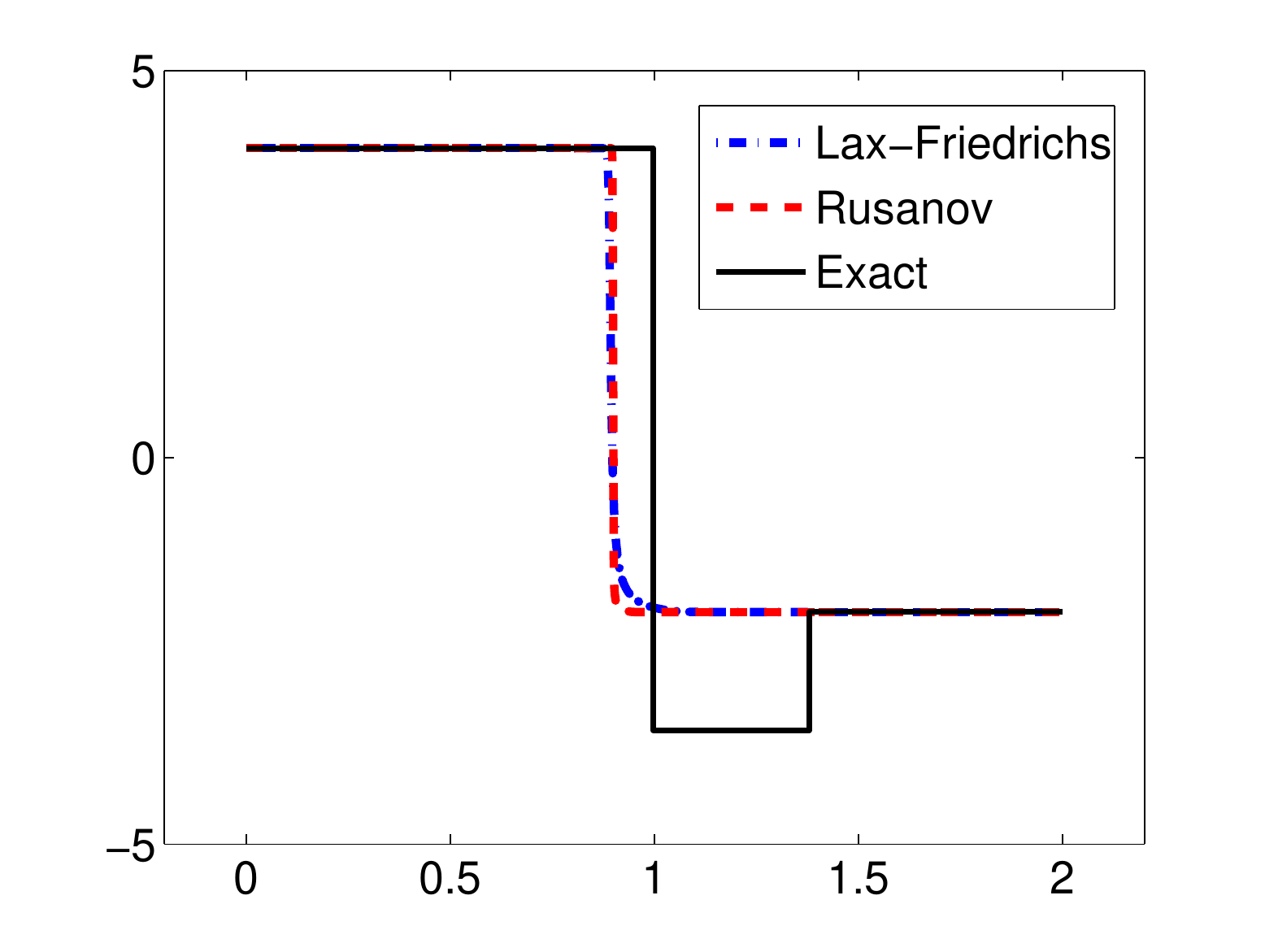}
\caption{Approximation of small-scale dependent shock waves for the cubic conservation law with vanishing diffusion and capillarity \eqref{eq:240} using the standard Lax-Friedrichs and Rusanov schemes}
\label{fig:1}
\end{figure}

As mentioned in the introduction, standard conservative and entropy stable schemes \eqref{eq:fdm1}--\eqref{eq:efdm1}
 fail to approximate nonclassical shocks (and other small-scale dependent solutions). As an illustrative example, we consider here the cubic scalar conservation law with linear diffusion and dispersion, that is, \eqref{eq:240} with $f(u) = u^3$ and a fixed $\delta >0$. The underlying conservation law is approximated with the standard Lax-Friedrichs and Rusanov schemes and the resulting solutions are plotted in Figure~\ref{fig:1}. The figure clearly demonstrates that Godunov and Lax-Friedrichs schemes, both, converge to the \emph{classical} entropy solution to the scalar conservation law and, therefore, do not approximate the nonclassical entropy solution, realized as the vanishing diffusion-dispersion limit of \eqref{eq:240} and also plotted in the same figure. The latter consists of three distinct constant states separated by two shocks, while the classical solution contains a single shock. 
\par As pointed out in the introduction, this failure of standard schemes in approximating small-scale dependent shocks (in various contexts) can be explained in terms of the equivalent equation of the scheme (as was first observed by Hou and LeFloch (1994) and Hayes and LeFloch (1996, 1998). The equivalent equation is derived via a (formal) Taylor expansion of the discrete scheme \eqref{eq:fdm1} and 
contains mesh-dependent terms and high-order derivatives of the solution. 

For a first-order scheme like \eqref{eq:fdm1}, the equivalent equation has the typical form   
\be
\label{eq:eqen}
\quad 
u_t + f(u)_x = \Dx (\overline{b}(u)u_x)_x + \delta \Dx^2 (\overline{c}_1(u)(\overline{c}_2(u)u_x)_x)_x + \Ord(\Dx^3).
\ee
Interestingly enough, this equation is of the augmented form \eqref{eq:129} with $\epsilon = \Dx$ being now the small-scale parameter. Standard schemes often have $\overline{b} \neq b$, $\overline{c}_1 \neq c_1$, and $\overline{c_2} \neq c_2$,
with $b, c_1, c_2$ being small-scale terms prescribed in the nonlinear physical model like \eqref{eq:129}. As the shocks realized as the $\epsilon \to 0$ limit of (\ref{eq:129}) depend explicitly on the expressions of the diffusion and dispersion terms, this difference in the diffusion and dispersion terms between (\ref{eq:eqen}) and (\ref{eq:129}) is the \emph{crucial} reason as to why standard schemes fail to correctly approximate small-scale dependent shocks.

%----------------------------------------------------------------------------------------------------------

\subsection{Finite difference schemes with controlled dissipation}
\label{sec:32}

The equivalent equation of the finite difference scheme \eqref{eq:eqen} also suggests a way for modifying the scheme such that the correct small-scale dependent shock waves can be approximated. Following Hayes and LeFloch (1997, 1998), LeFloch and Rohde (2000), and LeFloch and Mohamadian (2008), the key idea is to design finite difference schemes whose equivalent equation matches,  both, the diffusive and the dispersive terms in the augmented model 
\eqref{eq:129} (for instance). Namely, the schemes are designed so that their equivalent equation (\ref{eq:eqen}) has
$\overline{b} = b$, $\overline{c}_1 = c_1$, $\overline{c}_2 = c_2$. Such schemes are refered to as {\sl schemes with controlled dissipation}. 

In order to now proceed with the derivation of a class of schemes with controlled dissipation, we focus attention on the prototypical example of  the scalar conservation law, and we consider nonclassical shocks generated in the limit of balanced vanishing diffusion and dispersion; cf.~\eqref{eq:240}. On the uniform grid presented in the previous subsection and for any integer $p \geq 1$, we approximate the conservation law (\ref{eq:109}) with the following $2p$-th order consistent finite difference scheme:
\be
\label{eq:fds2p}
\qquad 
\frac{du_i}{dt} + \frac{1}{\Dx} \sum\limits_{j=-p}^{j=p} \alpha_j f_{i+j} 
 = \frac{c}{\Dx} \sum\limits_{j=-p}^{j=p} \beta_j u_{i+j}
 +  \frac{\delta c^2}{\Dx} \sum\limits_{j=-p}^{j=p} \gamma_j u_{i+j}.  
\ee
Here, $u_i(t) \approx u(x_i,t)$ is the cell nodal value, $f_{i} = f(u_i)$ is the flux, the constant $\delta$ is the coefficient of capillarity (given by the physics of the problem) and $c \geq 0$ being a positive constant. The coefficients $\alpha_j,\beta_j$ and $\gamma_j$ need to satisfy the following \emph{$2p$-order conditions}: 
\begin{eqnarray}
\label{eq:al}
\sum\limits_{j=-p}^{p} j \alpha_j  = 1, \qquad
\sum\limits_{j=-p}^{p} j^l \alpha_j = 0, \quad  \qquad l \neq 1,  \quad 0 \leq l \leq 2p. 
\end{eqnarray}  
These conditions define a set of $(2p+1)$ linear equations for $(2p+1)$ unknowns and can be solved explicitly. Similarly,
 the coefficients $\beta$ must satisfy
\begin{eqnarray}
\label{eq:ba}
\sum\limits_{j=-p}^{p} j^2 \beta_j  = 2, 
\qquad 
\sum\limits_{j=-p}^{p} j^l \beta_j = 0, \quad \qquad l \neq 2, \quad \quad 0 \leq l \leq 2p
\end{eqnarray}  
while, for the coefficients $\gamma$,  
\begin{eqnarray}
\label{eq:ga}
\sum\limits_{j=-p}^{p} j^3 \gamma_j  = 6, 
\qquad 
\sum\limits_{j=-p}^{p} j^l \gamma_j = 0, \quad \qquad l \neq 3, \quad \quad 0 \leq l \leq 2p. 
\end{eqnarray}  

The proposed finite difference scheme (\ref{eq:fds2p}) is a conservative and consistent discretization of the conservation law (\ref{eq:109}). It is formally only first-order accurate since the diffusive terms are proportional to $\Dx$. This scheme need not preserve the monotonicity of the solutions. 

The equivalent equation associated with the scheme (\ref{eq:fds2p}) reads 
\be
\label{eq:ee}
\aligned
\frac{du}{dt} + f(u)_x 
= \,  & c \, \Dx \, u_{xx} + \delta c^2 (\Dx)^2 \, u_{xxx} 
                                        - \sum\limits_{k=2p+1}^{\infty} \frac{(\Dx)^{k-1}}{k!} A^p_k (f(u))^{[k]} 
\\
                                        &+  c \sum\limits_{k=2p+1}^{\infty} \frac{(\Dx)^{k-1}}{k!} B^p_k u^{[k]} 
+  \delta c^2 \sum\limits_{k=2p+1}^{\infty} \frac{(\Dx)^{k-1}}{k!} C^p_k u^{[k]}. 
                                        \endaligned
                                        \ee
Here, 
$g^{[k]} = \frac{d^k g}{dx^k}$ denotes the $k$-th spatial derivative of a function $g$ and the above coefficients are 
\begin{eqnarray}
A^p_k = \sum\limits_{j=-p}^{p} \alpha_j j^k, 
\qquad 
B^p_k = \sum\limits_{j=-p}^{p} \beta_j j^k, 
\label{eq:coeff}
 \qquad
C^p_k = \sum\limits_{j=-p}^{p} \gamma_j j^k,
\end{eqnarray}
with $\alpha,\beta$ and $\gamma$ being specified by the relations (\ref{eq:al}), (\ref{eq:ba}), and (\ref{eq:ga}), respectively.

In view of the equivalent equation \eqref{eq:ee}, the numerical viscosity and dispersion terms are linear and match the underlying diffusive-dispersive equation \eqref{eq:240} with $\epsilon = c\Dx$. Hence, the numerical diffusion and dispersion are ``controlled'' in the sense that they match the underlying small-scale terms.  

\begin{figure}[htbp]
\centering
\includegraphics[height=4.5cm, width=8cm]{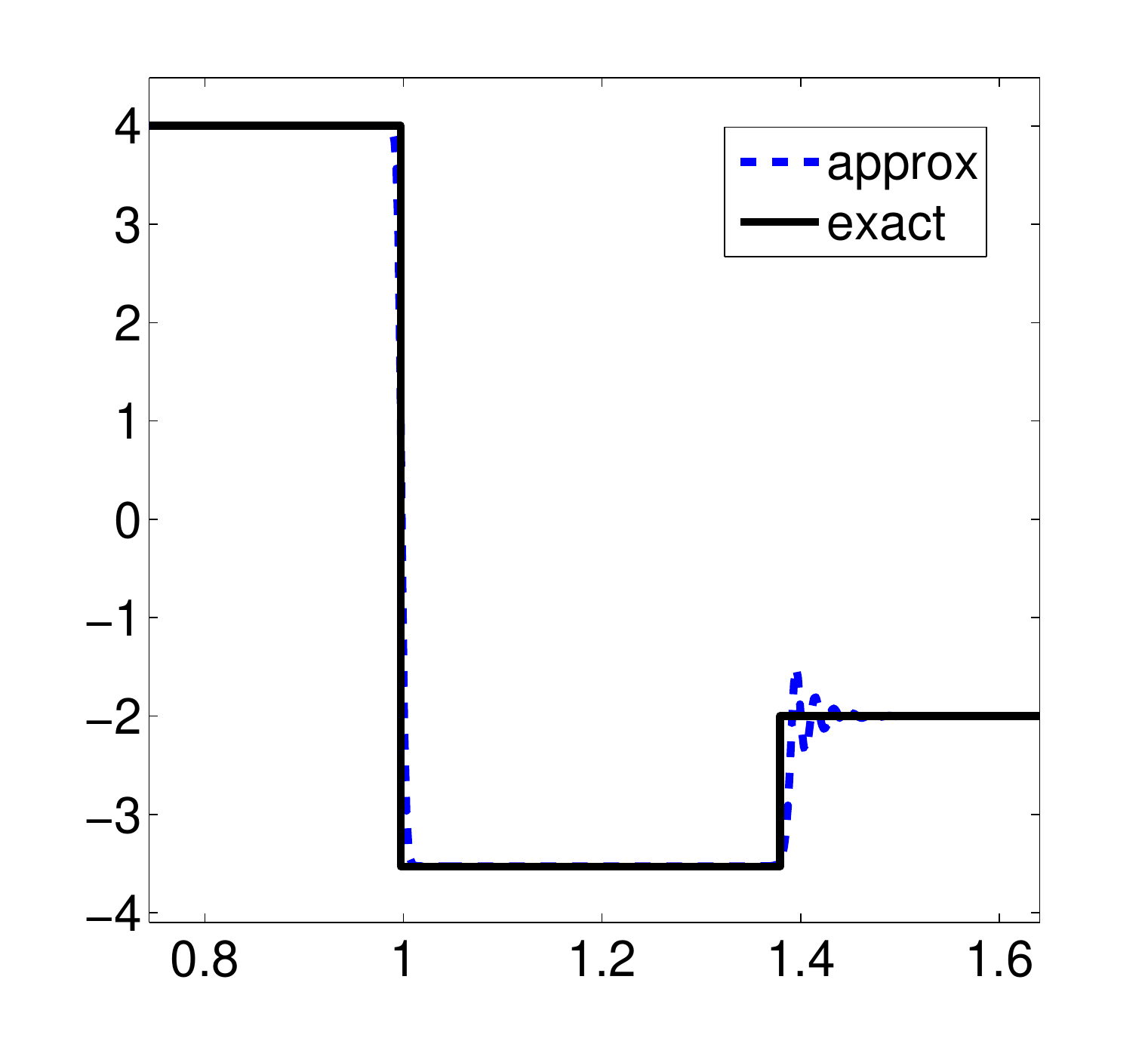}
\caption{Approximation of nonclassical shocks to the cubic conservation law (with vanishing diffusion and capillarity) \eqref{eq:240} using 
LeFloch-Mohamadian's sixth-order scheme with controlled dissipation \eqref{eq:fds2p} with $c = 5$}
\label{fig:2}
\end{figure}

The above schemes approximate small-scale dependent solutions very well. Indeed, let us illustrate their performance for a representative example of the cubic scalar conservation law with vanishing diffusion and dispersion \eqref{eq:240}. A sixth-order ($p=3$) finite difference scheme with controlled dissipation \eqref{eq:fds2p} was originally proposed by LeFloch and Mohamadian (2008) and we can take the coefficient $c=5$ and the same initial data, as in Figure~\ref{fig:1}. The results shown in Figure~\ref{fig:2} demonstrate that the proposed scheme with controlled dissipation converges the correct nonclassical shock wave. This is in sharp contrast with the failure observed with standard  schemes such as Rusanov and Lax-Friedrichs schemes which converge to classical solutions; see Figure~\ref{fig:1}. 

\subsubsection{The problem of strong shocks}

\begin{figure}[htbp]
\centering
\hskip-.75cm 
\begin{tabular}{lr}
 \includegraphics[height=4.5cm, width=0.45\linewidth]{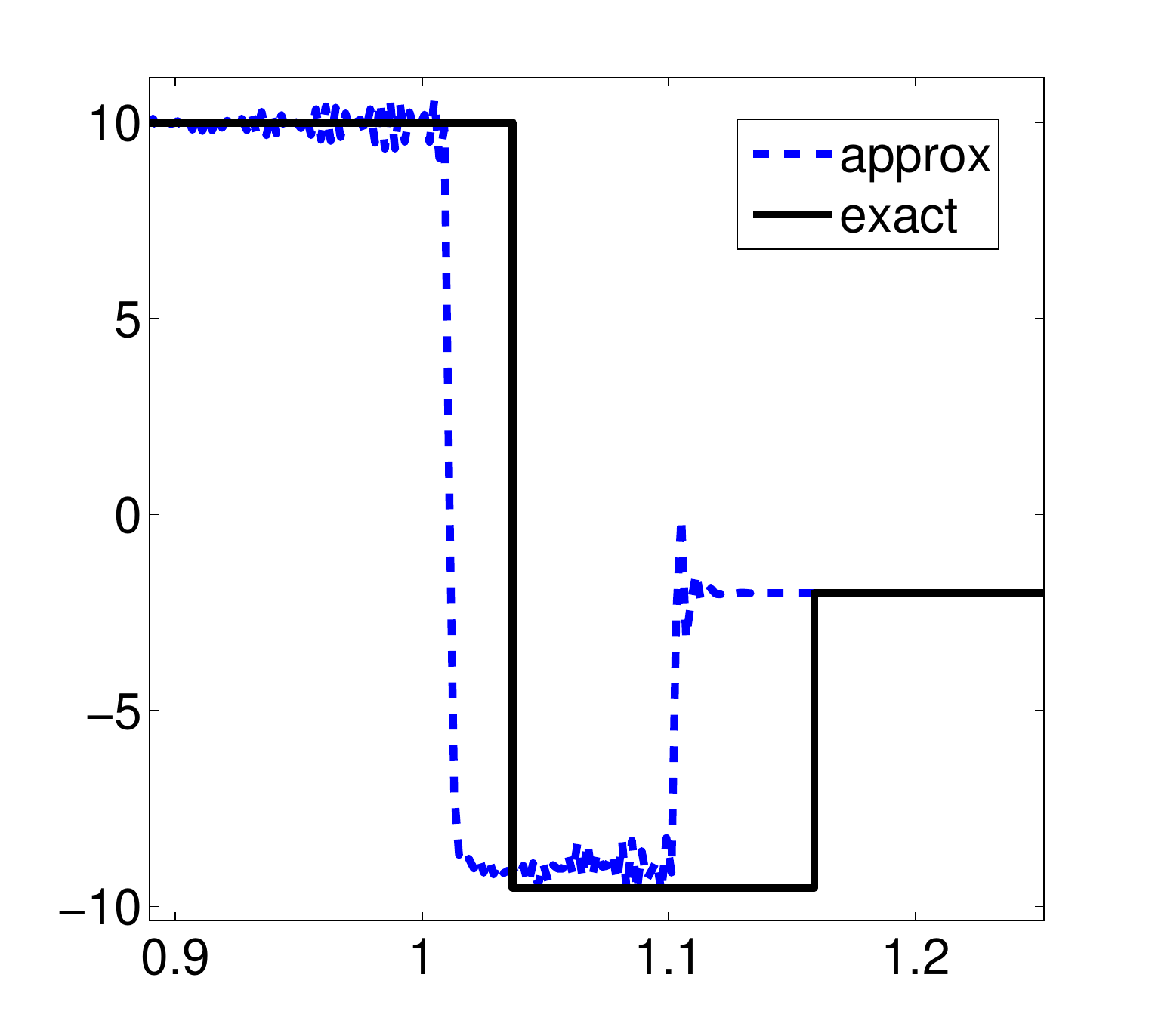} 
& \includegraphics[height=4.5cm, width=0.45\linewidth]{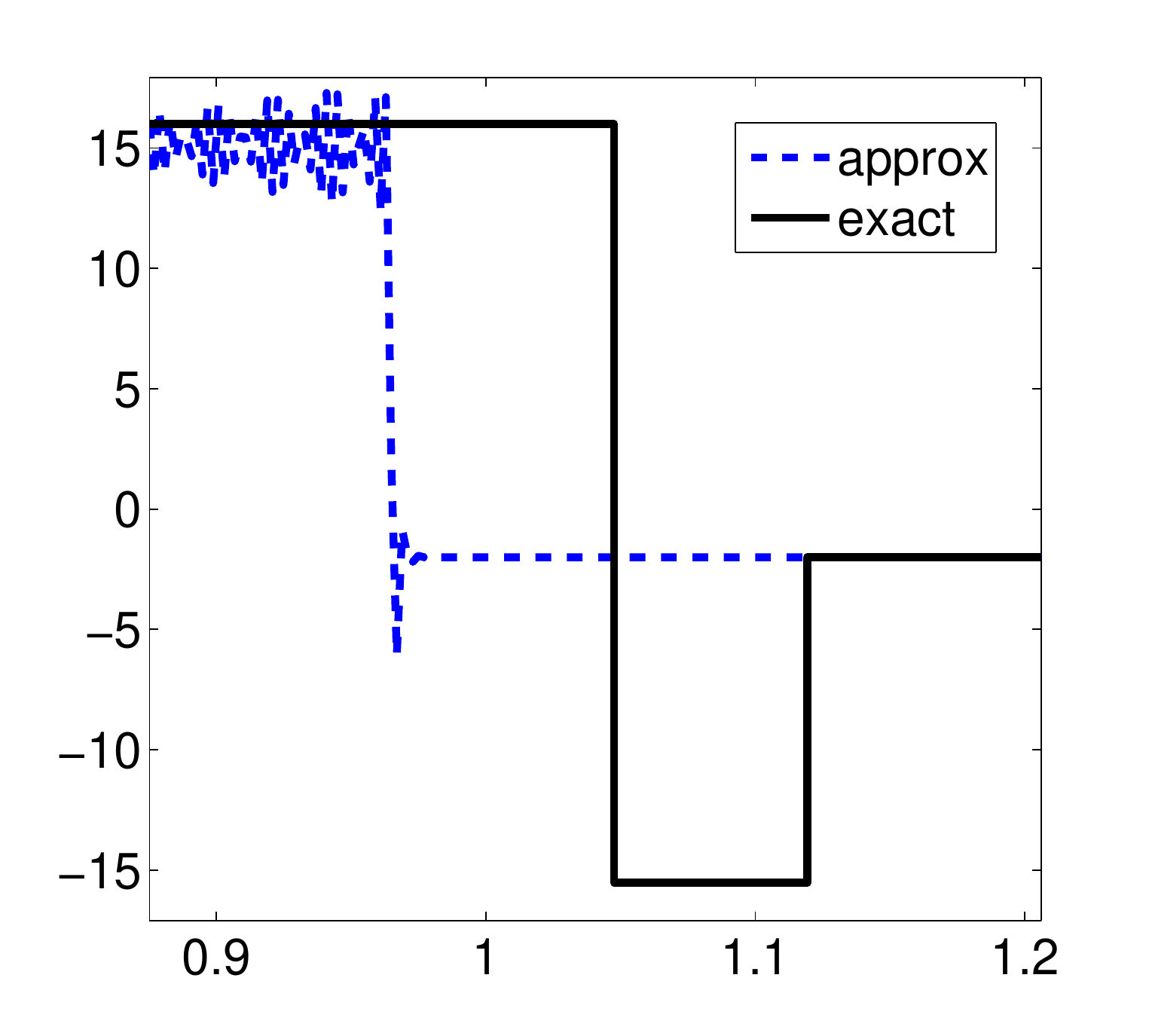} 
\end{tabular}
\caption{Approximation of strong nonclassical shocks  for the cubic conservation law \eqref{eq:240} using 
 LeFloch-Mohamadian's sixth-order scheme with controlled dissipation \eqref{eq:fds2p} with $c = 5$. Left: Initial jump $=10$. Right: Initial jump $=16$.} 
\label{fig:3}
\end{figure}

Schemes with controlled dissipation such as \eqref{eq:fds2p} approximate nonclassical solutions quite well in most circumstances. However, the  approximation significantly deteriorates for sufficiently strong shocks. As an example, we can consider strong nonclassical shocks for the cubic conservation law \eqref{eq:240} and attempt to approximate them with the sixth-order finite difference scheme with controlled dissipation \eqref{eq:fds2p}; the results are displayed in Figure~\ref{fig:3}. The figures clearly show that the sixth-order scheme with controlled dissipation fails at accurately resolving {\sl strong} nonclassical shocks. In particular, this scheme completely fails to approximate a large amplitude nonclassical shock with amplitude of around $30$.

%-----------------------------------------------------------------------------------------------

\subsection{WCD schemes for scalar conservation laws}
\label{sec:33-b}

In a recent work (Ernest, LeFloch, and Mishra 2013), the authors have identified the key reason for the failure of schemes with controlled dissipation at approximating nonclassical solutions with large amplitude, especially strong shocks. 
Again, the equivalent equation associated with the finite difference scheme \eqref{eq:fds2p} explains this behavior. As pointed out before, we have designed schemes with controlled dissipation in such a manner that the numerical diffusion and dispersion terms match the underlying diffusion and dispersion in the model of interest \eqref{eq:240}. However, as seen in the equivalent equation \eqref{eq:ee}, the higher-order terms (i.e.~of the order $\Ord(\Dx^3)$ and higher) do play a role, particularly when the approximated shock is strong. To illustrate this fact, we consider the equivalent equation \eqref{eq:ee} for a single shock wave.

Namely, consider a single shock connecting two states $u_L, u_R$ such that $[[u]] = u_L - u_R > 0$ and $[[f(u)]] > 0$ (the other cases being handled similarly). At this shock discontinuity, we formally find 
\begin{eqnarray*}
u^{[k]} \approx \frac{[[u]]}{\Dx^k}, \qquad
\qquad 
(f(u))^{[k]} \approx \frac{[[f(u)]]}{\Dx^k}.
\end{eqnarray*}
Substituting these formal relations into the equivalent equation (\ref{eq:ee}) at a single shock, we obtain
\be
\label{eq:s1}
\frac{du}{dt}+  \underbrace{\frac{[[f(u)]]}{\Dx} - \frac{c[[u]]}{\Dx} - \frac{\delta c^2 [[u]]}{\Dx}}_{l.o.t} \approx
 \underbrace{\frac{S_p^{D} c [[u]]}{\Dx} + \frac{S_p^{C} \delta c^2 [[u]]}{\Dx} - \frac{S_p^f [[f]]}{\Dx}}_{h.o.t.},
\ee
in which the coefficients read  
\be
\label{eq:s1-two}
S^f_p = \sum_{k=2p+1}^{\infty} \frac{A^p_k}{k!}, 
\qquad
S^D_p = \sum_{k=2p+1}^{\infty} \frac{B^p_k}{k!}, 
\qquad
S^C_p = \sum_{k=2p+1}^{\infty} \frac{C^p_k}{k!}, 
\ee
with $A^p_k,B^p_k,C^p_k$ being defined in (\ref{eq:coeff}).

The relation (\ref{eq:s1}) represents the balance of terms in the equivalent equation in the neighborhood of a single shock. Ideally, the higher-order error terms (h.o.t. in (\ref{eq:s1})) should be dominated in amplitude by the leading-order terms (l.o.t. in (\ref{eq:s1})).

\subsubsection{The WCD condition}

Ernest, LeFloch, and Mishra (2013) have sought to balance both sets of terms through a user-defined tolerance parameter $\tau << 1$. In order words, the condition 
$
\frac{|h.o.t|}{|l.o.t|} < \tau,
$
is imposed by a comparing, on one hand, the upper bound
$$
|h.o.t| \leq \left(\widehat{S}^D_p c + \widehat{S}^C_p |\delta| c^2 + \widehat{S}^f_p \sigma\right) \frac{|[[u]]|}{\Dx},
$$
$\sigma = \big|{[[f(u)]] \over  [[u]]} \big|$ being the shock speed and 
\be
\label{eq:s5}
\aligned
& \widehat{S}^f_p = \hskip-.15cm \sum_{k=2p+1}^{\infty} \left | \sum\limits_{j=-p}^{p} \frac{\alpha_j j^k}{k!}\right|,
\qquad 
\widehat{S}^D_p =  \hskip-.15cm \sum_{k=2p+1}^{\infty} \left | \sum\limits_{j=-p}^{p} \frac{\beta_j j^k}{k!} \right|,
\\
& \widehat{S}^C_p =  \hskip-.15cm \sum_{k=2p+1}^{\infty} \left |\sum\limits_{j=-p}^{p} \frac{\gamma_j j^k}{k!}\right|, 
\endaligned
\ee
with, on the other hand, the lower bound 
$$
|l.o.t| \geq \left(|\delta|c^2 + c - \sigma\right)\frac{|[[u]]|}{\Dx}.
$$
Therefore, we can achieve the condition $\frac{|h.o.t|}{|l.o.t|} < \tau$ {\sl provided} 
\be
\label{eq:wcd}
({\bf WCD}): \quad \left( |\delta|- \frac{\widehat{S}^C_p |\delta|}{\tau}\right) c^2 + \left(1 - \frac{\widehat{S}^D_p}{\tau}\right)c - \left(1 + \frac{\widehat{S}^f_p}{\tau}\right)\sigma > 0, 
\ee
which we refer to as the {\sl WCD condition} associated with the proposed class of schemes.  

Recall that, in \eqref{eq:wcd}, $\tau$ is a user-defined tolerance, $\delta$ is the coefficient of dispersion, $\widehat{S}_P^{f,D,C}$ are specified in (\ref{eq:s5}) and can be computed in advance (before the actual numerical simulation), while $\sigma$ is the shock speed (of the shock connecting $u_L$ and $u_R$) and depends on the solution under consideration. The only genuine parameter to be chosen is the numerical dissipation coefficient $c$. In contrast to schemes with controlled dissipation where $c$ was set to be a constant, it is now natural 
that this coefficient be {\sl time dependent} $c = c(t)$ and evaluated at each time step and chosen to satisfy the WCD condition (\ref{eq:wcd})

An important question is whether there exists a suitable $c$ such that the WCD condition (\ref{eq:wcd}) is satisfied for a given ordre $p$. In fact,  elementary properties of Vandermonde determinants (as observed by Dutta 2013) imply that the coefficients \eqref{eq:s5} satisfy
\be
\label{eq:lim}
\lim\limits_{p \to +\infty} \max \Big( \widehat{S}^f_p,\widehat{S}^D_p,\widehat{S}^C_p\Big) = 0.
\ee 
As $c$ is a coefficient of diffusion in (\ref{eq:fds2p}), we need that $c > 0$ and a simple calculation based on the quadratic relation (\ref{eq:wcd}) shows that $c > 0$ if and only if 
\be
\label{eq:con1}
\frac{\widehat{S}^C_p}{\tau} < 1. 
\ee
Recall that $\sigma$ (being the absolute value of the shock speed) is always positive. 

Given (\ref{eq:lim}), we can always find a sufficiently large $p$ such that the sufficient condition (\ref{eq:con1}) is satisfied for any given tolerance $\tau$. Hence, the ``order'' of the finite difference scheme (\ref{eq:fds2p}) needs to be increased in order to control the high-order terms in the equivalent equation, in terms of the leading order terms. The strategy of letting the exponent $p$ tend to infinity goes back to LeFloch and Mohamadian (2008), who established the convergence of the numerical kinetic function to the analytical kinetic function. 

Furthermore, in the limit of infinite $p$ and for any $\tau > 0$, the property \eqref{eq:lim} leads us to the following {\sl limiting version of the WCD condition:}
\be
\label{eq:con2}
|\delta| c^2 + c - \sigma > 0.
\ee
Solving the quadratic equation explicitly yields two real roots, one being negative and the other one positive. The convexity of the function implies that the choice $c > c_{2}$ (with $c_2$ being the positive root of the above quadratic equation) will satisfy the WCD condition. Thus for sufficiently large $p$ (that is, sufficiently high-order schemes), we can always choose a suitable numerical dissipation coefficient (depending on both $\delta$ and the wave speed $\sigma$) that yields the correct small-scale dependent solutions. 

Finally, we extend the above analysis at a single shock and we determine the diffusion coefficient $c$ in the finite difference scheme (\ref{eq:fds2p}) in the following manner. At each interface $x_{i+1/2} = \frac{1}{2}(x_i + x_{i+1})$, we use $u_L = u_{i}$ and $u_R = u_{i+1}$ in the WCD condition (\ref{eq:wcd}) and choose a coefficient $c_{i+1/2}$ such that this condition is satisfied. The coefficient for the entire scheme is then given by $c := c(t) = \max_{i} c_{i+1/2}$.

%-----------------------------------------------------------------------------------------------------------------

\subsection{Numerical experiments}

Following Ernest, LeFloch, and Mishra (2013), we test here the WCD schemes (\ref{eq:fds2p}) for the cubic scalar conservation law (\ref{eq:109}) (with $f(u) = u^3$) and, for definiteness, we set $\delta = 1$ in (\ref{eq:240}). The finite difference schemes defined above are semi-discrete, and we now also discretize the equation in time using a third-order, strong stability preserving, Runge-Kutta time stepping method. The time step is determined using a standard stability condition with CFL number $=.45$ for all numerical experiments.   

In order to compute the coefficient $c$, we need to choose $\tau$ as well as the order $2p$ of the scheme and, then, compute the dissipation coefficient suggested by the WCD condition. We use here the following Riemann initial data 
\begin{eqnarray*}
u(0,x) &=&  
               u_L \quad {\rm if} \quad x < 0.4; \qquad 
               \quad -2 \quad  {\rm if} \quad x > 0.4,
               \end{eqnarray*}
and we vary the state $u_L$ in order to cover various shock strengths. 

\subsubsection{Small shocks}

\begin{figure}[htbp]
\centering
\begin{tabular}{cc}
\includegraphics[height=4.5cm, width=0.48\linewidth]{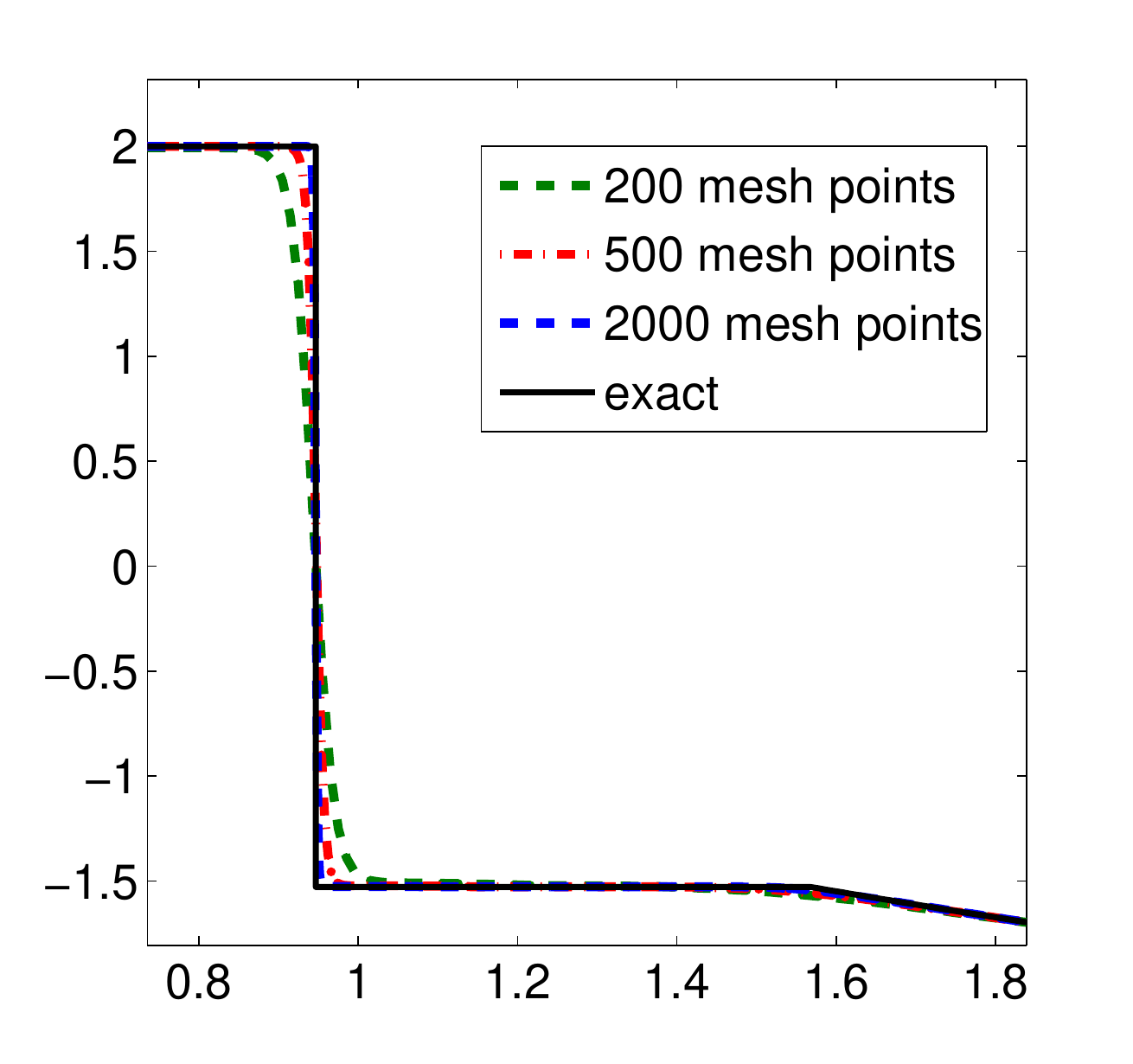} 
& \includegraphics[height=4.5cm, width=0.48\linewidth]{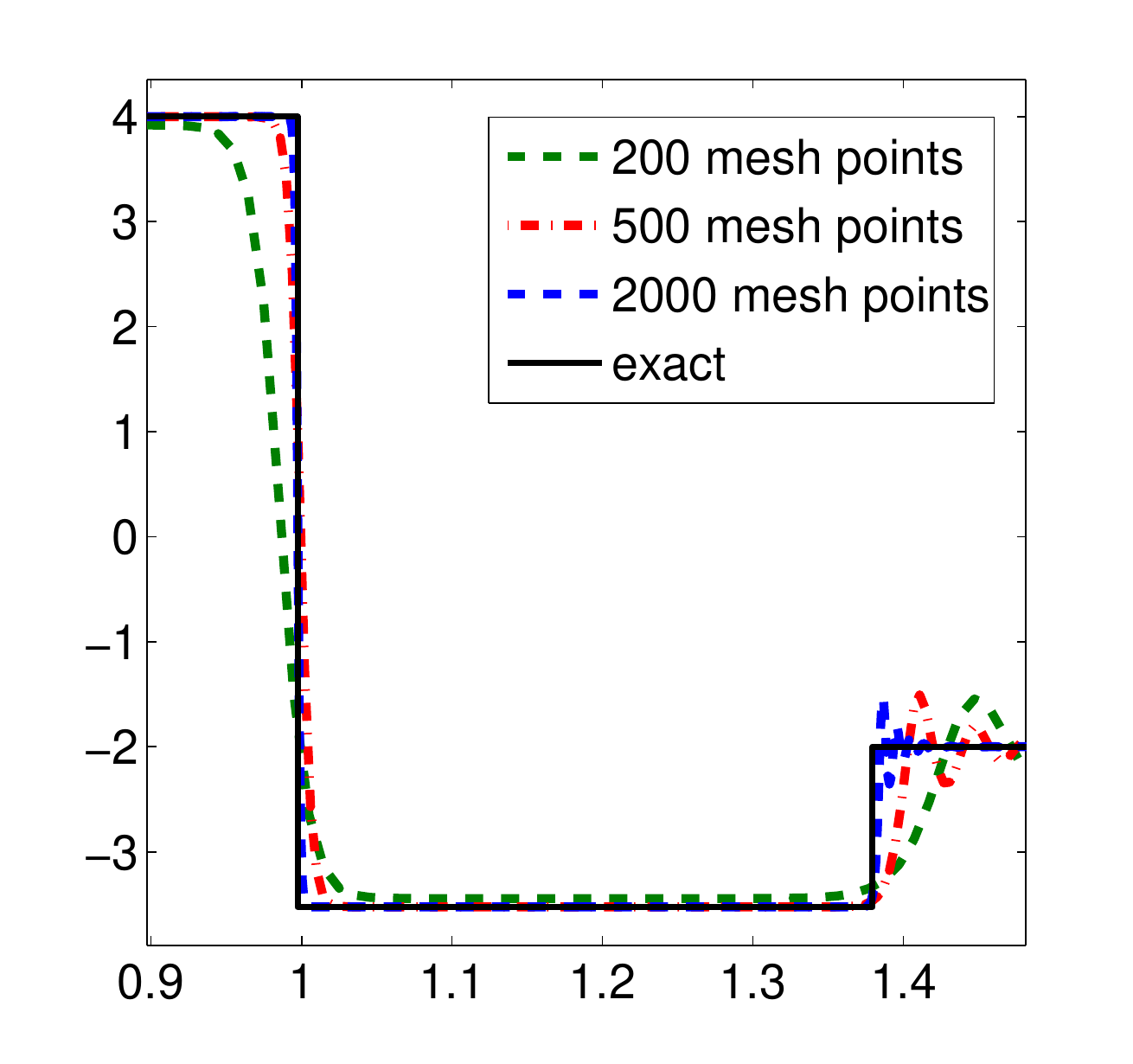}
\end{tabular}
\caption{Convergence (mesh refinement) for the WCD schemes and small shocks. All approximations are based on an order $8$ scheme and $\tau = 0.01$. Left: Riemann solution for $u_L=2$. Right:  Riemann solution for $u_L=4$.}
\label{fig:num1}
\end{figure}

We consider two different sets of $u_L = 2$ and $u_L = 4$ to represent shocks with \emph{small} amplitude. The numerical results are displayed in Figure \ref{fig:num1}, which presents approximate solutions for both sets of initial data, computed with a eighth-order WCD scheme and with $\tau = 0.01$ on a sequence of meshes. For $u_L = 2$, we see that the WCD scheme is able to approximate a nonclassical shock preceeded by a rarefaction wave. Similarly for $u_L = 4$, we see that the WCD scheme approximates both the leading classical shock and the trailing nonclassical shock quite well. In both cases, the quality of approximation improves upon mesh refinement. As expected, there are some oscillations near the leading shock. This is on account of the dispersive terms in the equivalent equation. As shown before, schemes with controlled dissipation were also able to compute small shocks (here the maximum shock strengh is around $7$) quite well.

\subsubsection{Large shocks}

\begin{figure}[htbp]
\centering
\begin{tabular}{cc}
\includegraphics[height=4.5cm, width=0.48\linewidth]{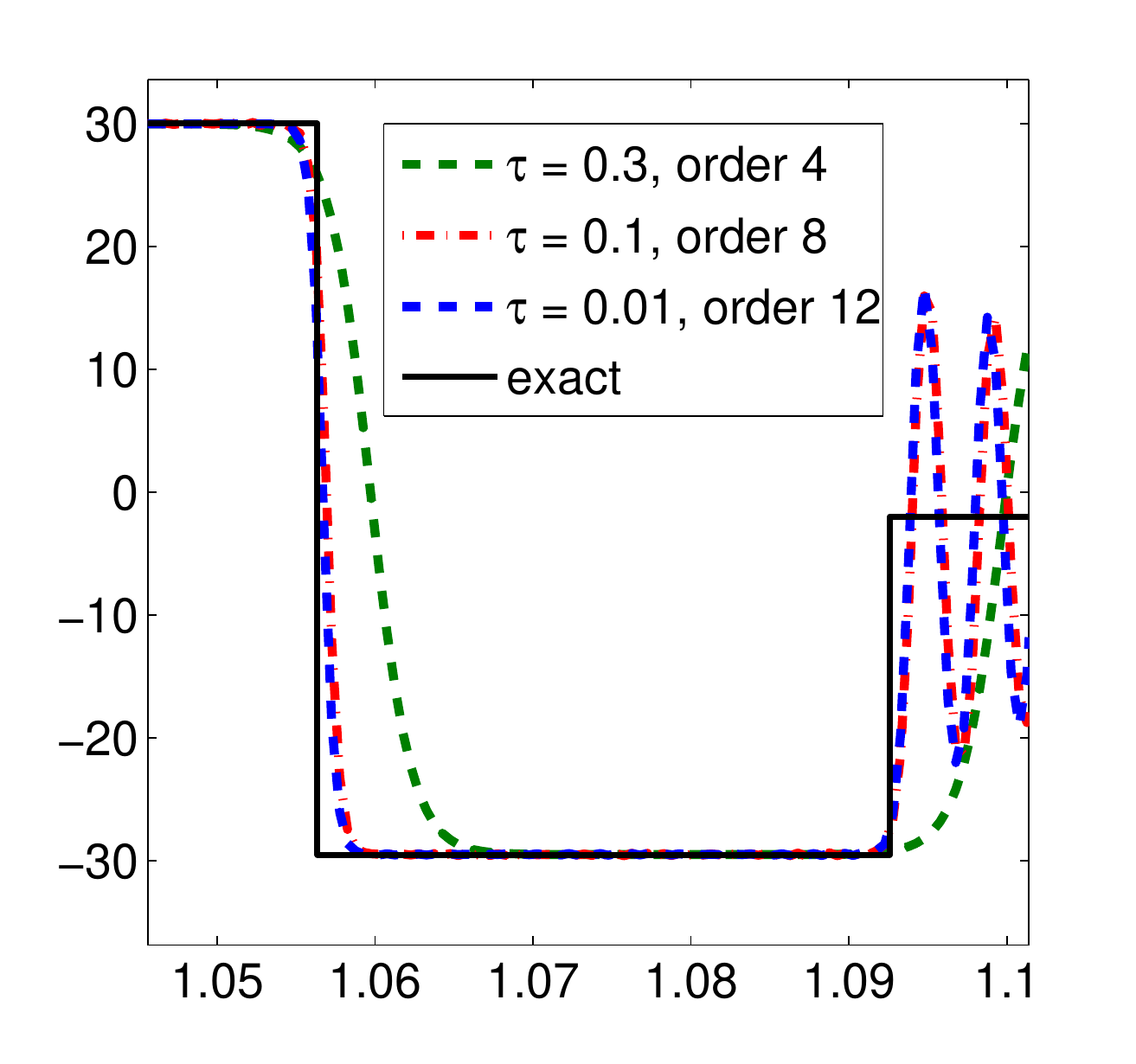} 
& \includegraphics[height=4.5cm, width=0.48\linewidth]{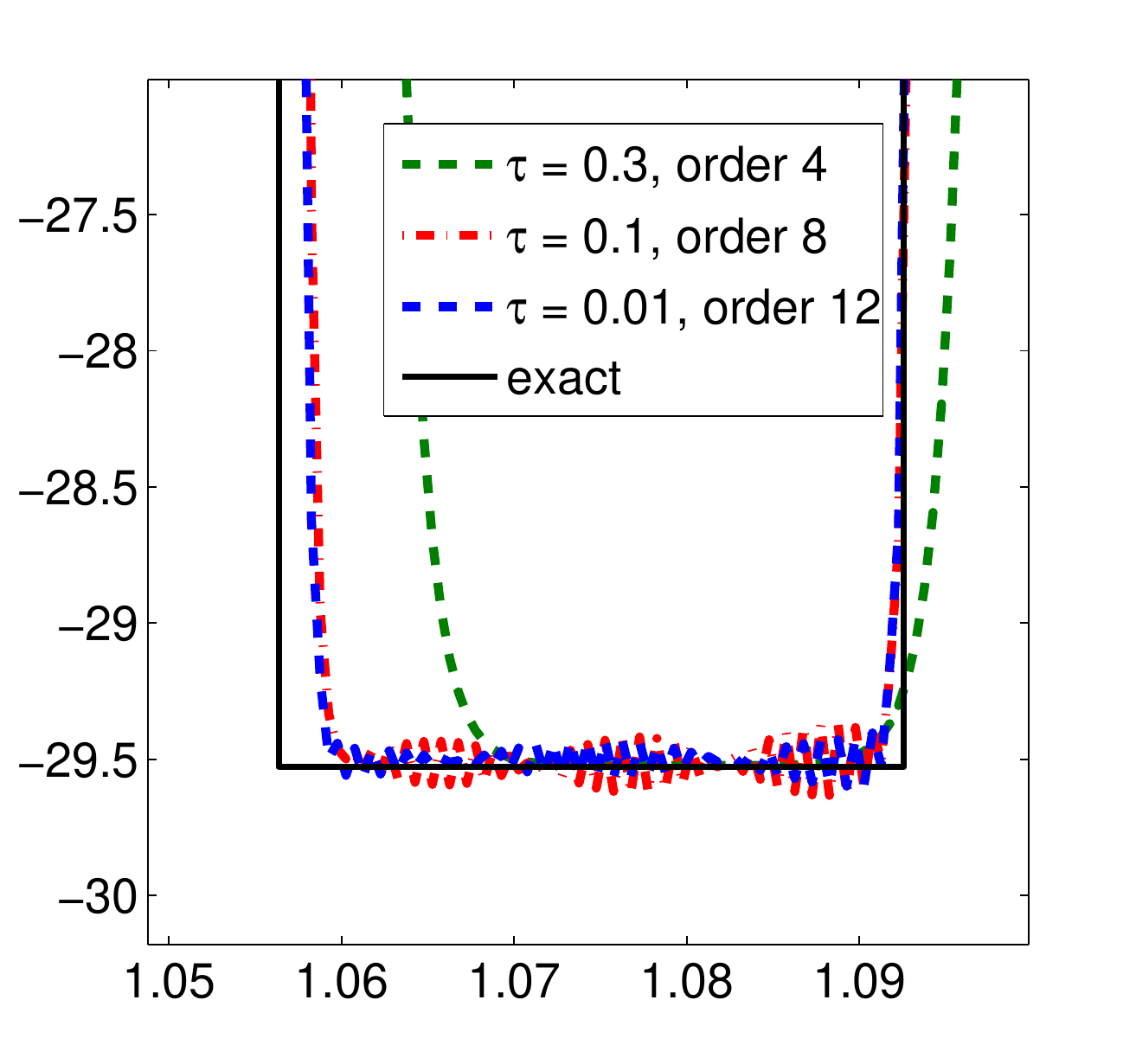}
\end{tabular}
\caption{Convergence $p \to +\infty$ (increasing order of the scheme) for the WCD scheme for a moderate shock. Left: Solution of the Riemann problem for $u_L=30$. Right:  Closer view of the middle state $u_M$.}
\label{fig:num2}
\end{figure}

In order to simulate nonclassical shocks of moderate to large strength, we now choose $u_L = 30$ and display the numerical results in Figure \ref{fig:num2}. The exact solution in this case consists of a leading shock and a trailing nonclassical shock of strength or around $60$ (far stronger than in the previous test with schemes with controlled dissipation that were found to fail for shocks with such strength). In Figure~\ref{fig:num2}, we illustrate how WCD schemes of different order approximate the solution for $4000$ mesh points. A related issue is the variation of the parameter $\tau$. Observe that $\tau$ represent how strong the high-order terms are allowed to be vis a vis the leading order diffusion and dispersion terms. Also, the order of the scheme depends on the choice of $\tau$. For instance, choosing $\tau = 0.1$ implies that fourth-order schemes ($p=2$) are no longer consistent with the WCD condition (\ref{eq:wcd}) for this choice of $\tau$ and one has to use a $sixth$- or even higher order scheme. In this particular experiment, we choose $\tau = 0.3$ (fourth-order scheme), $\tau = 0.1$ (eighth-order) scheme and $\tau = 0.01$ (twelfth-order) scheme. As shown in Figure~\ref{fig:num2}, all the three schemes approximate the nonclassical shock quite well. Also, increasing $\tau$ did not severely affect the shock-capturing abilities of the scheme. Clearly, the eighth and twelfth-order schemes were slightly better in this problem. 

\begin{figure}[htbp]
\centering
\begin{tabular}{cc}
\includegraphics[height=4.5cm, width=0.48\linewidth]{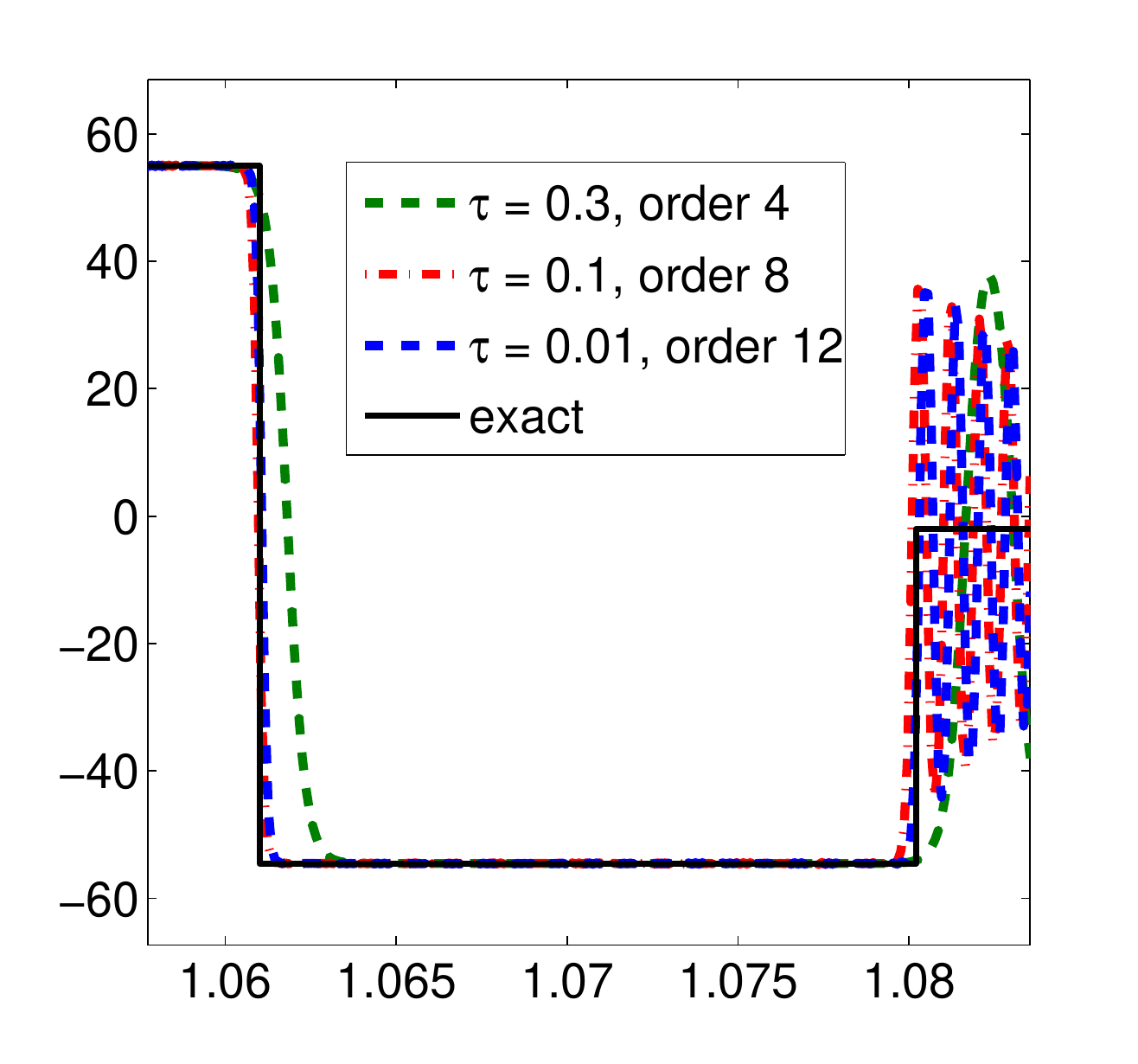} 
& \includegraphics[height=4.5cm, width=0.48\linewidth]{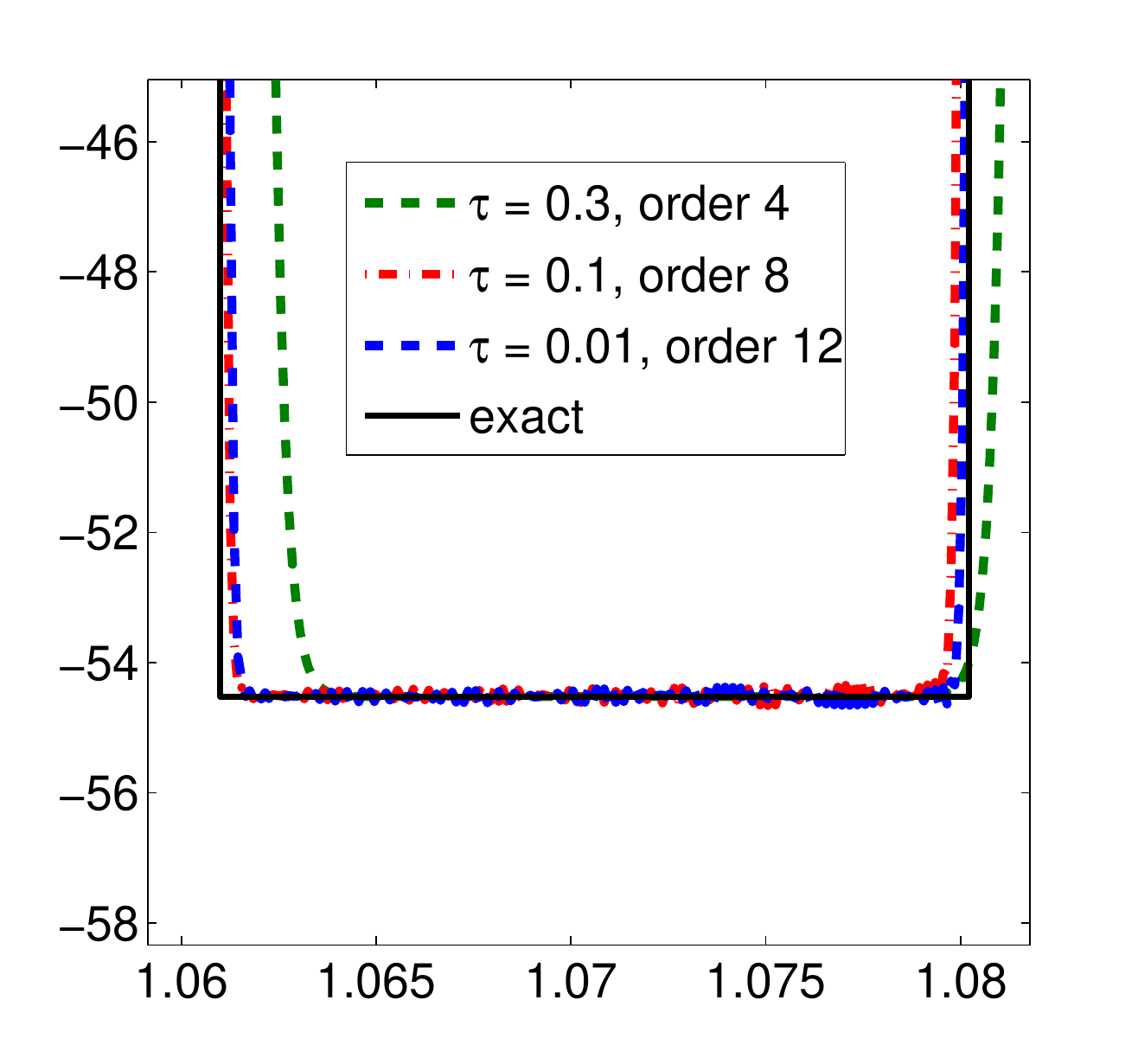}
\end{tabular}
\caption{Convergence $p \to +\infty$ for a large nonclassical shock for the cubic conservation law. Left: Solution of the Riemann problem for $u_L=55$. Right:  Closer view of the middle state $u_M$.}
\label{fig:num4}
\end{figure}

We simulate a very strong shock using $u_L = 55$ in the initial data. The numerical results are presented in Figure~\ref{fig:num4}. The exact solution consists of a strong nonclassical shock of magnitude around $110$ and a weaker (but still of amplitude $60$) leading shock wave. The results in the figure were generated with the fourth, eighth and twelfth order schemes. The mesh resolution is quite fine ($20000$ mesh points) as the difference in speeds for both shocks is quite small, the intermediate state is very narrow and needs to be resolved. It is important to emphasize that one can easily use a grid, adapted to the shock locations. The results clearly show that all the three schemes converge to the correct nonclassical shock, even for such a strong shock.

\subsubsection{Computing the kinetic relation}
\begin{figure}[htbp]
\centering
\includegraphics[height=4.5cm, width=0.48\linewidth]{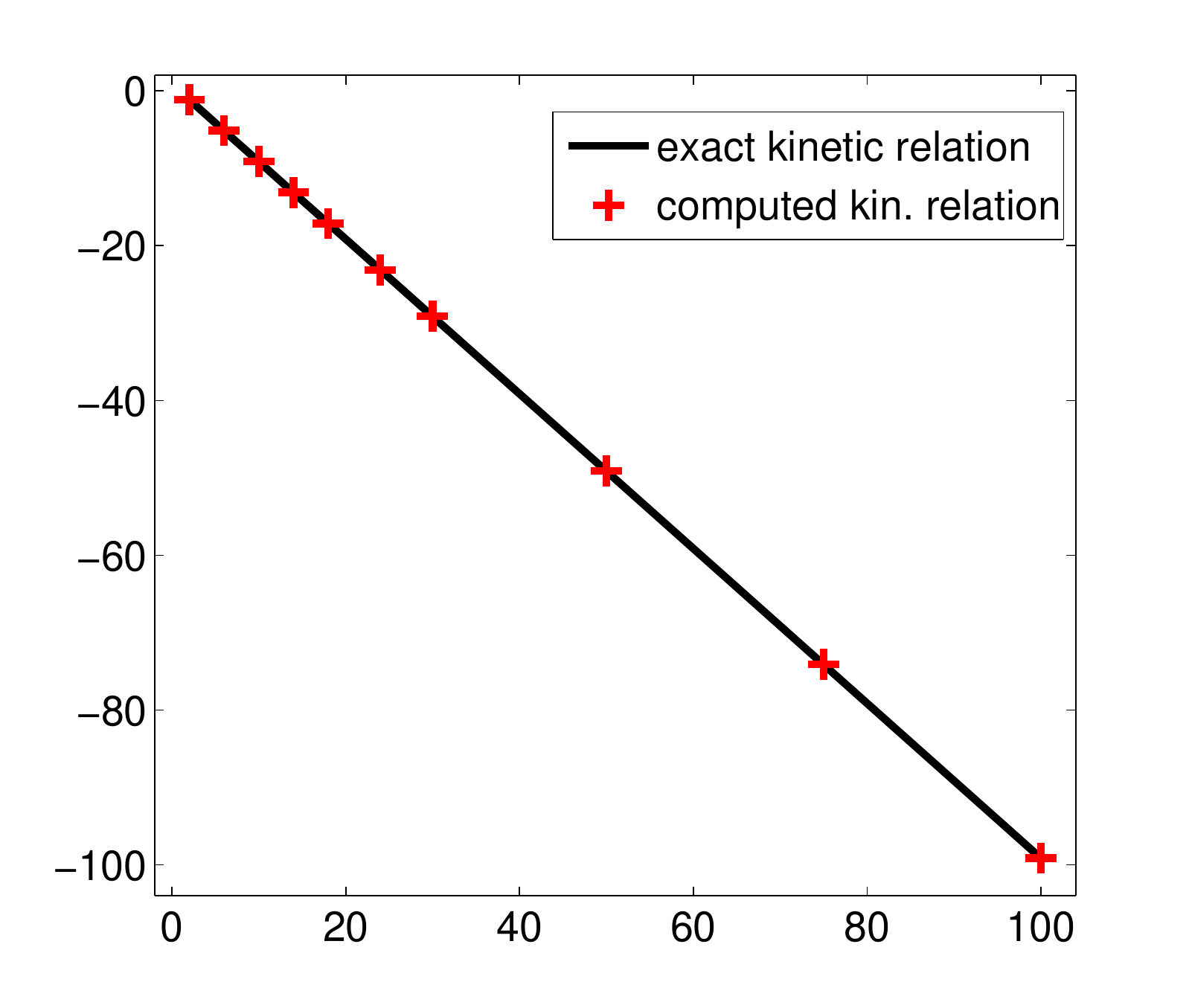}
\includegraphics[height=4.5cm, width=0.48\linewidth]{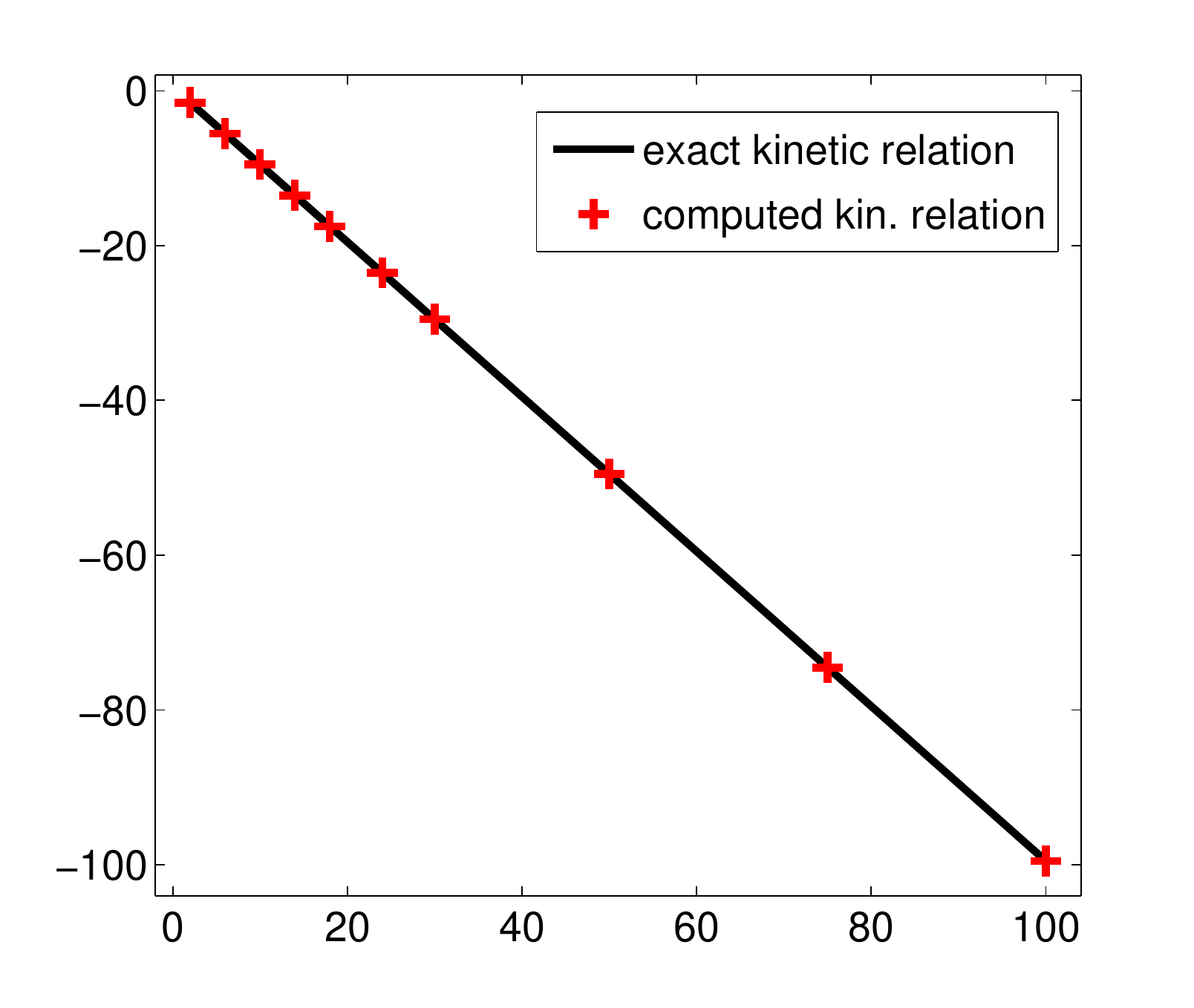}

\includegraphics[height=4.5cm, width=0.48\linewidth]{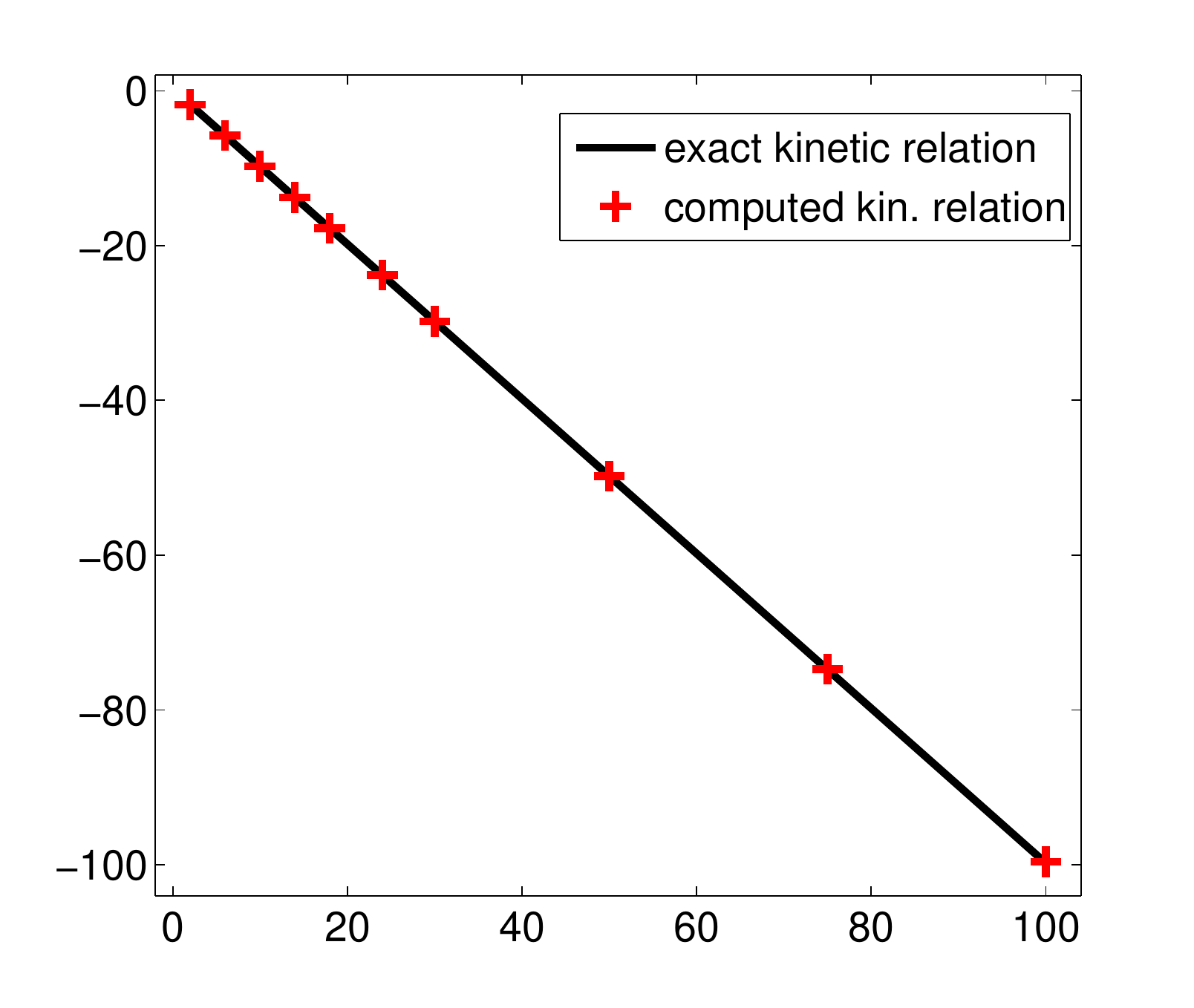}
\caption{Intermediate state (kinetic function) (Y-axis) for the cubic conservation law for varying left-hand state $u-L$ (X-axis), computed using an eighth-order WCD scheme with $\tau = 0.1$. Left: Kinetic function for $\delta=0.3$. Middle: Kinetic function for $\delta=1$. Right: Kinetic function for $\delta=5$.}
\label{fig:num5}
\end{figure}

 Since, the exact intermediate state of a Riemann problem is known for the cubic conservation law (for any given value of the dispersion parameter $\delta$), we can ascertain the  quality and accuracy of numerical approximation for a very large class of initial data by computing the numerical kinetic relation. We do so using the eighth-order WCD scheme for three different values of the dispersion parameter $\delta$. The results are presented in Figure~\ref{fig:num5} and clearly demonstrate that this WCD scheme is able to compute the correct intermediate state (kinetic relation) and, hence, the nonclassical shock wave, accurately for any given shock strength. In particular, very strong nonclassical shocks are captured accurately. Similar results were also obtained with WCD schemes of different orders. These results should be compared with earlier work on numerical kinetic functions (Hayes and LeFloch 1997, 1998, LeFLoch and Rohde 2000, LeFloch and Mohamadian 2008), where, despite the convergence $p \to +\infty$ being observed for each fixed shock strength, 
 the numerical kinetic function was found to significantly differ from the analytical one for large shocks. Furthermore, the WCD are very remarkable in that they even capture the {\sl correct asymptotic behavior} of the kinetic function in the limit of {\sl arbitrary large} shock strength. 
 
%------------------------------------------------------------------------------------------------------------- 

\subsection{WCD schemes for nonlinear hyperbolic systems}
\label{sec:33}

Finite difference schemes with either controlled or well-controlled dissipation can be readily extended to systems of conservation laws, now discussed. Again, for the sake of presentation, we focus on an example, specifically the nonlinear elasticity system \eqref{555} that models viscous-capillary flows of elastic materials. We rewrite this system in the general form 
\be
\textbf{U}_t + \textbf{F}_x = \epsilon D^{(1)} \textbf{U}_{xx} + \alpha \epsilon^2 D^{(2)} \textbf{U}_{xxx}, 
\label{eq:vdWsyst}
\ee
in which we have set 
$$
\textbf{U} = \left(\begin{array}{cc}  w \\  v \end{array} \right), \,\,\,\,\,\,\,\,\,\, D^{(1)} = \left( \begin{array}{cc}0&0\\0&1\end{array}\right), \,\,\,\,\,\,\,\,\,\, D^{(2)} = \left( \begin{array}{cc}0&0\\-1&0\end{array}\right)
$$
and the flux vector $\textbf{F}: R^2 \to R^2$ is given by
$
	\textbf{F}(\tau, u) = \left(\begin{array}{cc}-v\\ -\sigma(w) \end{array} \right).
$

We consider a uniform grid as was described in the previous section and a $2p$-th order accurate, finite difference scheme for  (\ref{eq:vdWsyst}) reads 
\begin{eqnarray}
\label{eq:fds2p_vdW}
&&
\frac{d \textbf{U}_i}{dt} + \frac{1}{\Delta x} \sum\limits_{j=-p}^{j=p} \alpha_j \textbf{F}_{i+j}
 \\
\nonumber                                         
&&=
 \frac{c}{\Delta x} \sum\limits_{j=-p}^{j=p} \beta_j D^{(1)} \textbf{U}_{i+j} 
+  \frac{\delta c^2}{\Delta x} \sum\limits_{j=-p}^{j=p} \gamma_j D^{(2)} \textbf{U}_{i+j},  
\end{eqnarray}
where $\textbf{U}_i = \textbf{U}(x_i,t)$, $\textbf{F}_i = \textbf{F}(\textbf{U}_i)$ and the coefficients $\alpha_j, \beta_j$ and $\gamma_j$ need to satisfy the order conditions (\ref{eq:al})-(\ref{eq:ga}). As explained in the previous subsection, 
the key tool in designing a scheme that can accurately approximate nonclassical shocks to (\ref{eq:vdWsyst}) 
is the equivalent equation associated the scheme (\ref{eq:fds2p_vdW}), which reads 
\be
\label{eq:eesystems}
\aligned
\frac{d \textbf{U}}{dt} &= 
\underbrace{- \textbf{F}_x + c\Delta x D^{(1)} \textbf{U}_{xx} + \delta c^2 \Delta x^2 D^{(2)} \textbf{U}_{xxx}}_{l.o.t} + \text{h.o.t.},
\\
\text{h.o.t.} &=
 - \sum\limits_{k=2p+1}^{\infty} \frac{\Delta x^{k-1}}{k!} A^p_k \textbf{F}^{[k]}  
                                        +  c \sum\limits_{k=2p+1}^{\infty} \frac{\Delta x^{k-1}}{k!} B^p_k D^{(1)} \textbf{U}^{[k]} 
\\
& \quad +  \delta c^2 \sum\limits_{k=2p+1}^{\infty} \frac{\Delta x^{k-1}}{k!} C^p_k D^{(2)} \textbf{U}^{[k]}, 
\endaligned
\ee
the coefficients $A^p_k, B^p_k$ and $C^p_k$ being defined as in (\ref{eq:coeff}). 

Following our discussion in the previous subsection, our design of a WCD scheme is based on the analysis of a single shock. Adapting from the scalar case, we impose a component-wise condition in order to balance high-order and low-order terms in the equivalent equation for a single shock. We thus assume a tolerance $\tau$ such that 
$|h.o.t.| \leq \tau \, |l.o.t.|$
holds componentwise. This analysis is carried out in Ernest, LeFloch, and Mishra (2013) and results in the following {\sl WCD condition for systems} ($i=1,2$): 
\be
\label{eq:wcdsystem}
\begin{aligned}
({\sl WCD})_i: \quad 
& \left( |\delta|- \frac{\widehat{S}^C_p |\delta|}{\tau}\right) |\left\langle D^{(2)}_{i},[[\textbf{U}]]\right\rangle| c_i^2 + \left(1 - \frac{\widehat{S}^D_p}{\tau}\right) |\left\langle D^{(1)}_{i},[[\textbf{U}]]\right\rangle| c_i 
\\
& - \left(1 + \frac{\widehat{S}^f_p}{\tau}\right)\sigma |[[\textbf{U}_i]]| > 0. 
\end{aligned}
\ee
Here, $\widehat{S}^D_p,\widehat{S}^C_p,\widehat{S}^f_p$ and $\sigma$ are defined as in (\ref{eq:s5}) and $\left\langle \cdot, \cdot\right\rangle $ denotes the product of two vectors and 
\be
\label{eq:defsig}
\sigma = \frac{|[[\textbf{F(U)}]]|}{|[[\textbf{U}]]|}, 
\ee
being a rough estimate on the maximum shock speed of the system with jump $[[\textbf{U}]]$ at the single shock. Observe that if $\left\langle D^{(1)}_{i},[[\textbf{U}]]\right\rangle = 0$, we set the corresponding $c_i = 0$ for $i = 1,2$. The scheme parameter $c$ in (\ref{eq:fds2p_vdW}) is defined as $c = \max(c_1, c_2)$. 
The global definition of the coefficient $c$ can be obtained by taking a maximum of the afore obtained $c$ over all cells.

%----------------------------------------------------------------------------

\subsubsection{Numerical experiments} 
\label{sec:34}

\begin{figure}[htbp]
\centering
\begin{tabular}{cc}
\includegraphics[width=0.46\linewidth]{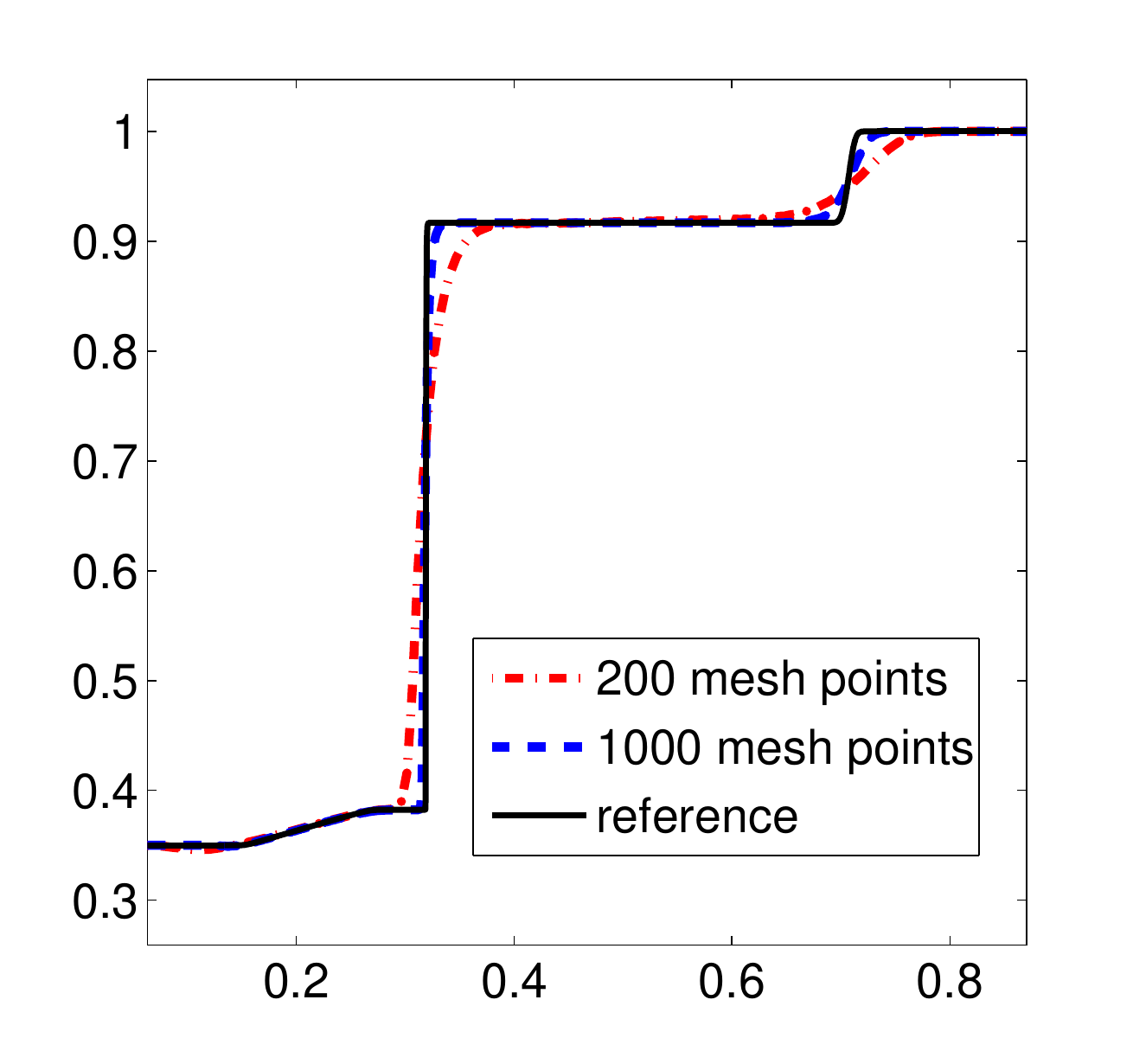} 
& \includegraphics[width=0.46\linewidth]{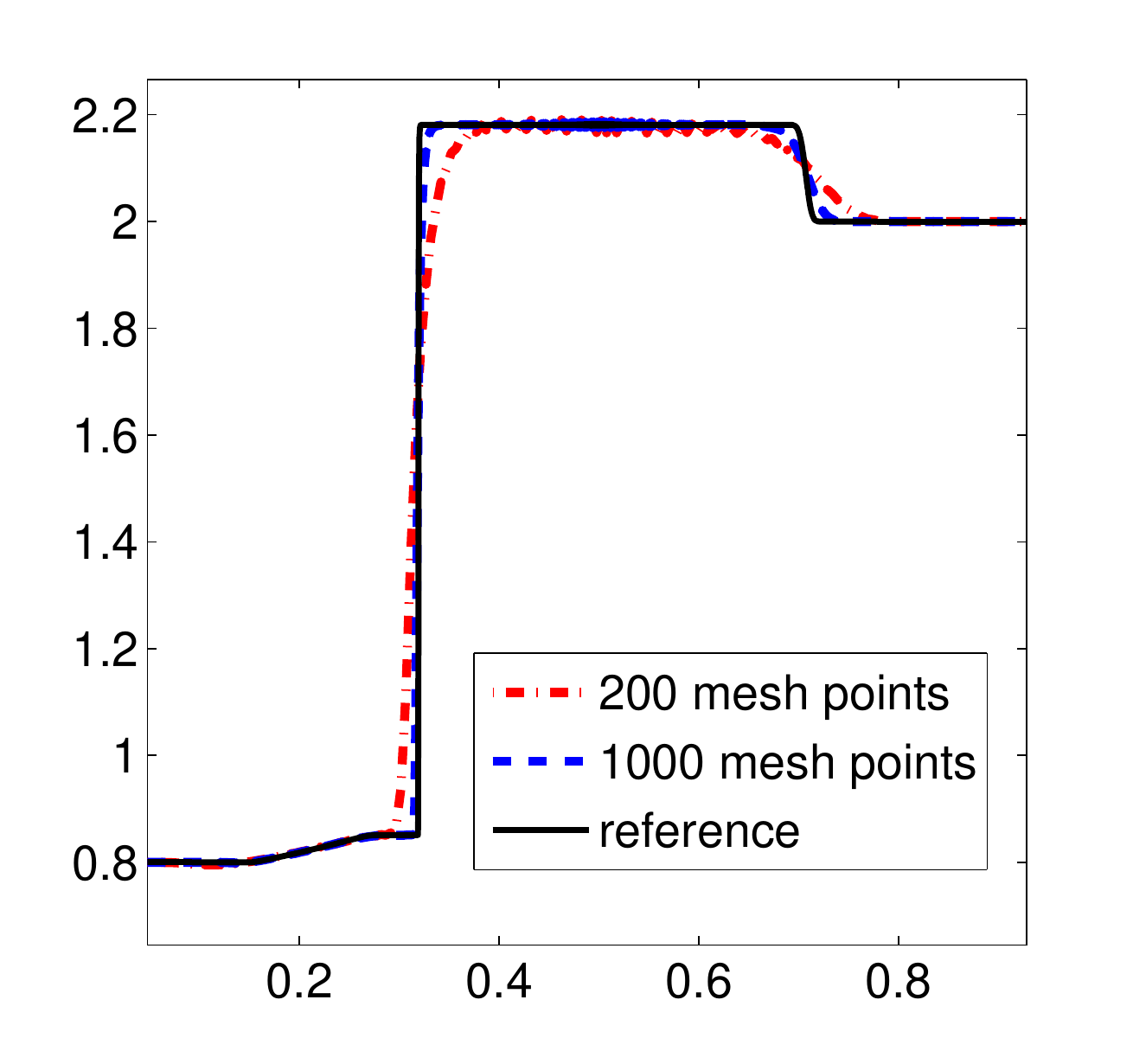}
\end{tabular}
\caption{Mesh-convergence for the WCD scheme for the dispersive limit of van der Waals fluid with initial data (\ref{eq:RP_vdW}) with a eighth-order WCD scheme. 
\newline 
Left: Velocity component $v$. Right: Volume component $w$.}
\label{fig:num_vdW3}
\end{figure}

\begin{figure}[htbp]
\centering
\begin{tabular}{cc}
\includegraphics[width=0.47\linewidth]{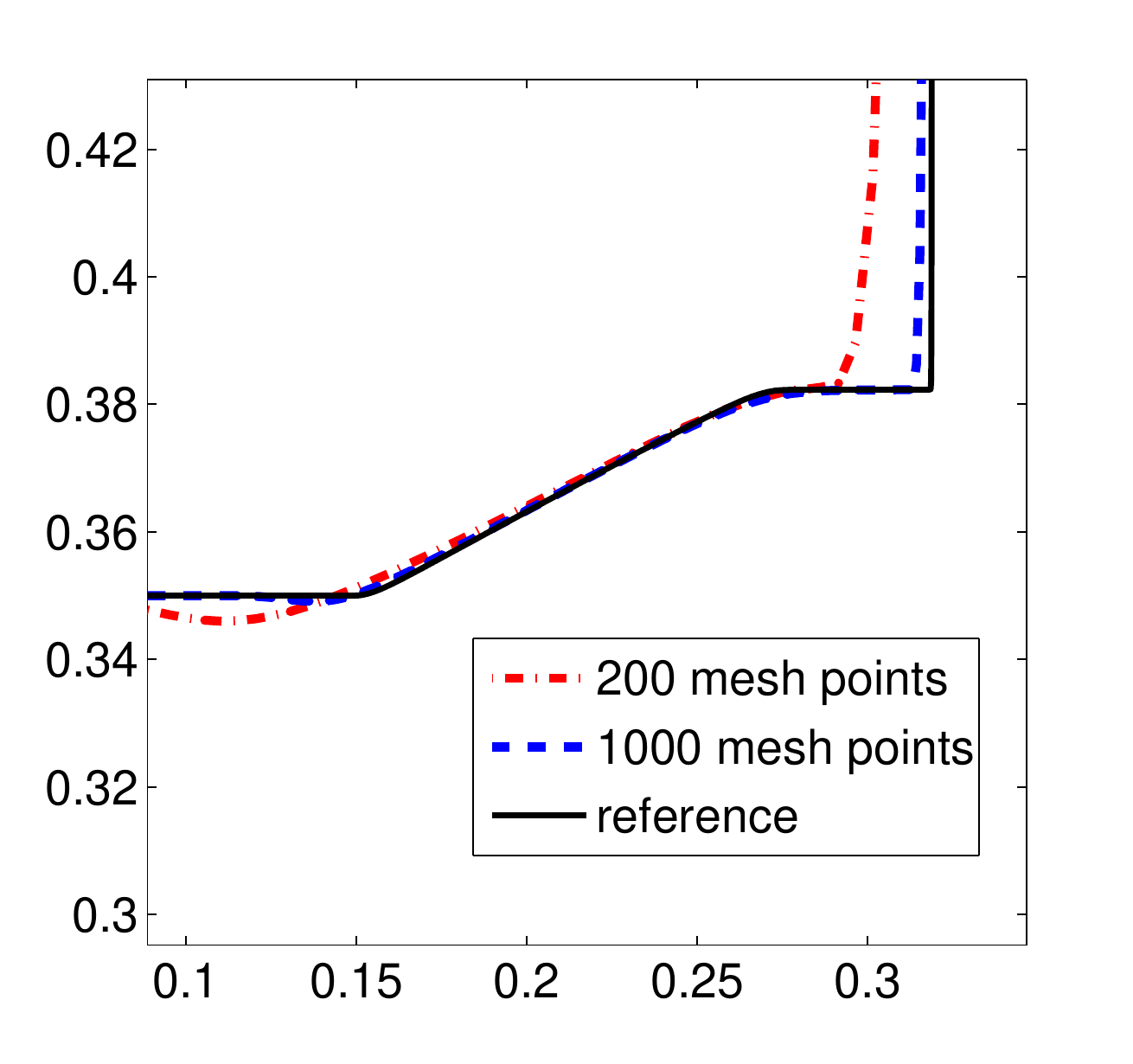} 
& \includegraphics[width=0.47\linewidth]{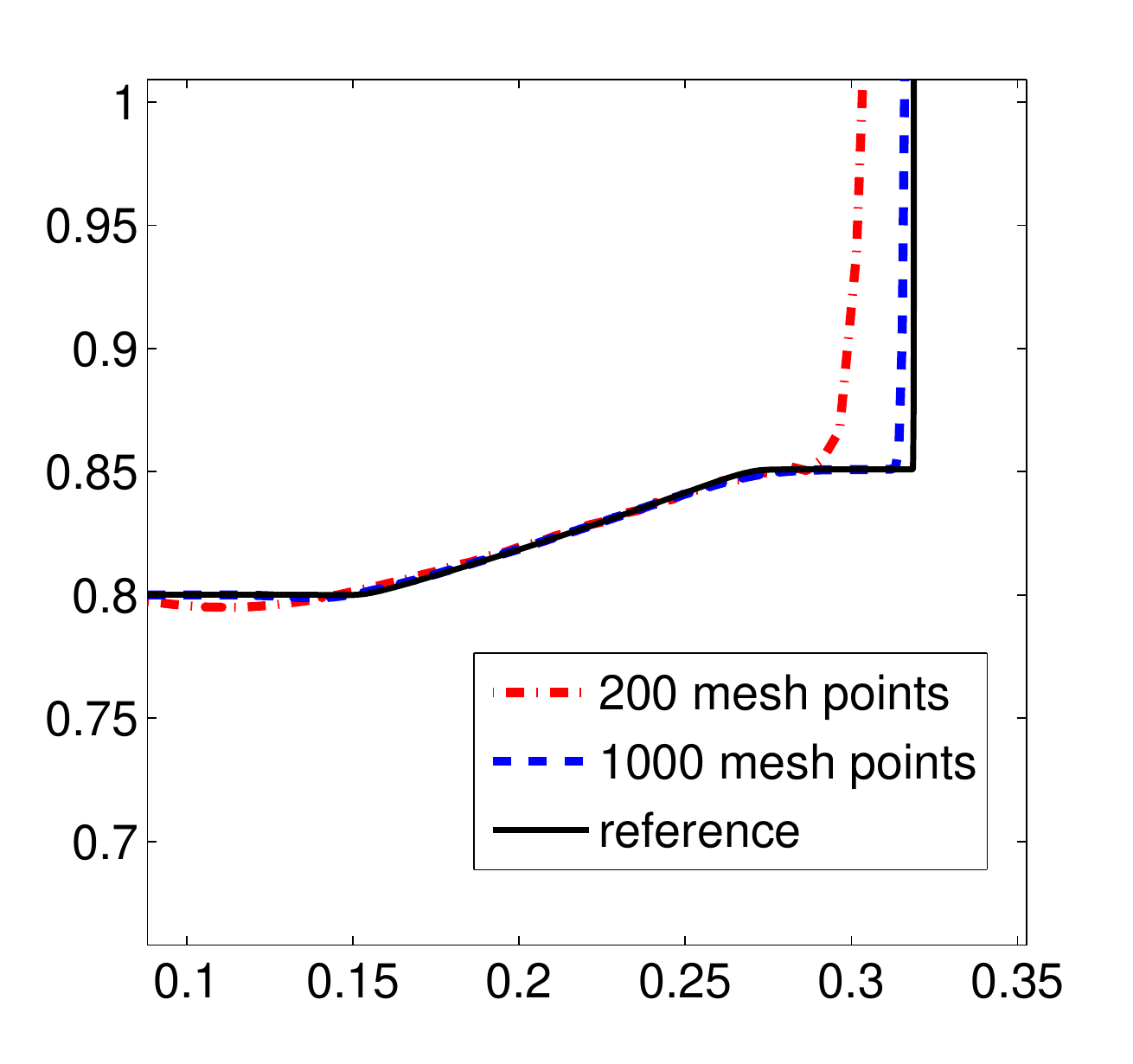}
\end{tabular}
\caption{Zoom near nonclassical states for the dispersive limit of van der Waals fluid with initial data (\ref{eq:RP_vdW}) and a eighth-order WCD scheme. 
\newline 
Left: Velocity component $v$. Right: Volume component $w$.}
\label{fig:num_vdW4}
\end{figure}

In our numerical tests, we use the normalized van der Waals flux given by
$$
\sigma(w) := -\frac{RT}{\left(w - \frac{1}{3} \right) } - \frac{3}{w^2}
$$
with $R = \frac{8}{3}$ and $T = 1.005$ which has two inflection points at $1.01$ and $1.85$. With this choice of parameter, the elasticity system \eqref{555} is strictly hyperbolic. We let $\alpha = 1$ and consider the initial Riemann data 
\be
	v(0, x) = \left\lbrace \begin{array}{cc} 0.35, & x < 0.5,
 \\ 1.0, & x > 0.5,	 \end{array} \right. \qquad \qquad
	w(0, x) = \left\lbrace \begin{array}{cc} 0.8, & x < 0.5,
 \\ 2.0, & x > 0.5,	 \end{array} \right.
\label{eq:RP_vdW}
\ee
and the scheme parameter is set to $c = c_{WCD}$. Figures \ref{fig:num_vdW3} and  \ref{fig:num_vdW4} show a nonclassical state in both variables $v$ and $w$ and displays mesh convergence of the scheme as the mesh is refined. Furthermore, the eighth-order WCD scheme approximates the nonclassical state quite well for both variables, even at a very coarse mesh resolution. 

\subsubsection{Large shocks}

\begin{figure}[htbp]
\centering
\begin{tabular}{cc}
\includegraphics[height=4.5cm, width=0.48\linewidth]{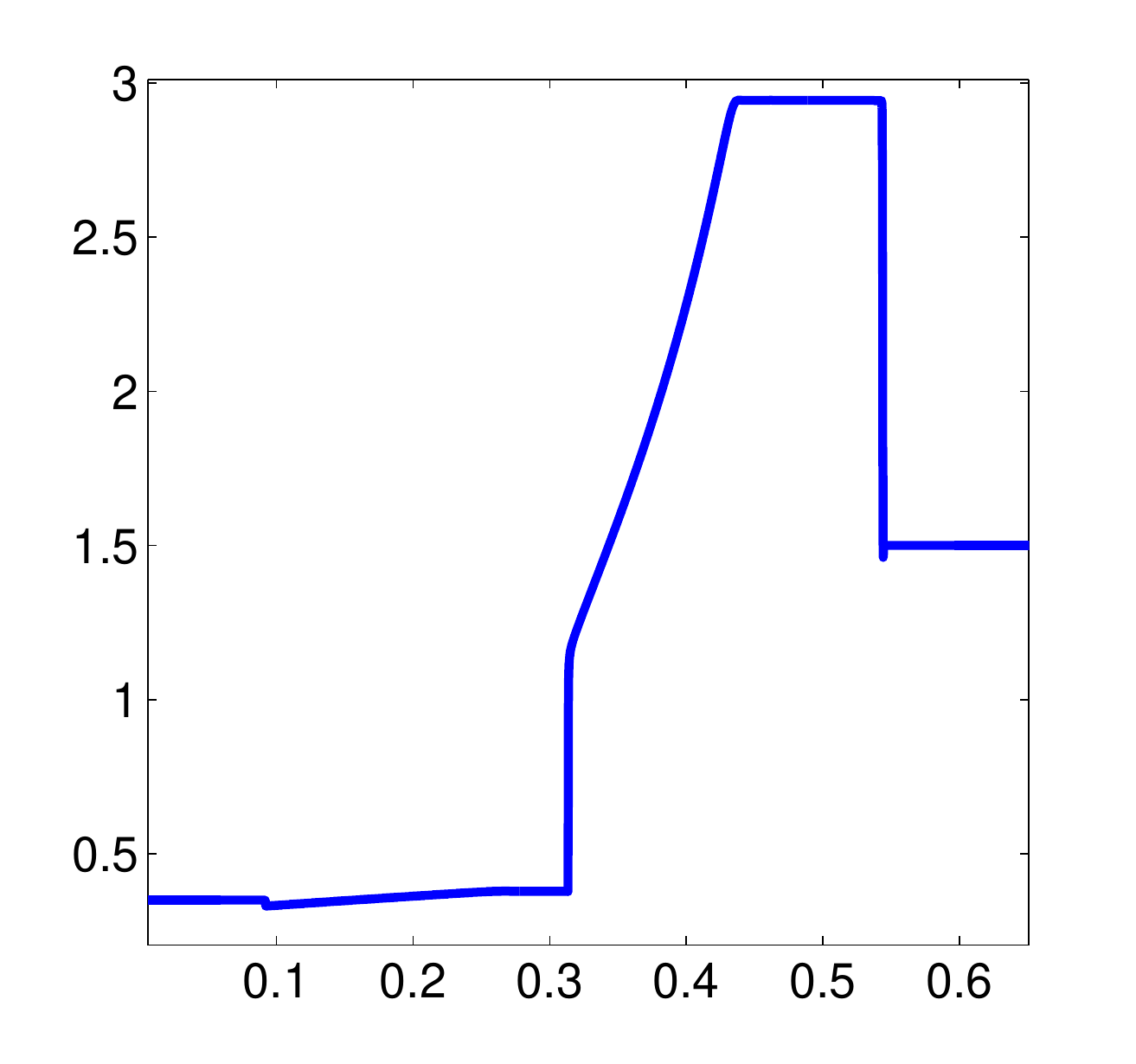} 
& \includegraphics[height=4.5cm, width=0.48\linewidth]{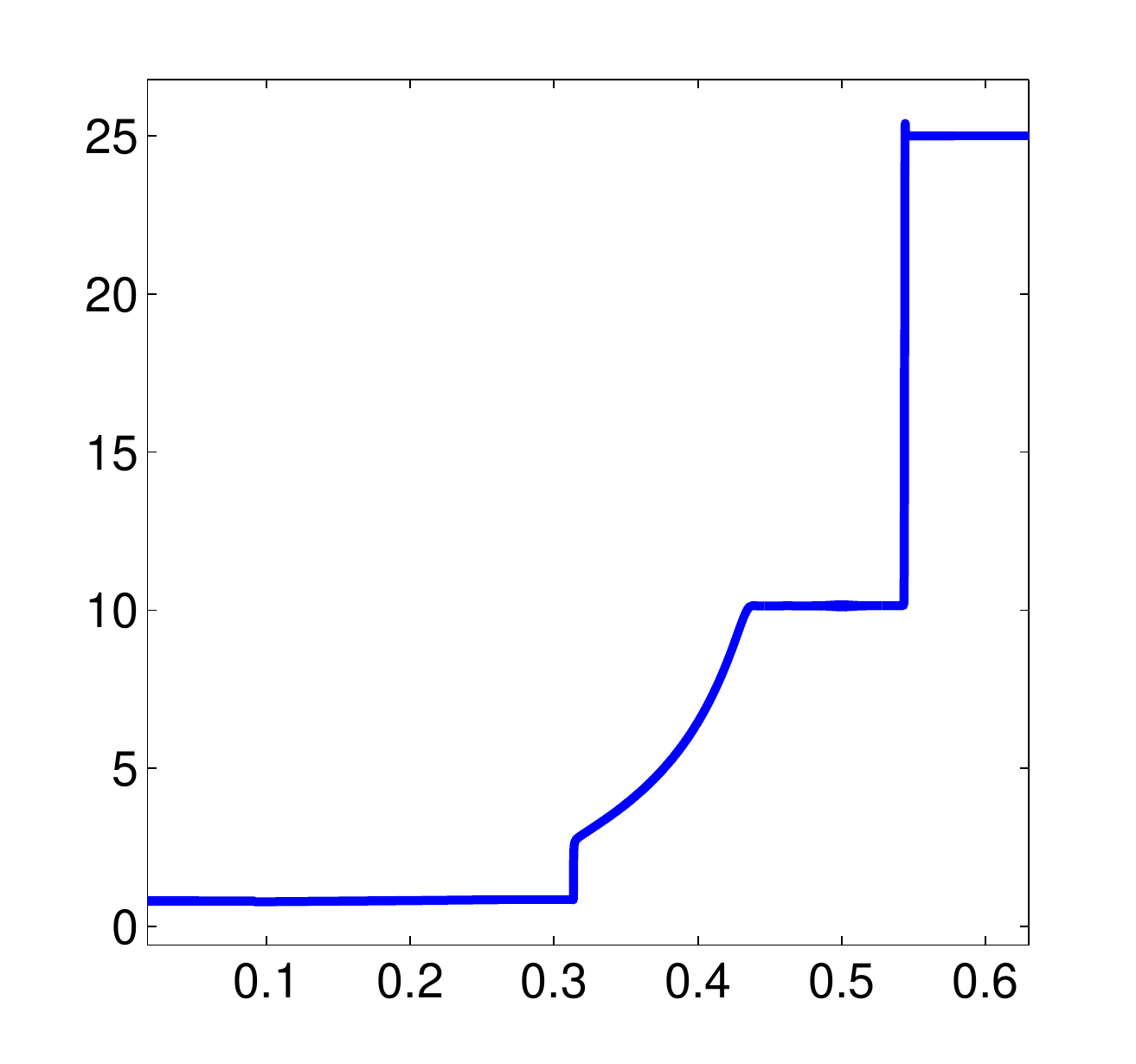}
\end{tabular}
\caption{Approximation of the Riemann solution with initial data (\ref{eq:RP}) resulting in a large jump in volume $\tau$. 
\newline 
Left: Velocity component $u$. Right: Volume component $\tau$.
\newline 
Results obtained with an eighth order WCD scheme with $25000$ points.}
\label{fig:num_vdW5}
\end{figure}

\begin{figure}[htbp]
\centering
\begin{tabular}{cc}
\includegraphics[height=4.5cm, width=0.48\linewidth]{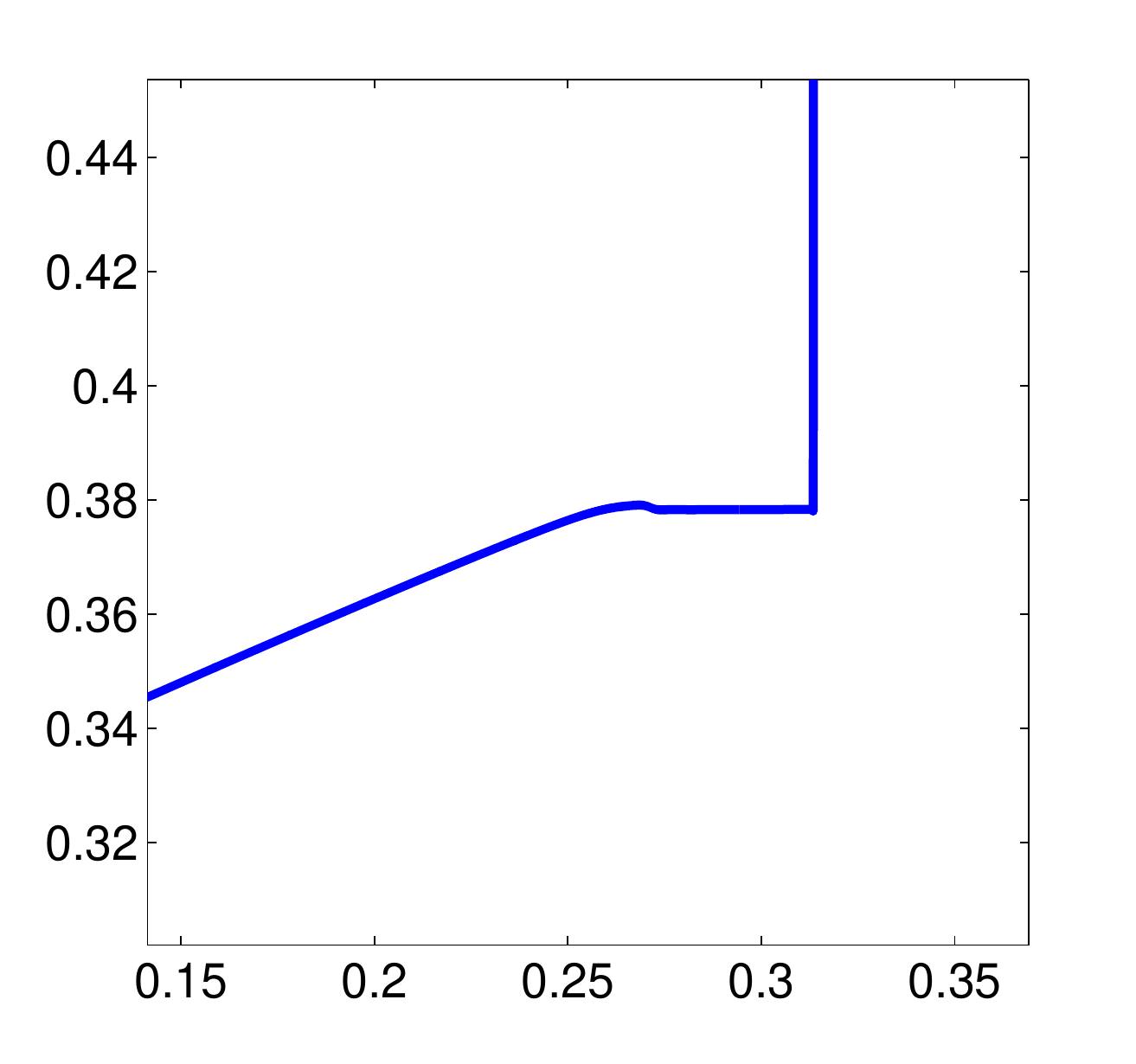} 
& \includegraphics[height=4.5cm, width=0.48\linewidth]{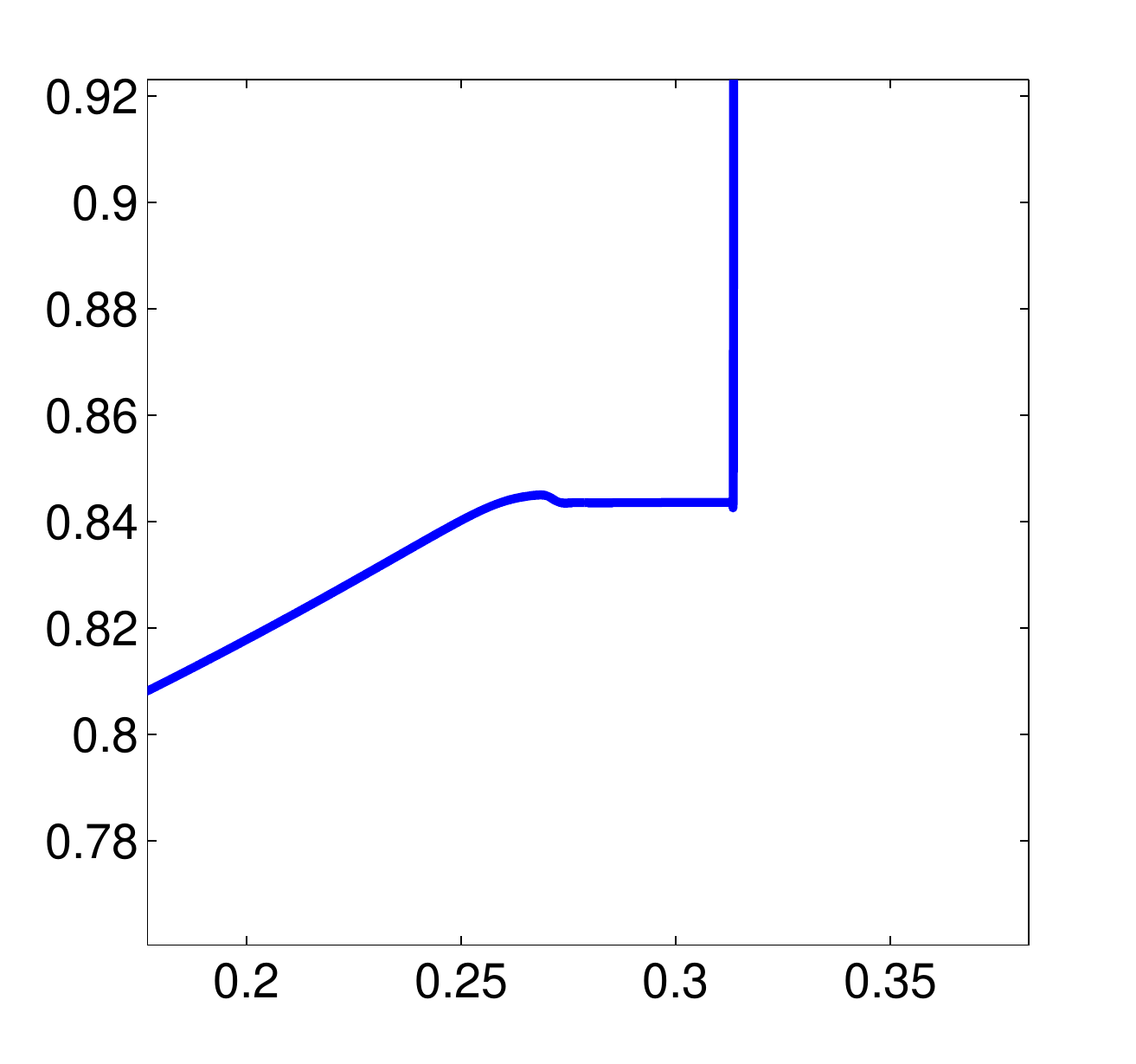}
\end{tabular}
\caption{Zoom at nonclassical states in the approximation of the Riemann problem with initial data (\ref{eq:RP}) resulting in a large jump in volume $\tau$. 
\newline 
Left: Velocity component $u$. Right: Volume component $\tau$.
\newline  Results obtained with an eighth order WCD scheme with 25000 points.}
\label{fig:num_vdW6}
\end{figure}

Next, we approximate large nonclassical shocks associated with the Riemann initial data
\be
	v(0, x) = \left\lbrace \begin{array}{cc} 0.35, & x < 0.5, \\ 1.5, & x > 0.5, 	 \end{array} \right. \qquad \qquad
	w(0, x) = \left\lbrace \begin{array}{cc} 0.8, & x < 0.5,
 \\ 25.0, & x > 0.5.	 \end{array} \right.
	\label{eq:RP}
\ee
The results with a eighth-order scheme are shown in Figures~\ref{fig:num_vdW5} and \ref{fig:num_vdW6} and clearly show that the WCD scheme is able to approximate the nonclassical shock of large amplitude (in the volume) quite well. 

%------------------------------------------------------------------------------------------------------

\subsection{A model of magnetohydrodynamics with Hall effect} 

Next, we consider the simplified model of ideal magnetohydrodynamics 
\be
\label{834}
\begin{aligned}
{v_t + \big( (v^2 + w^2) \, v\big)_x} & = \eps \, v_{xx} + \alpha \, \epsilon \, w_{xx}, 
\\
{w_t + \big( (v^2 + w^2) \,w \big)_x} & = \eps \, w_{xx} - \alpha \, \epsilon \, v_{xx},
\end{aligned}
\ee
where $v,w$ denote the transverse components of the magnetic field, $\eps$ the magnetic resistivity, and $\alpha$ the so-called Hall parameter. The Hall effect is relevant in order to investigate, for instance, the solar wind interaction with the Earth's magnetosphere. The viscosity-only regime $\alpha=0$ was studied by Brio and Hunter (1990), Freist\"uhler (1992), and Freist\"uhler and Pitman (1992, 1995).  The left-hand part of the equations \eqref{834} form a hyperbolic but a non-strictly hyperbolic system of conservation laws.   
 Furthermore, observe that solutions to \eqref{834} satisfy the identity 
$$
\aligned
{1 \over 2} \big(v_\eps^2 + w_\eps^2\big)_t + {3 \over 4} \big((v_\eps^2 + w_\eps^2)^2\big)_x 
=
& - \eps \, \big( (v^\eps_x)^2 + (w^\eps_x)^2 \big) + \eps \, C^\eps_x,  
\endaligned
$$ 
so that in the formal limit $\eps \to 0$ the following {\sl quadratic entropy inequality} holds for 
 the magnetohydrodynamic model:  
\be
{1 \over 2} \big(v^2 + w^2\big)_t + {3 \over 4} \big((v^2 + w^2)^2\big)_x \leq 0. 
\ee

Following Ernest, LeFloch, and Mishra (2013), we design finite difference schemes with well-controlled dissipation for \eqref{834}. We first 
 rewrite it in the general form:
\be
\textbf{U}_t + \textbf{F}_x = \epsilon D^{(1)} \textbf{U}_{xx} + \alpha \epsilon D^{(2)} \textbf{U}_{xx}, 
\label{eq:MHDsyst}
\ee
where the vector of unknowns is ${\sl U} = \{v,w\}$ and
$$
 D^{(1)} = \left( \begin{array}{cc}1&0\\0&1\end{array}\right), \,\,\,\,\,\,\,\,\,\, D^{(2)} = \left( \begin{array}{cc}0&1\\-1&0\end{array}\right)
$$
and the flux $\textbf{F}: R^2 \to R^2$ reads 
$\textbf{F}(v, w) = \left(\begin{array}{cc}(v^2+w^2)v\\ (v^2 + w^2)w \end{array} \right).
$
We then introduce a uniform grid as in previous subsections and, for any integer $p \geq 1$, we approximate (\ref{eq:MHDsyst}) with the $2p$-th order consistent, finite difference scheme 
\be
\label{eq:fds2p_MHD}
\aligned
& \frac{d \textbf{U}_i}{dt} + \frac{1}{\Delta x} \sum\limits_{j=-p}^{j=p} \alpha_j \textbf{F}_{i+j} 
\\
&= \frac{c}{\Delta x} \sum\limits_{j=-p}^{j=p} \beta_j D^{(1)} \textbf{U}_{i+j} 
+  \frac{\alpha c}{\Delta x}\sum\limits_{j=-p}^{j=p} \beta_j D^{(2)} \textbf{U}_{i+j}, 
\endaligned
\ee
where $\textbf{U}_i = \textbf{U}(x_i,t)$, $\textbf{F}_i = \textbf{F}(\textbf{U}_i)$ and the coefficients $\alpha_j $ and $\beta_j$ need to satisfy the order conditions (\ref{eq:al})-(\ref{eq:ba}). 

The equivalent equation associated with the scheme (\ref{eq:fds2p_MHD}) reads 
\be
\label{eq:eeMHD}
\aligned
\frac{d \textbf{U}}{dt} &= \underbrace{- \textbf{F}_x +c\Delta x D^{(1)} \textbf{U}_{xx} + \alpha c \Delta x^2 D^{(2)} \textbf{U}_{xx}}_{l.o.t} + \text{h.o.t.},
\\
\text{h.o.t.} &=  
 - \sum\limits_{k=2p+1}^{\infty} \frac{\Delta x^{k-1}}{k!} A^p_k \textbf{F}^{[k]} 
                                        +  c \sum\limits_{k=2p+1}^{\infty} \frac{\Delta x^{k-1}}{k!} B^p_k D^{(1)} \textbf{U}^{[k]} 
\\
& \qquad +  \alpha c \sum\limits_{k=2p+1}^{\infty} \frac{\Delta x^{k-1}}{k!} B^p_k D^{(2)} \textbf{U}^{[k]}.
\endaligned
\ee
The coefficients $A^p_k$ and $B^p_k$ are defined as in (\ref{eq:coeff}). As in the cases of scalar conservation laws and nonlinear elasticity models, we can analyze the equivalent equation at a single shock with jump $[[\textbf{U}]]$ and then follow the steps of the previous subsections. This analysis was carried out by Ernest, LeFloch, and Mishra (2013) and and led them to the following WCD condition:
\be
\label{eq:wcdMHD} 
\begin{aligned}
({\sl WCD})_i: 
& \left(\left( |\alpha|- \frac{\widehat{S}^D_p |\alpha|}{\tau}\right) |\left\langle D^{(2)}_{i},[[\textbf{U}]]\right\rangle| + \left(1 - \frac{\widehat{S}^D_p}{\tau}\right) |\left\langle D^{(1)}_{i},[[\textbf{U}]]\right\rangle| \right) c_i
\\
& - \left(1 + \frac{\widehat{S}^f_p}{\tau}\right)\sigma |[[\textbf{U}_i]]| > 0.
\end{aligned}
\ee
Observe that it is particularly simple to satisfy the WCD condition in this case as it is set of independent linear relations. The scheme parameter $c$ in (\ref{eq:fds2p_MHD}) is defined as
$
	c = \max (c_1, c_2). 
$

\subsubsection{Numerical experiments}

We set $ v = r \cos (\theta)$ and $w = r \sin (\theta)$
and consider the following class of Riemann intial data 
\be
	r(0, x) = \left\lbrace \begin{array}{cc} r_L, & x < 0.25, 
\\ 
0.6r_L, & x > 0.25,	
 \end{array} \right. \qquad \qquad
	\theta(0, x) = \left\lbrace \begin{array}{cc} \frac{3}{10} \pi, & x < 0.25,
 \\ \frac{13}{10} \pi, & x > 0.25,	 \end{array} \right. 
\ee
in which the state values of $r_L$ and $\alpha$ will be varied in our experiments.

\begin{figure}[htbp]
\centering
\begin{tabular}{cc}
\includegraphics[width=0.48\linewidth]{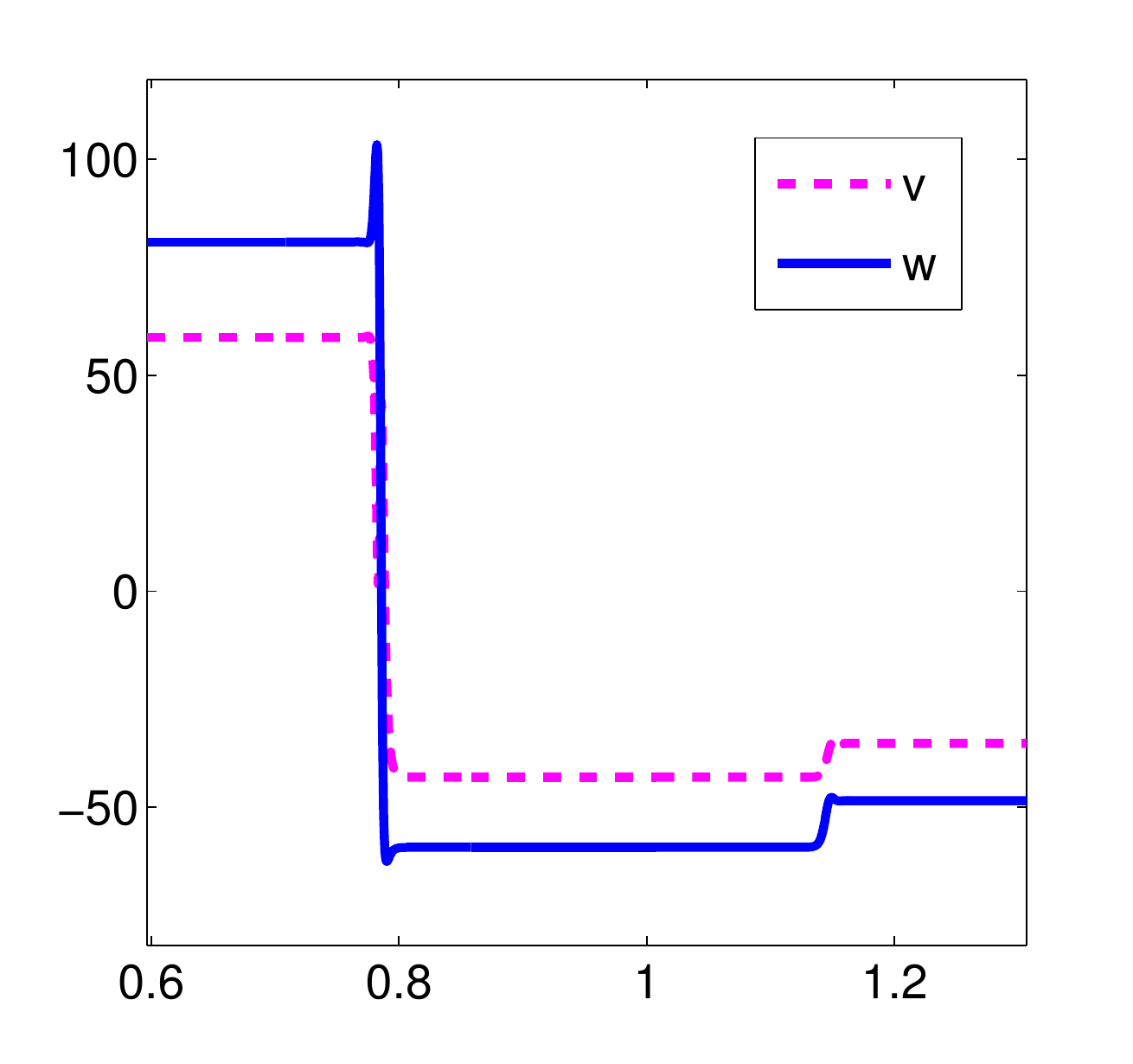} 
& \includegraphics[width=0.48\linewidth]{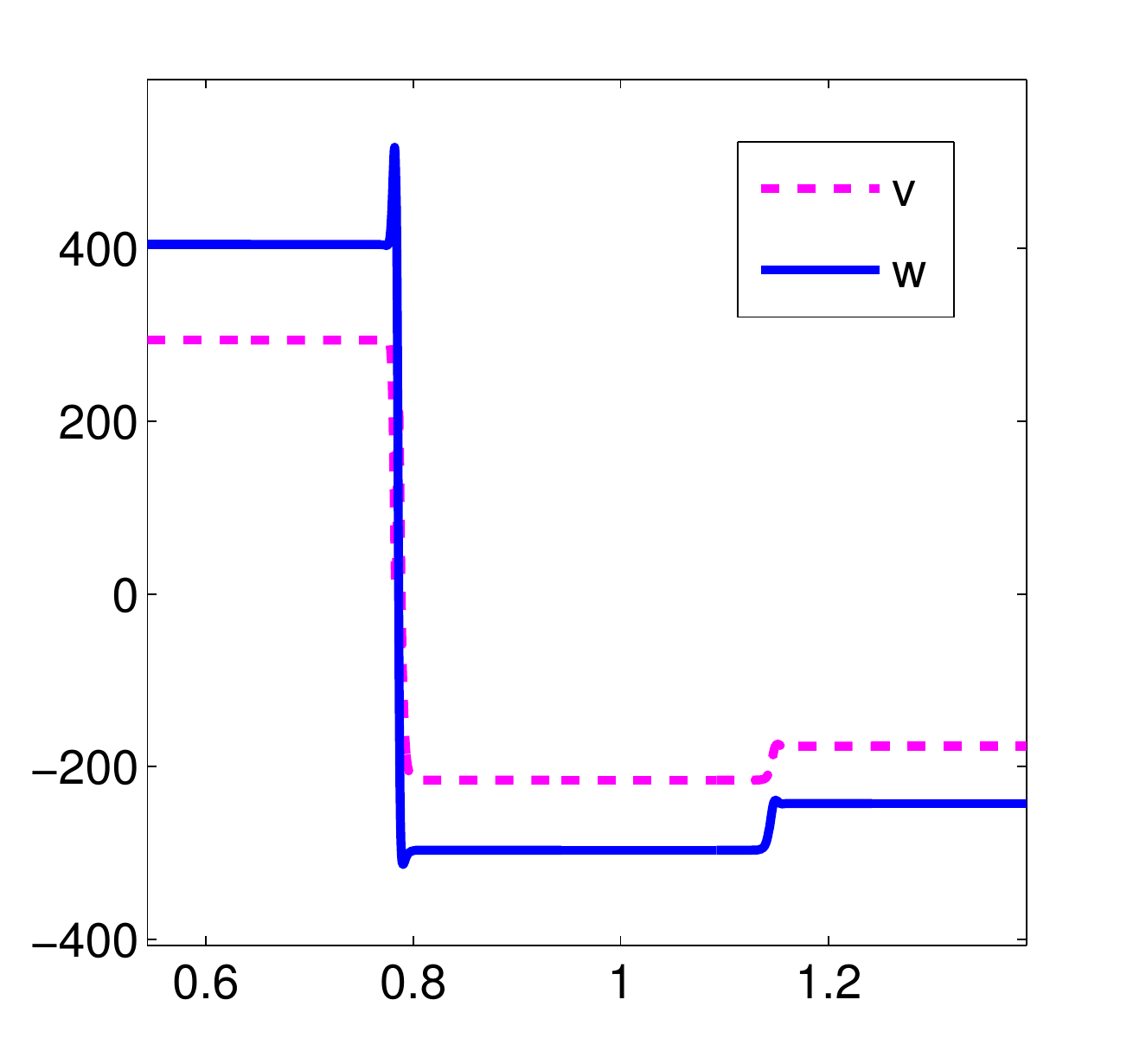}
\end{tabular}
\caption{Large shocks in $v$ and $w$-variable for the Hall MHD system using a fourth order WCD scheme with $4000$ mesh points.}
\label{fig:numMHD6}
\end{figure}

We consider the approximation of strong shocks by setting $r_L =100$ and $r_L = 500$. The approximations, performed with an fourth-order $WCD$ scheme on a mesh of $4000$ uniformly spaced points is shown in Figure~\ref{fig:numMHD6}. The results clearly show that the solution consists of very large amplitude ${\mathcal O}(10^3)$ nonclassical shocks in both the unknowns. The fourth-order WCD scheme is able to approximate these  very large nonclassical shocks quite well.

\subsubsection{Computing the kinetic relation}

\begin{figure}[htbp]
\centering
\begin{tabular}{cc}
\includegraphics[height=4.5cm, width=0.48\linewidth]{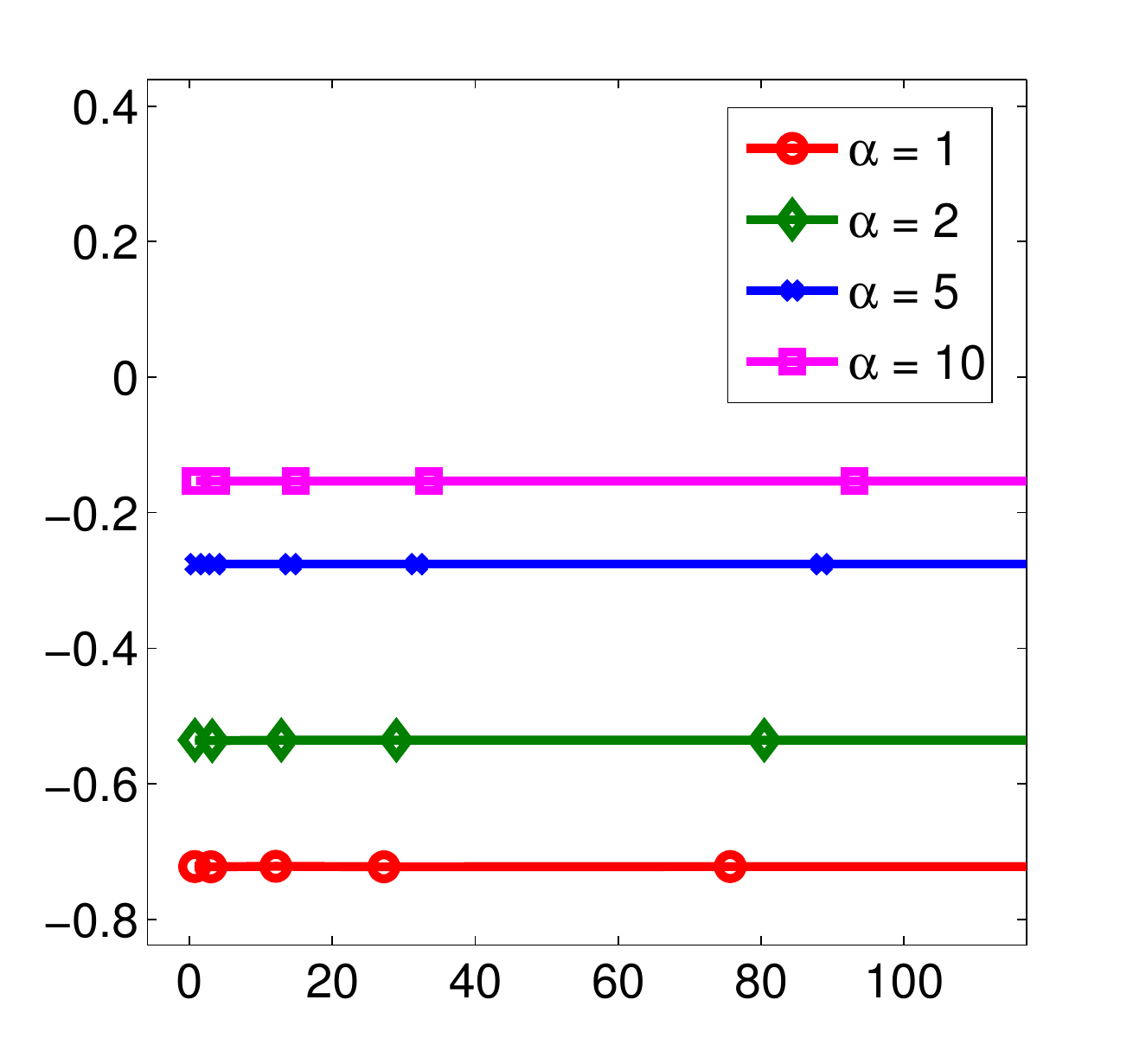} 
& \includegraphics[height=4.5cm, width=0.48\linewidth]{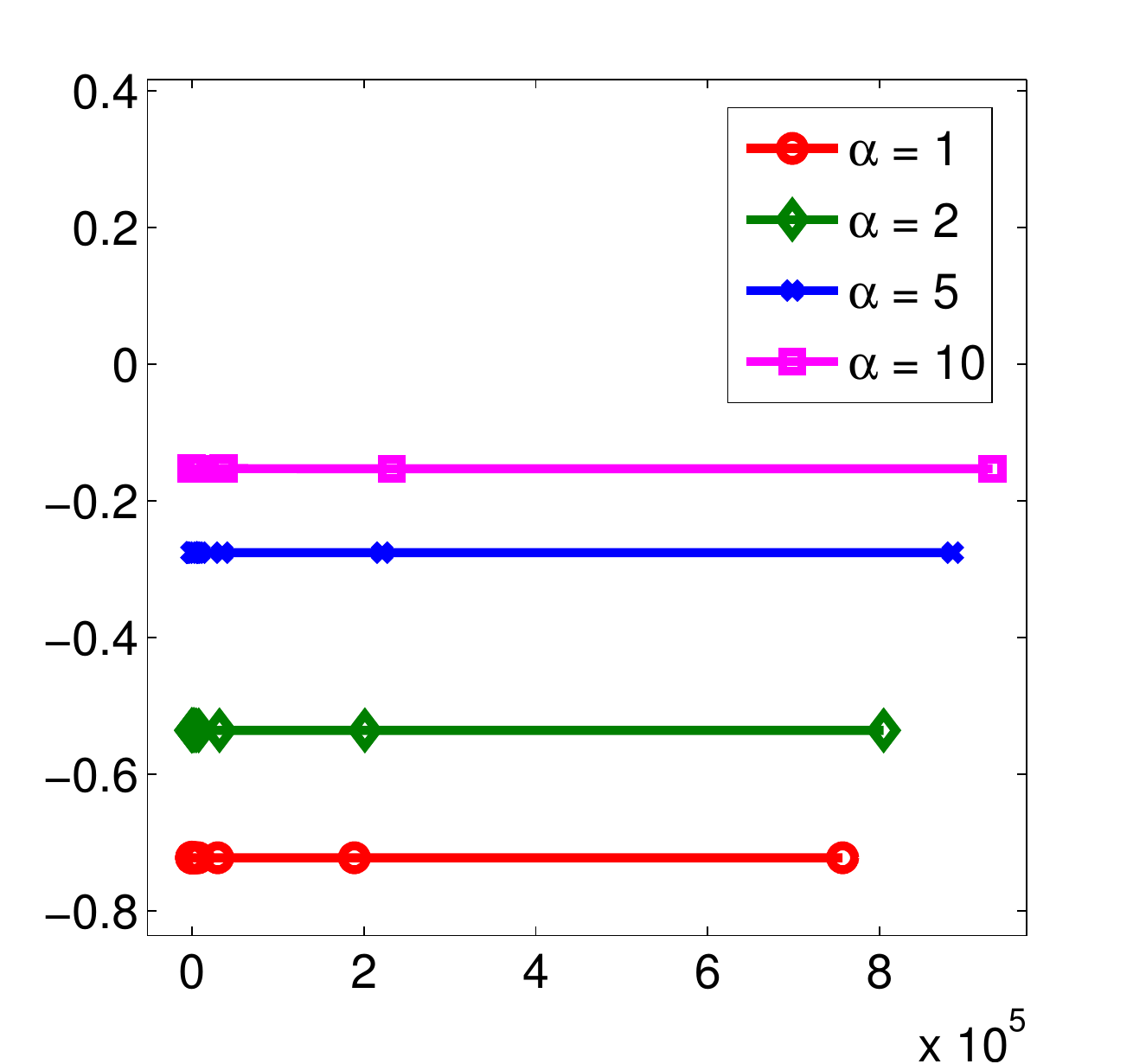}
\end{tabular}
\caption{Scaled entropy dissipation $\frac{\phi(s)}{s^2}$ vs. shock speed $s$
\newline 
 for the Hall MHD model (\ref{834}) with a fourth-order WCD scheme.} 
\label{fig:numMHD9}
\end{figure}

We now examine the kinetic relation 
$$
	\phi(s) = -s \frac{1}{2} [[v^2 + w^2]] + \frac{3}{4} [[(v^2 + w^2)^2]]
$$
at nonclassical shocks 
numerically for $\alpha = 1, 2, 10$ using the $4$-th order WCD scheme with $\tau = 0.1$ on a grid with $N=4000$ mesh points. The numerical kinetic relation is plotted as the 
scaled entropy dissipation versus the shock speed and is
 shown in Figure~\ref{fig:numMHD9}. These results suggest that the kinetic relation for the simplified MHD model 
with Hall effect has the {\sl quadratic expression} 
\be
\phi(s) = k_{\alpha} s^2,
\ee
for some constant $k_{\alpha}$ which depends upon the value of the Hall coefficient $\alpha$. Our results demonstrate the ability 
of the WCD schemes to compute nonclassical shocks to \eqref{834} with arbitrary strength.

%---------------------------------------------------------------------------------------------------------------

\subsection{Entropy stable WCD schemes}

The WCD schemes \eqref{eq:fds2p} constructed in the previous section may not satisfy a discrete version of the entropy inequality \eqref{eq:efdm1}, but can be modified as we now explain. 
Namely, there has been considerable progress in the last two decades regarding the construction of `entropy stable schemes' of arbitrary order. In this context, we refer to the pioneering contributions of Tadmor (1987, 2003)  who was able to characterize entropy stable, first-order, numerical fluxes for systems of conservation laws, as those which are more diffusive than an entropy conservative flux. Recall here that an entropy conservative flux for the finite difference (or finite volume) scheme \eqref{eq:fdm1} is a consistent, numerical flux $g^{\ast}_{i+1/2}$ such that the resulting scheme satisfies a 
{\sl discrete entropy identity}, of the form 
\be
\label{eq:efdm2}
\frac{d}{dt} U(u_i) + \frac{1}{\Dx}\left(G^{\ast}_{i+1/2}(t) - G^{\ast}_{i-1/2}(t) \right) = 0,
\ee
which is associated with an  entropy $U$ and some numerical entropy flux $G^{\ast}$. Tadmor (1987, 2003) showed the existence of such entropy conservative fluxes (consistent with the entropy flux $F$)
and provided a recipe for the construction of entropy stable schemes, which were 
determined by adding a diffusive part to the entropy conservative flux (in terms of the entropy variables $v := \partial_{u} U$)
and resulted in a discrete entropy inequality \eqref{eq:efdm1}. 
 
The next major step in the construction of entropy conservative fluxes was taken by LeFloch and Rohde (2000) and, later, LeFloch, Mercier and Rohde (2002). Therein, the authors were able to construct {\sl arbitrarily high-order} accurate entropy conservative fluxes for systems of conservation laws. Furthemore, fully discrete schemes were also designed in LeFloch, Mercier and Rohde (2002). In more recent years, the evaluation of entropy conservative fluxes in terms of explicit and simple to implement formulas has been proposed by Tadmor (2003), Ismail and Roe (2009), Fjordholm, Mishra, and Tadmor (2009, 2010), Kumar and Mishra (2012), and in references therein. Furthermore, the design of arbitrarily high-order and computationally efficient numerical diffusion operators was proposed by Fjordholm, Mishra and Tadmor (2012). This numerical diffusion operator was based on an ENO reconstruction of the scaled entropy variables that satisfied a subtle sign property, as later shown in Fjordholm, Mishra, and Tadmor (2013). An entropy stable spacetime discontinuous Galerkin finite element method that also utilizes entropy conservative fluxes was also introduced by Hiltebrand and Mishra (2013).

The use of entropy conservative fluxes in the context of computation of nonclassical shock waves was pioneered by Hayes and LeFloch (1996) and, later, LeFloch and Rohde (2000). We observe here that their technique can be readily adapted to the context of our WCD schemes. We illustrate the approach by considering scalar conservation laws, realized as a limit of vanishing viscosity and dispersion \eqref{eq:240} and discretized on a uniform mesh, by the following finite difference scheme:
\be
\label{eq:fds2pent}
\qquad 
\frac{du_i}{dt} + \frac{1}{\Dx} \sum\limits_{j=-p}^{j=p} \zeta_j g^{\ast}(u_i,u_{i+j}) 
 = \frac{c}{\Dx} \sum\limits_{j=-p}^{j=p} \beta_j u_{i+j}
 +  \frac{\delta c^2}{\Dx} \sum\limits_{j=-p}^{j=p} \gamma_j u_{i+j}.  
\ee
Here, the coefficients $\beta,\gamma$ were specified in \eqref{eq:fds2p} and the coefficients $\zeta$ are given by 
$$
\zeta_0 = 0, \quad \zeta_j := 2 \alpha_j, \qquad  j \neq 0,
$$
with $\alpha_j$ being the coefficients specified in \eqref{eq:fds2p}. Furthermore, the numerical flux $g^{\ast}$ is an entropy conservative flux given by
\be
\label{eq:entcon}
[[v]]_{i+1/2} g^{\ast}_{i+1/2} := [[vf - F]]_{i+1/2},
\ee
in which an entropy function $U$, an entropy flux function $F$, and the corresponding entropy variable $v$ have been fixed. Recall that the above formula determines a unique entropy conservative flux for scalar conservation laws. The resulting scheme is both $2p$-th order accurate (formally) as well as entropy stable, that is, satisfies a discrete entropy inequality of the form \eqref{eq:efdm1}. Furthermore, the entire analysis of imposing the WCD condition, as explained in the previous sections, can be easily modified to cover this situation.  Preliminary numerical experiments suggest very similar performance of the entropy stable scheme in comparison to the corresponding WCD scheme \eqref{eq:fds2p}. 

The construction of entropy stable WCD schemes can be extended to systems of conservation laws, like the MHD model \eqref{834}, by requiring that the entropy conservative flux satisfy
\be
\label{eq:entcon1}
\langle[[v]]_{i+1/2}, g^{\ast}_{i+1/2} \rangle:= [[\langle v,f \rangle - F]]_{i+1/2}
\ee
for the entropy variables $v$ and entropy flux $F$. Explicit formulas for the entropy conservative flux $g^{\ast}$ for \eqref{834} were derived in LeFloch and Mishra (2009).  

%===============================================================================

\section{Nonconservative hyperbolic systems}
\label{sec:5}

\subsection{Models from continuum physics}

We now turn our attention to nonlinear hyperbolic systems in nonconservative form, that is, \eqref{eq:111} when the matrix $A$ cannot be written in terms of the Jacobian of a flux. (In other words, $A(u) \neq D_u f(u)$ for any $f$.)
As mentioned in the introduction, many interesting systems in physics and engineering sciences takes this nonconservative form. 
The rigorously mathematical study of such systems was undertaken in LeFloch (1988, 1990), whose initial motivation came from 
the theory of two-phase flows and the dynamics of hypo-elastic materials, and has now been built upon the so-called DLM theory proposed by Dal Maso, LeFloch, and Murat (1990, 1995).  

As a first prototypical example, we consider the coupled Burgers' equation, proposed by Castro, Macias, and Pares (2001) and Berthon (2002):
\bel{eq:coupBurg}
\begin{split}
\partial_t w+ w\,\partial_x(w+ v) &= 0, \\
\partial_t v+ v\,\partial_x(w +v) &= 0.
\end{split}
\ee
The system has the nonconservative form \eqref{eq:111} with
$$
u= \begin{pmatrix}
w\\
v
\end{pmatrix},
\qquad
A(u) = \begin{pmatrix}
w & w\\
v & v
\end{pmatrix}.
$$
By adding the two components of this system, we obtain Burgers equation for the dependent $\omega := w+v$:
$$
\partial_t \omega + {1 \over 2} \partial_x \omega^2  = 0.
$$

A second prototypical example for the class of nonconservative hyperbolic systems is provided by the system of four equations governing a (one-dimensional, say) flow of two superposed immiscible shallow layers of fluids: 
\be\label{eq:tlsw}
\begin{split}
(h_1)_t + (h_1u_1)_x &= 0, \\
(h_2)_t + (h_2u_2)_x &= 0, \\
(h_1u_1)_t + \left(\frac{1}{2}gh_1^2 + h_1u_1^2\right) &= -gh_1(b+h_2)_x, \\
(h_2u_2)_t + \left(\frac{1}{2}gh_2^2 + h_2u_2^2\right) &= -gh_2(b+rh_1)_x.
\end{split}
\ee
Here,  $u_j=u_{j}(t,x)$ and $h=h_{j}(t,x)$ (for $j=1,2$) represent the depth-averaged velocity and thickness of the $j$-th layer, respectively, while $g$ denotes the acceleration due to gravity and $b = b(x)$ is the bottom topography. In these equations, the index $1$ and $2$ refer to the upper- and lower-layers. Each layer is assumed to have a constant density $\rho _{j}$ (with $\rho _{1}<\rho _{2}$), and $r = \rho_1/\rho_2$ represents the density ratio.

Other examples for the class of nonconservative hyperbolic systems include the equations governing multiphase flow (studied in Saurel and Abgrall, 1999) and a version of the system of gas dynamics in Lagrangian coordinates (LeFloch 1988, Karni 1992, Abgrall and Karni 2010).

%------------------------------------------------------------------------------------------------------------------

\subsection{Mathematical framework}

Solutions to nonlinear hyperbolic systems in the nonconservative form \eqref{eq:111} can contain discontinuities such as shock waves and, therefore, an essential mathematical difficulty is to define a suitable notion of weak solutions to such systems. Namely, the nonconservative product $A(u)u_x$ in \eqref{eq:111} cannot be defined in the distributional sense, by integrating by parts, since this term does not have a divergence form. The theory introduced by Dal Maso, LeFloch, and Murat (1990, 1995) allows them to define the nonconservative product ${A}(u)u_x$  as a {\sl bounded measure} for all functions $u$ with bounded variation, provided a {\sl family of Lipschitz continuous paths} $\Phi: [0,1] \times \Ucal \times \Ucal \to \Ucal$ is prescribed, which must satisfy certain regularity and  compatibility conditions, in particular
\be \label{cond1}
\Phi(0;u_l, u_r) = u_l, \qquad  \Phi(1;u_l, u_r) = u_r, \qquad \Phi(s; u, u) = u. 
\ee 
We refer to Dal Maso, LeFloch, and Murat (1990, 1995)  for a full presentation of the theory. 

Once the nonconservative product has been defined, one may define the weak solutions of \eqref{eq:111}.  According to this theory, across a discontinuity a weak solution has to satisfy the generalized Rankine-Hugoniot condition
\be \label{gen_R-H}
\sigma [[u]] = \int_0^1 {A}(\Phi(s;u_-,u_+)) \partial_s\Phi(s;u_-,u_+)\, ds,
\ee
where $\sigma$ is the speed of propagation of the discontinuity, $u_-$ and $u_+$ are the left- and right-hand limits of the solution at the discontinuity, and $[[u]] = u_+ - u_-$. Notice that, if $A(u)$ is the Jacobian matrix of some function $u$, then \eqref{gen_R-H} reduces to the standard Rankine-Hugoniot conditions for the conservation law \eqref{eq:109}, regardless of the chosen family of paths.

Analogous to the theory of conservation laws, many interesting {\sl nonconservative} hyperbolic systems arising in continuum physics are equipped with an entropy formulation (LeFloch, 1988).  More specifically, \eqref{eq:111} is equipped with an entropy pair $(\eta, q)$, i.e.~a convex function $\eta: \Ucal \to \mathbb{R}$ and a function $q: \Ucal \to \mathbb{R}$ such that $\nabla q(u)^\top = v^\top A(u)$, where $v := \nabla \eta(u)$ is the so-called entropy variable. Then entropy solutions of \eqref{eq:111} satisfy the following entropy inequality (in the sense of distributions)
\be\label{eq:entr): neq}
\eta(u)_t + q(u)_x \leq 0, 
\ee
while {\sl smooth} solutions to \eqref{eq:111} satisfy the entropy equality:
\be
\label{eq:entreq}
\eta(u)_t + q(u)_x = 0.
\ee

%---------------------------------------------------------------------------------------------------------------

\subsection{Small-scale dependent shock waves}

Following an idea by LeFloch (1990), we now explain how the family of paths is derived in applications, from an augmented model associated with a nonlinear hyperbolic system. Observe 
that the concept of entropy weak solutions, as outlined above, depends on the chosen family of paths. Different families of paths lead to different jump conditions, hence different weak solutions. A priori, the choice of paths is arbitrary. Thus, the crucial question is how to choose the ``correct'' family of paths so as to recover the physically relevant solutions.

In practice, any hyperbolic system like \eqref{eq:111} is obtained as the limit of a regularized problem when the high-order terms (corresponding to small-scale effects) tend to $0$. For instance, it may be the vanishing-viscosity limit of a family of hyperbolic-parabolic problems (as in the conservative case):
\be
u_t^\eps + A(u^\eps) u_x^\eps = \eps \, (B(u^\eps)u^\eps_x)_x, 
\label{NC-p}
\ee
where the right-hand side is elliptic in nature on account of the viscosity matrix $B$. Then, the correct jump conditions (corresponding to the physically relevant solutions) should be consistent with the \emph{viscous profiles}, that is, with traveling wave solutions $u^\eps(t,x) = V\left( \frac{x - \sigma t}{\eps}\right)$ of \eqref{NC-p} 
satisfying $\lim_{\xi \to \pm\infty} V(\xi ) = u_\pm, \lim_{\xi \to \pm \infty} V'(\xi) = 0.$
A single-shock solution
\bel{limitdisc}
u(t,x) =
\begin{cases}
u_-, \qquad & x < \sigma t, \\
u_+, &  x > \sigma t, 
\end{cases}
\ee
is considered `admissible' if $u = \lim_{\epsilon \rightarrow 0} u^\epsilon$ (almost everywhere).

It is easily checked that the viscous profile $V$ has to satisfy the system of ordinary differential equations  
$$
-\sigma V' + A(V)V' = (B(V)V')'.
$$
By integrating these equations in $\xi \in \mathbb{R}$, we obtain the jump condition (LeFloch, 1990) 
\be
 \label{vpjump}
\sigma [[u]] = \int_{-\infty}^\infty A(V(\xi)) \, V' (\xi) \, d\xi.
\ee
By comparing this jump condition with \eqref{gen_R-H}, we conclude that the correct choice for the path connecting the states $u_-$ and $u_+$ should be (after re-parameterization) the trajectory of the viscous profile $V$.  

From the above discussion, we see that different choices of viscous term $B$ in \eqref{NC-p} may lead to different viscous profiles and, consequently, different jump conditions. This dependence upon the jump conditions (and thus of the definition of weak solutions) on the explicit form of the neglected small-scale effects has profound implications on the design of efficient numerical methods, as we already explained it in the case of non-convex conservative problems in Section~\ref{sec:3}.

%--------------------------------------------------------------------------------------------------------

\subsection{Error in formally path consistent schemes}

Unlike for systems of conservation laws, where consistency of finite difference (or finite volume) schemes was established by Lax and Wendroff (1960) by requiring that the schemes has a discrete conservative form,  no such requirement is known about discrete schemes for nonconservative systems. This {\sl lack of convergence} of nonconservative schemes was discovered by Hou and LeFloch (1994) and turned out to difficult to observe, as the error may be quite small in certain applications; error estimates were derived therein, which involve
 the strength of shocks and the order of the schemes between the numerical and the exact solutions. 
Hou and LeFloch's theory was motivated by a work by Karni (1992), who advocated, and demonstrated a definite advantage of, using nonconservative schemes in certain application to fluid dynamics. 

A large literature is now available on the problems and numerical methods presented in this section and we will not be able to present here 
a fully exhaustive review. For further reading we thus refer to Audebert and Coquel (2006), 
Berthon and Coquel (1999, 2002, 2006, 2007), and Chalons and Coquel (2007), and Berthon, Coquel, and LeFloch (2002, 2012). 

More recently, Pares (2006) introduced a notion of \emph{path consistent schemes} for the nonconservative systems \eqref{eq:111}. Given a uniform grid as in the previous section, numerical schemes for nonconservative systems can be written in the following \emph{fluctuation} form:
\be\label{eq:ncfds}
\frac{d}{dt}u_i + \frac{1}{\Dx}\left(D^+_{i-1/2} + D^-_{i+1/2}\right) = 0, 
\ee
where 
$D_{i+1/2}^{\pm}(t) = D^\pm\left(u_i(t), u_{i+1}(t) \right)$
and $D^\pm: \Ucal \times \Ucal \mapsto \Ucal$ are Lipschitz continuous functions satisfying  
\be
D^\pm(u, u) = 0.\label{eq:pc1}
\ee
As we just explained, weak solutions to \eqref{eq:111} require the specification of a DLM family of paths. In Pares (2006), 
the path is explicitly introduced into the scheme \eqref{eq:ncfds} by imposing the (formal) \emph{path consistency} condition  
 \be
D^-(u_l, u_r) + D^+(u_l, u_r) = \int_0^1 A(\Phi(s; u_l, u_r) )
\partial_s \Phi(s; u_l, u_r) ds.
\label{eq:pc2}
\ee
\label{eq:pc}
Of course, a family of paths must be specified in advance in \eqref{eq:pc2} before writing a path consistent scheme. 

Assuming that a suitable path is selected (for instance from an augmented model and by deriving its viscous profiles), 
it is natural to investigate whether the approximate solutions to \eqref{eq:111} by the path consistent scheme \eqref{eq:ncfds} converges to the correct (physically relevant) solution of the nonconservative system \eqref{eq:111}. 
Unfortunately, the answer to this fundamental question is negative in most cases, as follows from 
Hou and LeFloch (1994), which was revisited by Castro, LeFloch, Munoz-Ruiz, and Pares (2008) and, more recently, Abgrall and Karni (2010).

\begin{figure}[t]
\subfigure[Top: $u$. Bottom: $v$]{\includegraphics[width=0.48\linewidth]{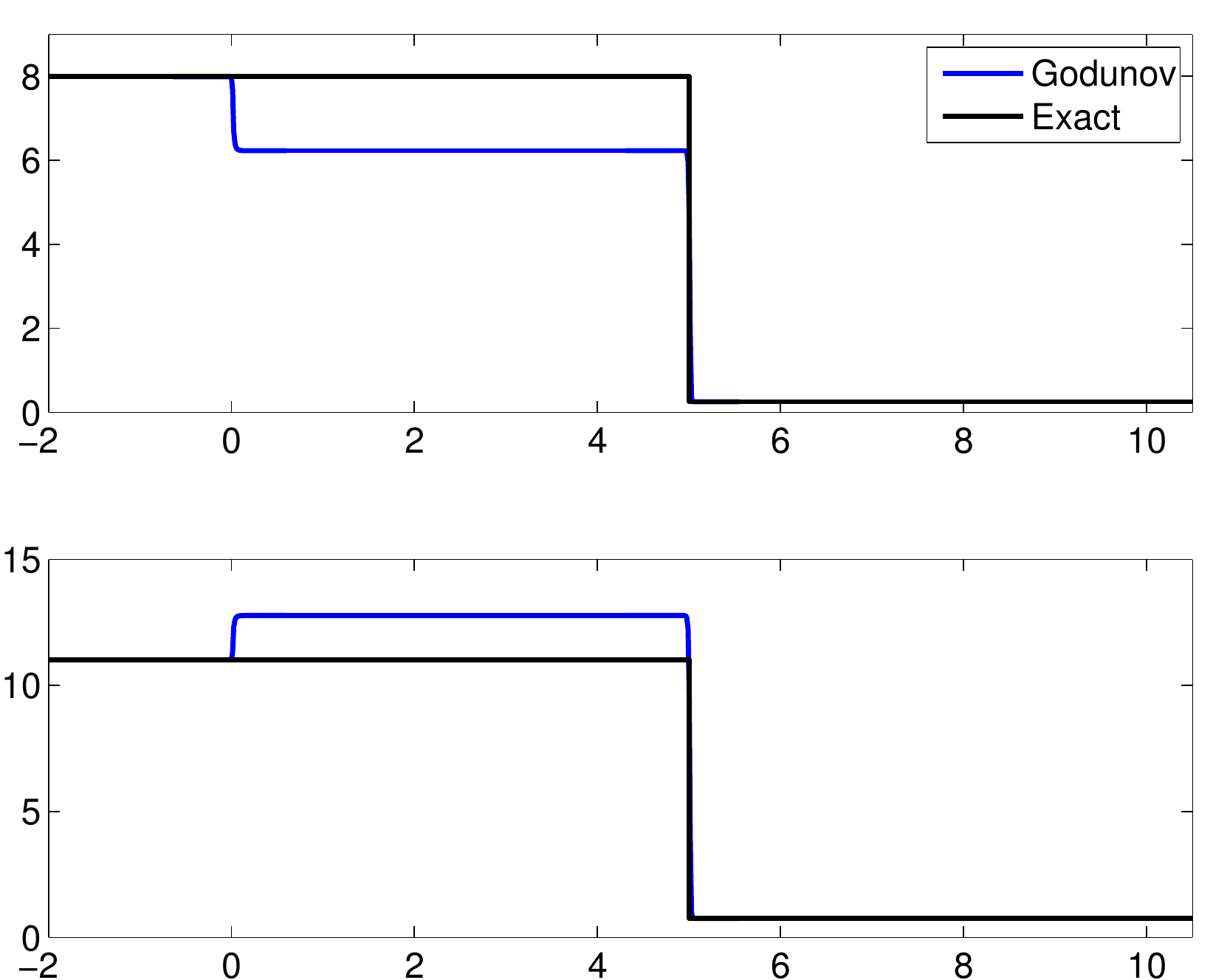}} 
\caption{Godunov method of Munoz and Pares (2007) for the coupled Burgers system \eqref{eq:coupBurg} with CFL= 0.4 and 1500 grid points. Comparison with the exact solution computed from the viscous regularization \eqref{NC-p-coupBurg}.}
\label{fig:godex}
\end{figure}

Here, we illustrate this deficiency of path consistent schemes by considering a very simple nonconservative system -- the coupled Burgers system \eqref{eq:coupBurg} of Berthon (2002). A Godunov-type scheme was derived in Munoz and Pares (2007)  and was shown to be (formally) consistent with the path computed from viscous profiles to the corresponding parabolic regularization \eqref{NC-p-coupBurg}. A numerical example (with further details are provided in the following subsection) is shown in Figure~\ref{fig:godex}. The results show that although the Godunov-type path-consistent scheme converges as the mesh is refined, it \emph{does not converge to the physically relevant} (correct) solution, which is computed explicitly from the parabolic regularization. 

As in the case of all small-scale dependent shocks to (conservative or nonconservative)) hyperbolic systems, an explanation for this lack of convergence of path-consistent schemes lies in the \emph{equivalent equation} of the scheme \eqref{eq:ncfds}, say, 
\be
u_t^{\Delta x} + A(u^{\Delta x})u_x^{\Delta x} = \Delta x
\big( \widetilde{B}(u^{\Delta x})u^{\Delta x}_x \big)_x + \mathcal{H}.
\label{NC-p-disc}
\ee
Here, $\mathcal{H}$ includes the higher-order terms that arise from a formal Taylor expansion in the scheme \eqref{eq:ncfds} and $\widetilde{B}$ is the (implicit) numerical viscosity. Assuming that the high-order terms are relatively \emph{small} (which is valid for shocks with small amplitude), we can expect that jump conditions associated with the numerical solutions to be, at best, consistent with the viscous profiles of the regularized equation
\be
u_t^{\Delta x} + A(u^{\Delta x})u_x^{\Delta x} = \Delta x
\big( \widetilde{B}(u^{\Delta x})u^{\Delta x}_x \big)_x.
 \label{NC-p-disc-tr}
\ee
But, in general,  $B \neq \widetilde{B}$. As was discussed before, the solutions to the nonconservative system \eqref{eq:111} depend \emph{explicitly} on the underlying viscosity operator. Therefore, the numerical solutions generated by the scheme \eqref{eq:ncfds} may not converge to the physically relevant solutions of \eqref{eq:111}. Thus, the (implicit) numerical viscosity that is added by any finite difference scheme (as observed first in Hou and LeFloch, 1994) is responsible for the observed lack of convergence to the physically relevant solutions.

%------------------------------------------------------------------------------------------------------------------

\subsection{Schemes with controlled diffusion}

As in the case of conservative hyperbolic systems with small-scale dependent shock waves, the equivalent equation provides the appropriate tool to design finite difference schemes that can approximate small-scale dependent shocks to nonconservative systems.  The main idea is to design a finite difference scheme, such that the numerical viscosity (the leading term in its equivalent equation) matches that of the underlying physical viscosity in \eqref{NC-p}. In addition, we would like this scheme to be entropy stable, that is, to satisfy a discrete version of the entropy inequality. One may proceed by requiring the (formal) path consistent condition introduced in Pares (2006), while proving the corrections suggested by Hou and LeFloch (1995) and Castro, LeFloch, Munoz-Ruiz, and Pares (2008).  

Following Castro, Fjordholm, Mishra, and Pares (2013), let us present here a class of schemes in the fluctuation form \eqref{eq:ncfds} with fluctuations satisfying 
the (entropy) consistency condition 
\be
\label{eq:ncECcond}
v_l ^\top D^-(u_l,u_r) +  v_r^\top D^+(u_l,u_r) = q(u_r) - q(u_l),       \qquad u_l, u_r \in \Ucal,
\ee
where $v = \partial_{u} \eta$ is the entropy variable. Furthermore, the above specified fluctuations are also required to satisfy the path-consistency condition \eqref{eq:pc2} for any choice of path. It is shown in Castro et. al (2013) that the resulting scheme \eqref{eq:ncfds} is \emph{entropy conservative}, that is, satisfies a discrete version of the entropy identity \eqref{eq:entreq}. Furthermore, the existence of such fluctuations was also proved therein. In particular, the condition \eqref{eq:ncECcond} was shown to be a natural extension of the definition of Tadmor's entropy conservative flux \eqref{eq:entcon1} to a nonconservative system. 
 
Following the construction of the entropy stable scheme \eqref{eq:fds2pent}, we need to add numerical viscosity to stabilize the entropy conservative path consistent scheme. Castro et al (2013) proposed to add the following viscosity operator to their fluctuations:
 \bea 
\widetilde{D}^+(u_i, u_{i+1}) &= D^+(u_i, u_{i+1}) +\frac{\eps}{\Delta x} \widehat{B}(v_{i+1} - v_i), \label{eq:DpmNV1}\\
\widetilde{D}^-(u_i, u_{i+1}) & = D^- (u_i, u_{i+1}) - \frac{\eps}{\Delta x} \widehat{B}(v_{i+1} - v_i)\label{eq:DpmNV2}
\eea
Here, $D^{\pm}$ are the entropy conservative fluctuations defined in \eqref{eq:ncECcond} and the numerical viscosity $\widehat{B}$ is defined as
\be
\label{eq:nv1}
\widehat{B} := B u_v,
\ee
with $B$ being the physical viscosity in \eqref{NC-p}. The corresponding scheme \eqref{eq:ncfds} was then shown to be 
\begin{itemize}
\item formally path-consistent, 

\item entropy stable, i.e, it satisfied a discrete version of the entropy inequality, and 

\item the leading-order term of the equivalent equation of this scheme matched the underlying parabolic equation \eqref{NC-p}.
\end{itemize}
Thus, this scheme is the correct extension of finite difference schemes with controlled dissipation which approximates small-scale dependent shocks to nonconservative hyperbolic systems. 

%-------------------------------------------------------------------------------------------------------------------

\subsection{Numerical experiments}

We present a few numerical experiments from the recent article (Castro et al 2013), which serve to illustrate the relevance of entropy stable schemes with controlled dissipation for approximating nonconservative hyperbolic systems. 

\subsubsection{Coupled Burgers system}

First, we consider the system \eqref{eq:coupBurg} and note that this system is equipped with the entropy-entropy flux pair:
$$
\eta(u) = \frac{\omega^2}{2}, \qquad  q(u) = \frac{\omega^3}{3}, \quad \omega = v + w.
$$
Recall that Berthon (2002) computed the exact viscous profiles of the regularized system, i.e. 
\be
\begin{aligned}
\partial_t w + w\, \partial_x (w+v) = \varepsilon_1 \partial^2_{xx}(w+v),\\
\partial_t v + v\,\partial_x (w+v) = \varepsilon_2 \partial^2_{xx}(w+v).\\
\end{aligned}
\label{NC-p-coupBurg}
\ee
In the limit $\epsilon_1,\epsilon_2 \rightarrow 0$, this gives the correct (physically relevant) entropy solutions to the Riemann problem for the coupled Burgers equation. In the rest of this section, we choose $\epsilon_1 = \epsilon_2 = \epsilon$.

In Munoz and Pares (2007), the authors devised a Godunov path-conservative method \eqref{eq:pc} with
\begin{align*}
D^-_{i+1/2} = & \int_0^1 A(\Phi(s; u_i^n, u_{i+1/2}^{n,-}))\partial_s \Phi(s; u_i^n, u_{i+1/2}^{n,-})\,ds,
 \\
D^+_{i+1/2} = & \int_0^1 A(\Phi(s;u_{i+1/2}^{n,+}, u_{i+1}^{n}))\partial_s \Phi(s; u_{i+1/2}^{n,+}, u_{i+1}^{n})\,ds,
\end{align*}
where $u_{i+1/2}^{n, \pm}$ are the limits to the left- and to the right-hand sides at $x = 0$ of the Riemann solution with initial data $(u_i^n, u_{i+1}^n)$.

To test the performance of the Godunov scheme, we approximate the Riemann problem for \eqref{eq:coupBurg} with initial data $u_l = [7.99, 11.01]^\top$, $u_r = [0.25, 0.75]^\top$ and we compare the exact solution with the numerical one provided by  the Godunov method in the interval $[-2,10.5]$ with 1500 points and CFL = 0.4. As shown in Figure \ref{fig:godex}, the location of the discontinuities is correctly approximated, whereas the intermediate states approximated by the Godunov scheme are {\sl incorrect.}
 The error in these
intermediate states does vanish as $\Delta x$ tends to 0.  Even though the Godunov method takes into account the exact expression of the viscous profiles (paths) and  the exact solutions of the Riemann problems, the numerical solutions
provided by the method {\sl do not converge} to the expected weak solutions,  due to the numerical viscosity added in the projection step.

\begin{figure}
\includegraphics[width=0.5\linewidth]{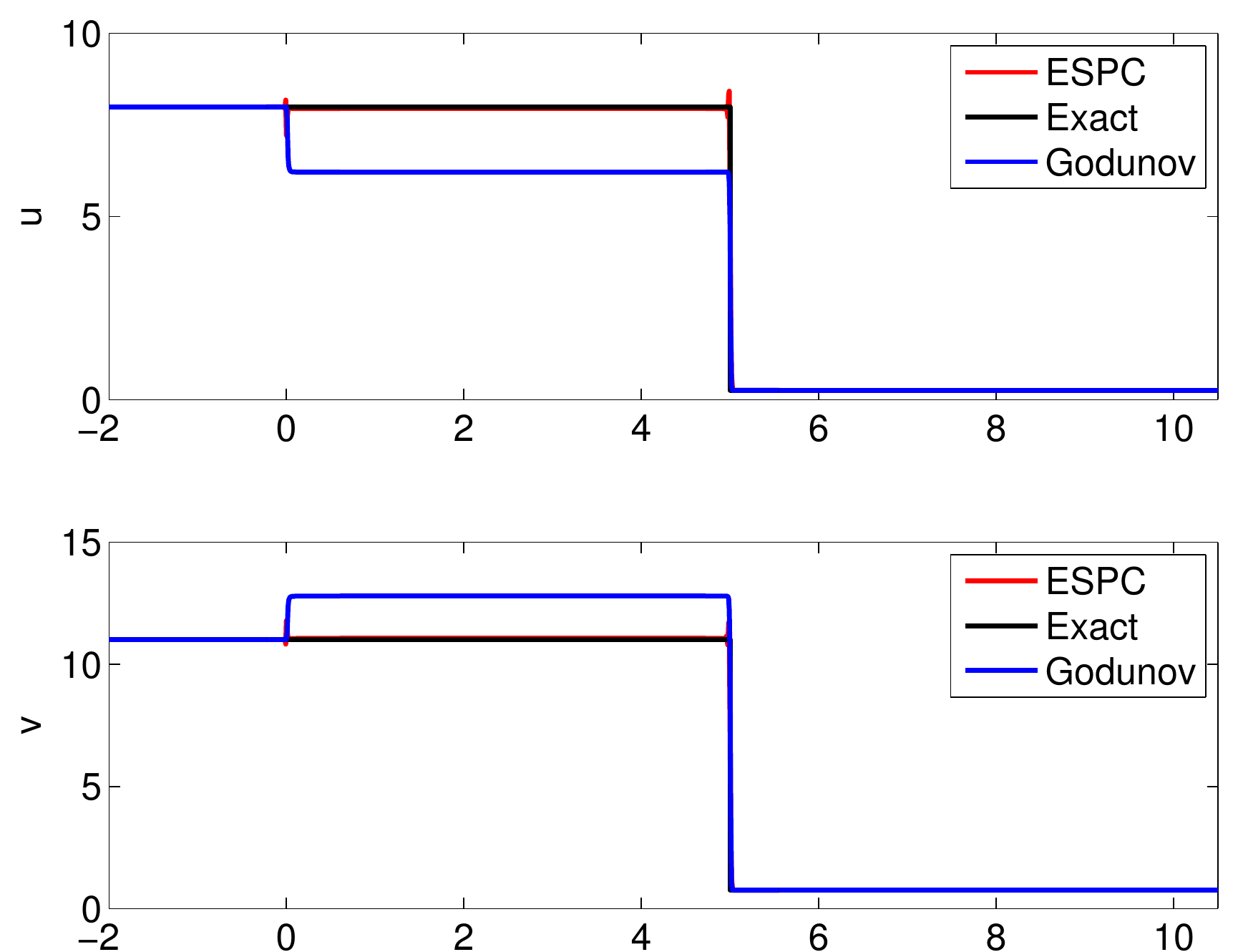}
\caption{Comparison of the ESPC and Godunov schemes for the coupled Burgers system \eqref{eq:coupBurg} with the exact solution.}
\label{fig:cp1}
\end{figure}

Following the procedure outlined in the construction of the entropy stable scheme above, the following entropy stable fluctuations were derived in Castro et al (2013) for the coupled Burgers equation:
\be
\label{eq:cpd}
\begin{aligned}
\widetilde{D}_{i+1/2}^- &= \frac{1}{6} [[\omega]]_{i+1/2} \begin{pmatrix}2w_i + w_{i+1} \\ 2v_i + v_{i+1}\end{pmatrix} - \frac{\epsilon}{\Dx}\begin{pmatrix}[[\omega]]_{i+1/2}\\ [[w]]_{i+1/2} \end{pmatrix}, \\
\widetilde{D}_{i+1/2}^+ &= \frac{1}{6} [[\omega]]_{i+1/2} \begin{pmatrix}2w_i + w_{i+1} \\ 2v_i + v_{i+1}\end{pmatrix} + \frac{\epsilon}{\Dx}\begin{pmatrix}[[\omega]]_{i+1/2}\\ [[w]]_{i+1/2} \end{pmatrix}.
\end{aligned}
\ee
In order to validate these numerical schemes, we consider again the Riemann problem with initial data $u_l = [7.99, 11.01]^\top$, $u_r = [0.25, 0.75]^\top$ 
and compare the exact solution with the numerical one provided by the ESPC scheme (the relevant scheme with controlled dissipation) in the interval $[-2,10.5]$ with 1500 grid points. The results are shown in Figure \ref{fig:cp1}. 
\begin{figure}
\includegraphics[width=0.45\linewidth]{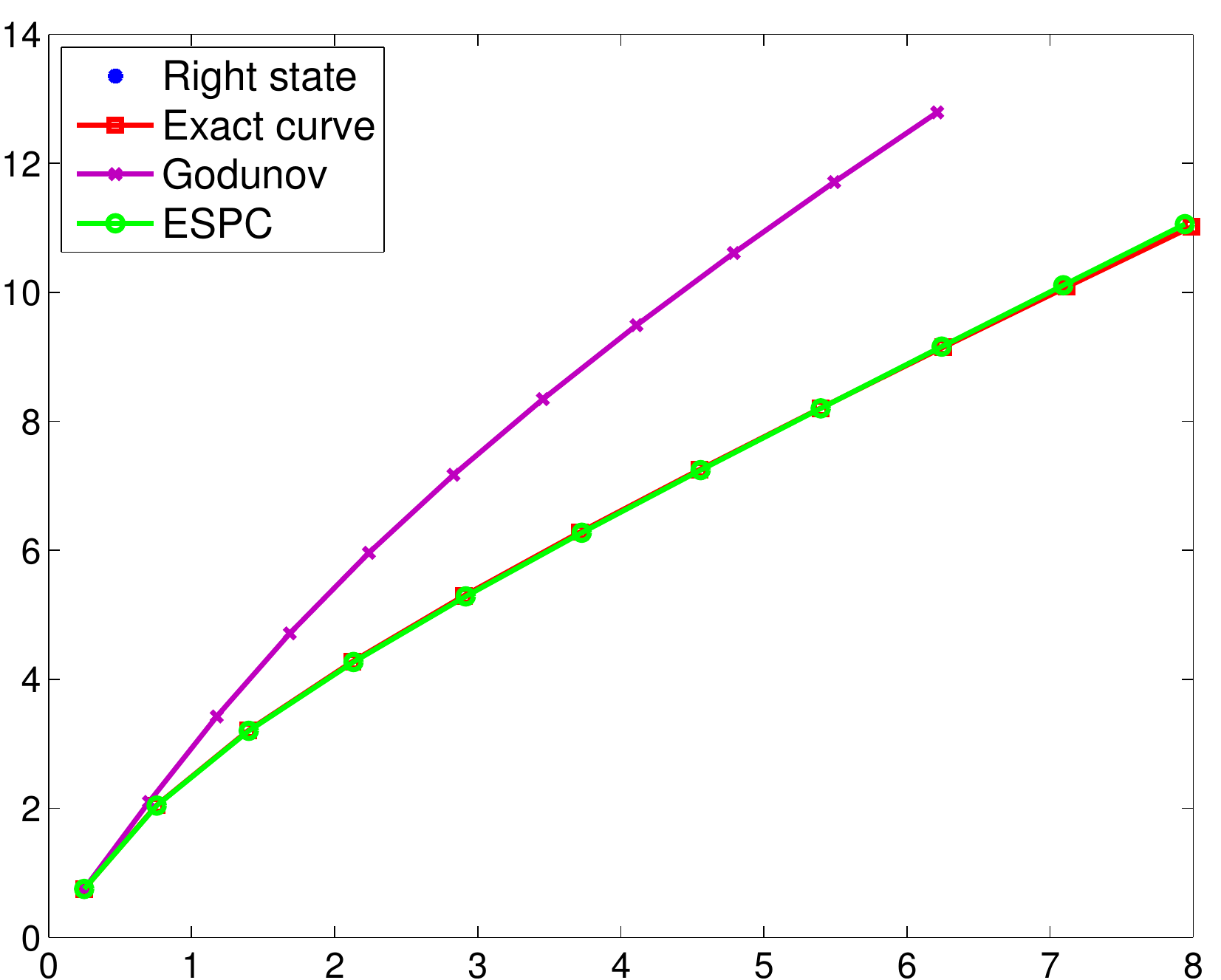}
\caption{Numerical Hugoniot locus for the coupled Burgers equation \eqref{eq:coupBurg} generated by the ESPC and Godunov schemes, compared with the exact Hugoniot locus.}
\label{fig:cp2}
\end{figure}
In order to compare the ESPC scheme with the Godunov scheme, we computed the numerical Hugoniot locus by approximating a family of Riemann problem whose initial data are given by $u_r = [0.75, 0.25]^\top$  and a series of left-hand states belonging to the exact shock curve.  The Riemann problem is solved in the interval $[-2, 10]$ and the corresponding left-hand state (at the shock) is used to compute the numerical Hugoniot locus. The results are presented in Figure \ref{fig:cp2} and show that the Godunov scheme does a poor job of approximating the exact solution. The numerical Hugoniot locus for this scheme starts diverging even for shocks with small amplitude. On the other hand, the ESPC schemes approximates the correct weak solution. The numerical Hugoniot locus coincides with the exact locus for a large range of shock strengths. Only for very strong shocks does the Hugoniot locus show a slight deviation. This is to be expected as the high-order terms in the equivalent equation \eqref{NC-p-disc} become larger with increasing shock strength and may lead to deviations in the computed solution. However, the gain in accuracy with the ESPC scheme over the Godunov scheme is considerable.

\subsubsection{Two-layer shallow water system}

Next, let us consider the two-layer shallow water system \eqref{eq:tlsw}. It is widely accepted that the correct regularization mechanism for this system is provided by the \emph{eddy viscosity} resulting in the following hyperbolic-parabolic system:
\be\label{eq:tlswv}
\begin{split}
(h_1)_t + (h_1u_1)_x &= 0, \\
(h_2)_t + (h_2u_2)_x &= 0, \\
(h_1u_1)_t + \left(\frac{1}{2}gh_1^2 + h_1u_1^2\right) &= -gh_1(b+h_2)_x + \nu(h_1(u_1)_x)_x,\\
(h_2u_2)_t + \left(\frac{1}{2}gh_2^2 + h_2u_2^2\right) &= -gh_2(b+rh_1)_x + \nu(h_2(u_2)_x)_x.
\end{split}
\ee
Here, $\nu \ll 1$ is the coefficient of eddy viscosity. An entropy-entropy flux pair for the two-layer shallow water system is given by
\begin{subequations}\label{eq:entroptlsw}
\begin{align} 
\eta &= \sum_{j=1}^2 \rho_j \left(h_j \frac{u_j^2}{2}+g \frac{h_j^2}{2} + gh_j b \right) + g\rho_1 h_1 h_2,
 \label{eq:entroptlsw1} \\
q &= \sum_{j=1}^2 \rho_j\left(h_j\frac{u_j^2}{2}+gh_j^2 + g h_j b \right)u_j + \rho_1 g h_1 h_2 (u_1+u_2). 
\label{eq:entroptlsw2} 
\end{align}
\end{subequations}
The corresponding entropy variables read 
$$
v =  \left( \begin{array}{c} \rho_1 \left(  -\frac{1}{2} u_1^2 + g(h_1 + h_2 + b) \right)  \\ \rho_1 u_1  \\\rho_2 \left(  -\frac{1}{2} u_2^2 + g( h_2 + b) \right)  + \rho_1gh_1 \\
\rho_2 u_2 \\ \rho_1 g h_1 + \rho_2 g h_2  \end{array} \right).
$$
An entropy stable scheme with controlled dissipation for this system was proposed in Castro et al. (2013) to which the reader is refered  for an explicit expression. 

\begin{figure} 
\subfigure[$h_2$ and $h_1+h_2$]{\includegraphics[width=0.45\linewidth]{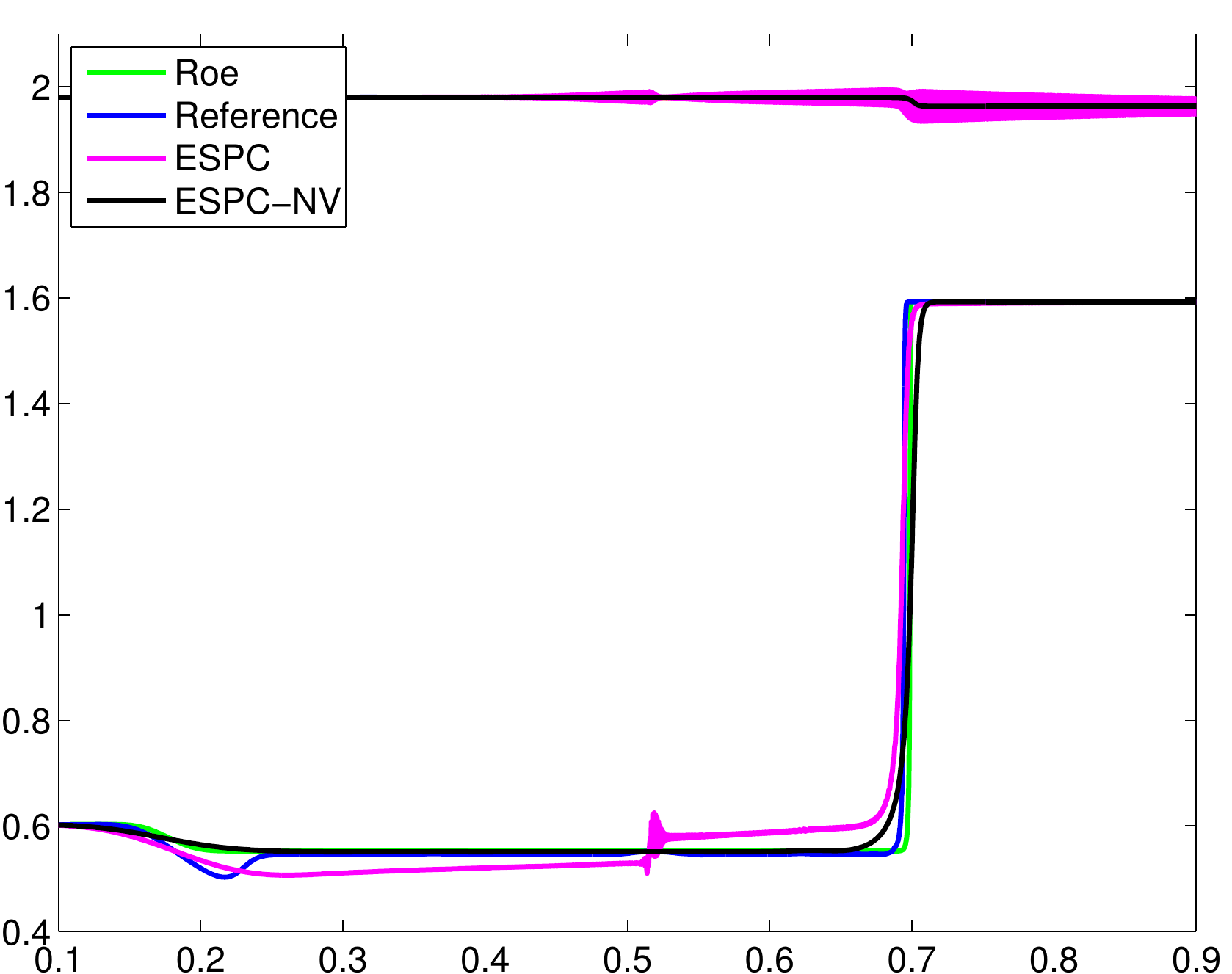}}
\subfigure[Closeup]{\includegraphics[width=0.45\linewidth]{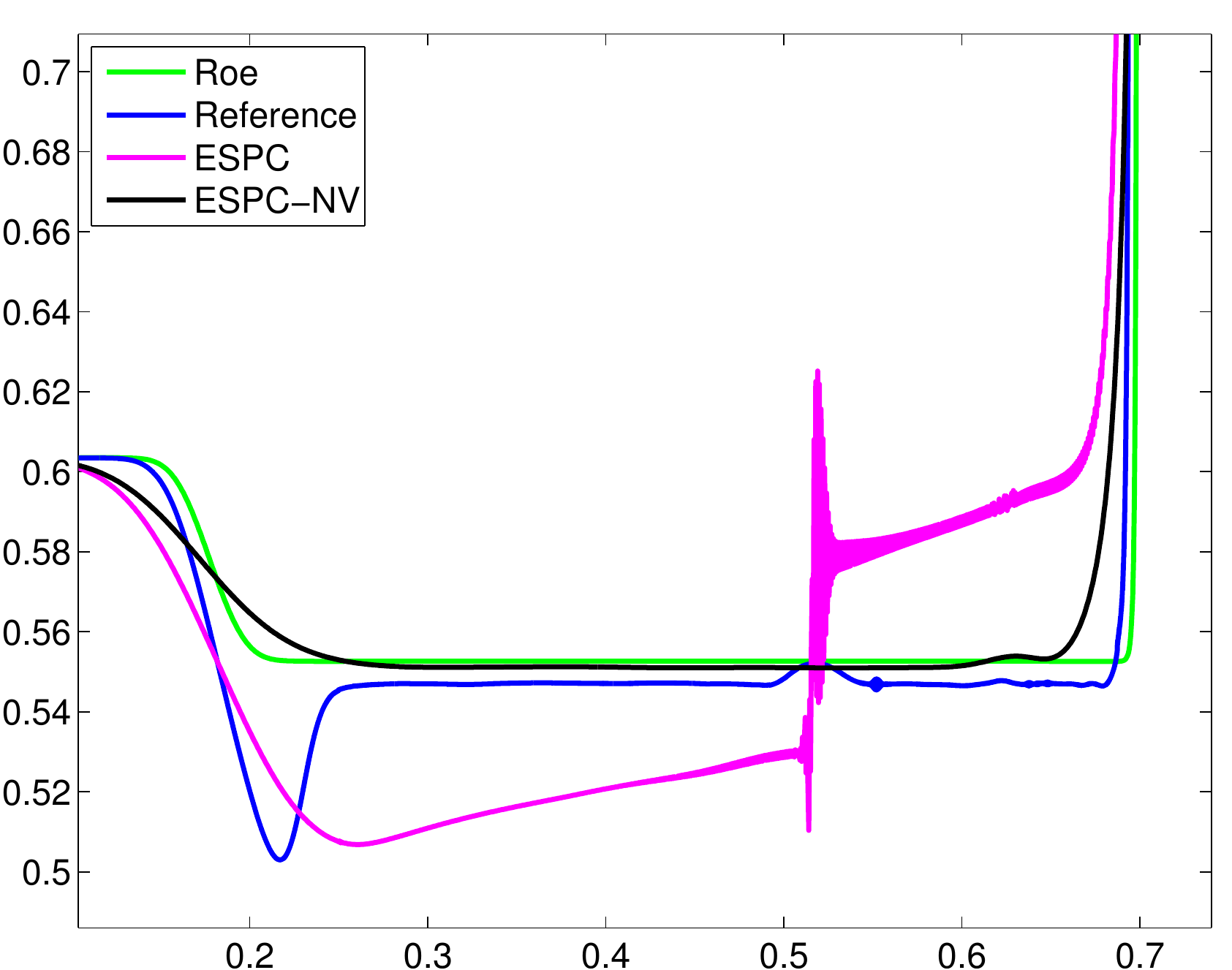}}
\caption{Approximate solutions for height of bottom layer ($h_2$) and total height ($h_1 = h_2$)  for the two-layer shallow water system \eqref{eq:tlsw} with the ESPC, ESPC-NV and path-consistent Roe schemes. A reference solution, computed from the viscous shallow water system \eqref{eq:tlswv} is also displayed.}
\protect \label{fig:tl1}
\end{figure}

\begin{figure}
\includegraphics[width=0.45\linewidth]{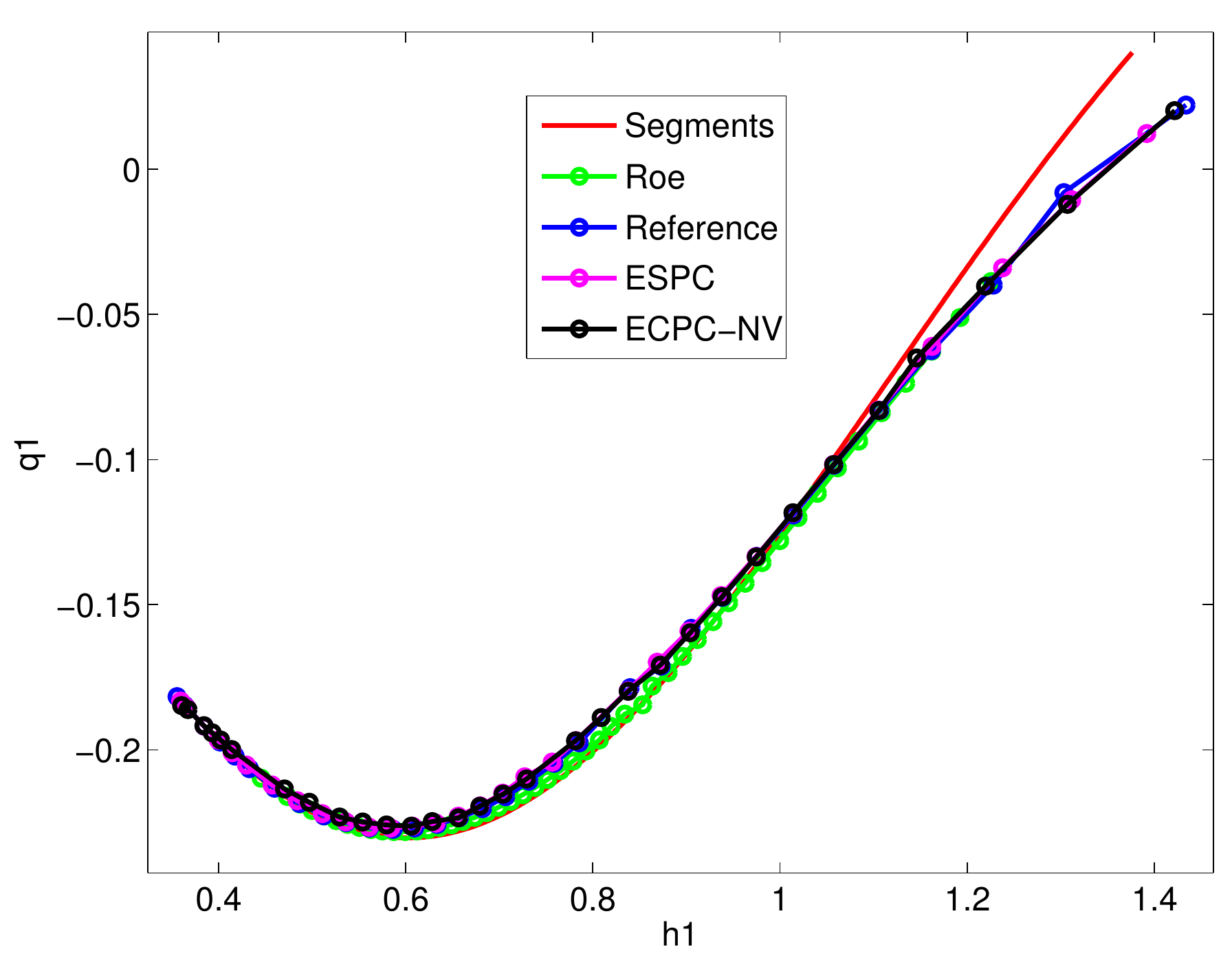}
\includegraphics[width=0.45\linewidth]{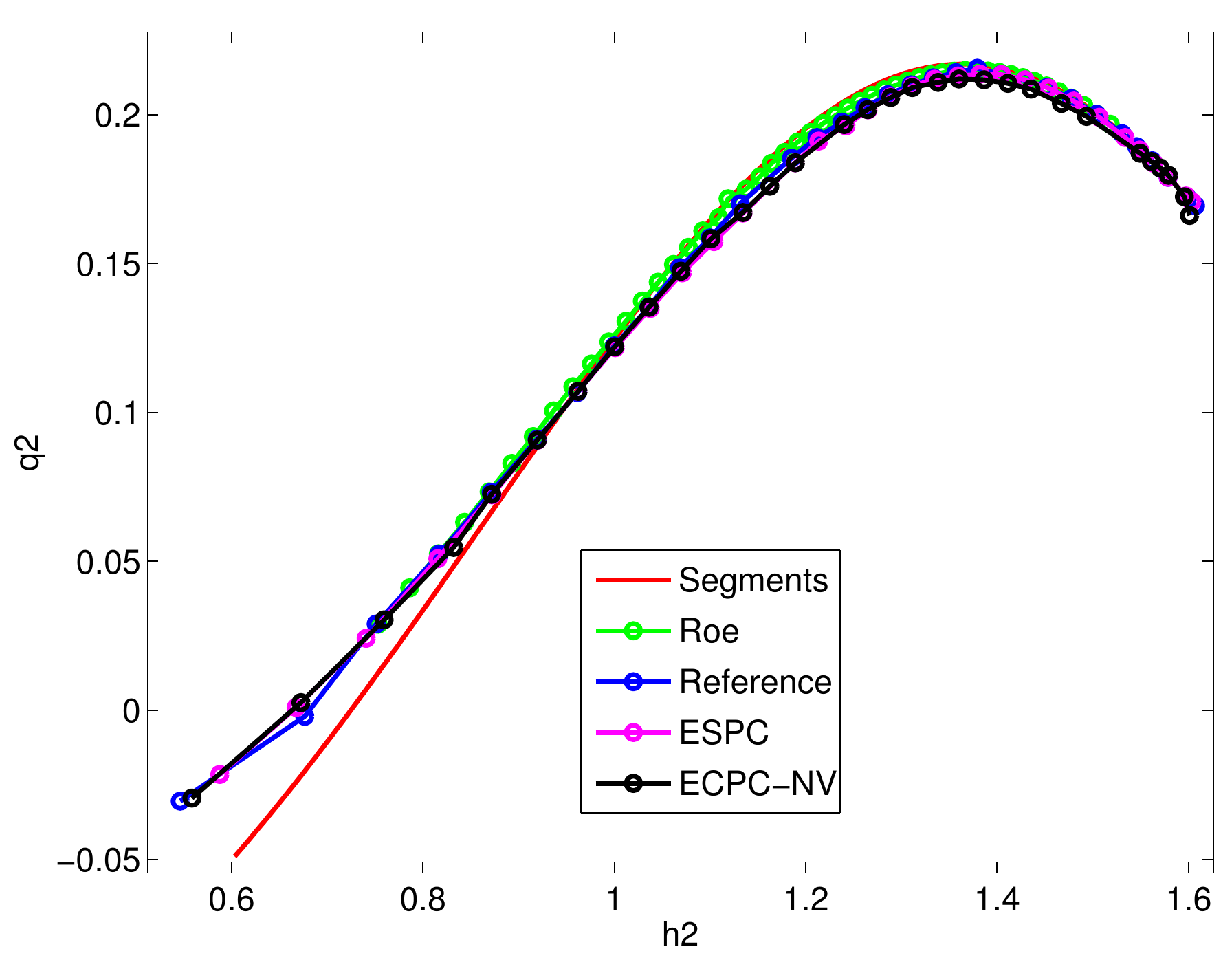}
\includegraphics[width=0.45\linewidth]{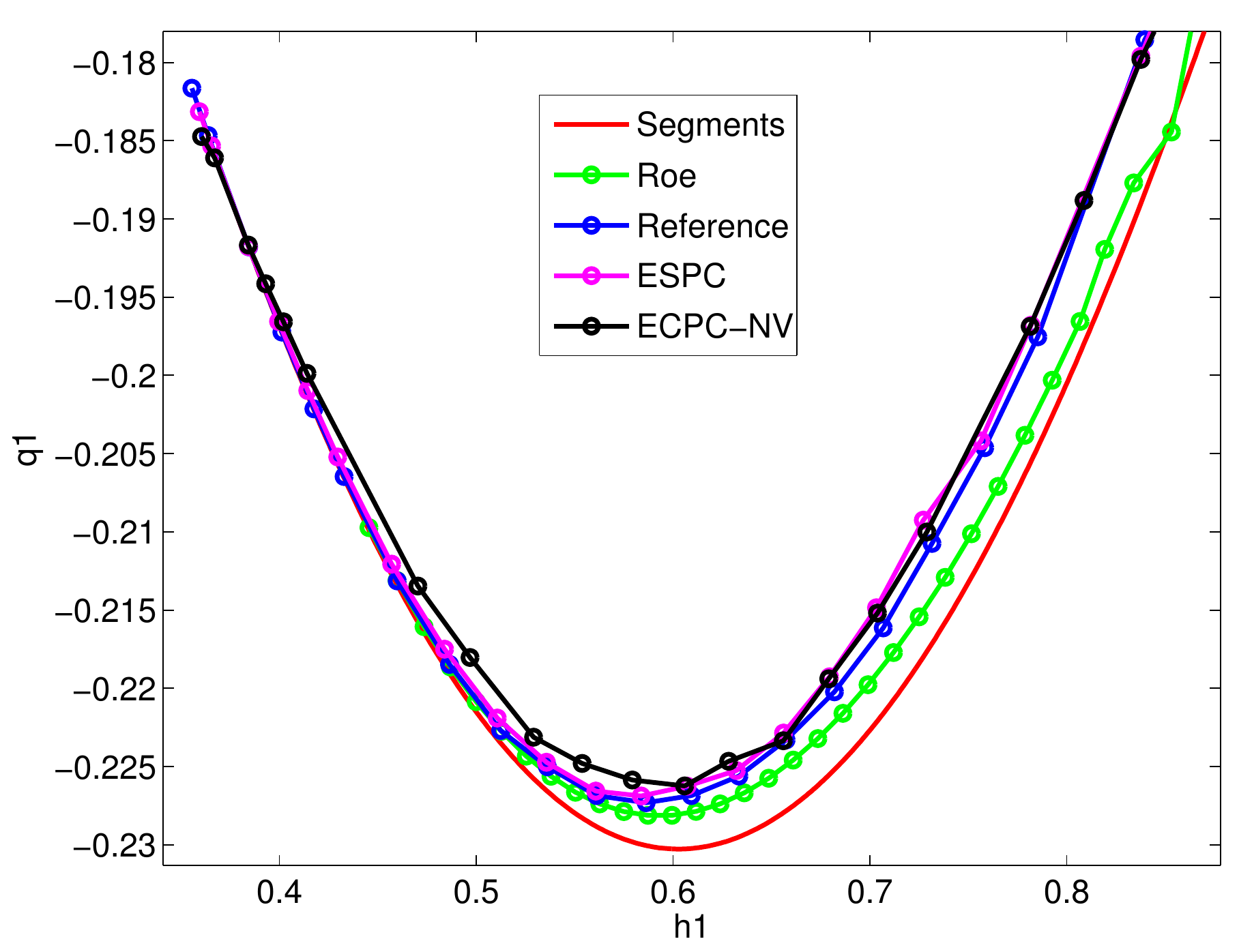}
\includegraphics[width=0.45\linewidth]{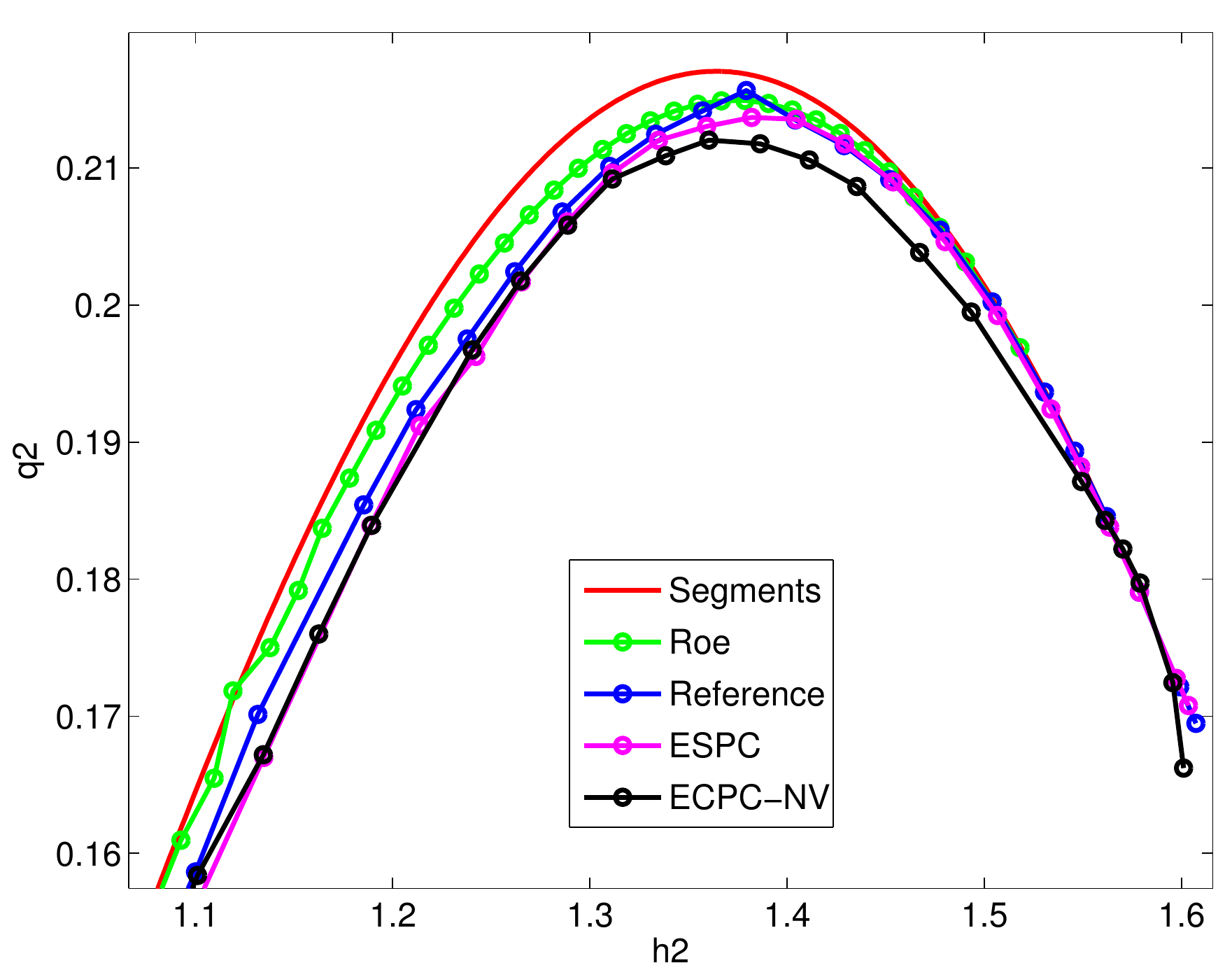}
\caption{Hugoniot loci in the $h_1-(h_1u_1)$ and $h_2-(h_2u_2)$ planes for the two-layer shallow water equations \eqref{eq:tlsw}, computed with the ESPC, ESPC-NV and Roe schemes. A reference Hugoniot locus computed from the viscous shallow water system \eqref{eq:tlswv} and a Hugoniot locus computed using straight line paths are also shown.}
\protect \label{fig:tl2}
\end{figure}

As it is very difficult to compute the viscous profiles explicitly from the viscous system \eqref{eq:tlswv}, we compute the viscous profiles numerically by taking a fixed $\nu \ll 1$ in the ESPC scheme. Note that as the parameter $\nu$ is fixed, the scheme approximates the parabolic system \eqref{eq:tlswv}. The corresponding solution and Hugoniot locus are computed and are labeled \emph{'reference'} in Figures \ref{fig:tl1} and \ref{fig:tl2}.

In order to demonstrate the dependence of the weak solutions to the two-layer shallow water equations on the choice of paths, an alternative path is chosen by fixing a left-hand state $u_l$ and computing numerically the states that can be linked to this state by a shock satisfying the Rankine-Hugoniot conditions associated to the family of \emph{straight line segments}. The computed Hugoniot locus, following Castro et. al (2008), is labeled \emph{`segments'} in Figure \ref{fig:tl2}.

In addition to the ESPC scheme of Castro et al (2013) which lacks any viscosity in the mass balance equations, additional numerical viscosity is added to the mass equations resulting in a  scheme termed as ESPC-NV in the subsequent experiments. To provide a further comparison, we compute the solutions of the two layer shallow water equations with the Roe scheme (consistent with a straight line paths) of Castro, Macias and Pares (2001).

In Figure \ref{fig:tl1}, we plot the solutions obtained with the ESPC, ESPC-NV and Roe schemes for a Riemann problem with initial data 
\be\label{eq:RP-2}
u_l = \begin{pmatrix} 1.376 \\ 0.6035 \\ 0.04019 \\ -0.04906 \end{pmatrix}, \qquad u_r = \begin{pmatrix} 0.37 \\ 1.593 \\ -0.1868 \\ 0.1742 \end{pmatrix}
\ee
and homogeneous Neumann boundary conditions on the computational domain $[0,1]$. All the simulations are performed with $2000$ mesh points. For the sake of comparison, a reference solution computed with the eddy viscosity system \eqref{eq:tlswv} and a fixed $\nu = 2\times 10^{-4}$ on a very fine mesh of $2^{16}$ mesh points is also shown. As seen in this figure, the solutions computed with all the schemes are quite close to the reference solution. As seen in the closeup, there is a minor difference in the intermediate state computed by the ESPC-NV and Roe schemes. The ESPC scheme contains oscillations. This is to be expected as the mass conservation equations contains no numerical viscosity. However, the approximate solution computed by this scheme is still quite close to the reference solution. 

In order to compare the performance of the schemes for a large set of initial data, we compute a \emph{numerical} Hugoniot locus by fixing the same right-hand state as in \eqref{eq:RP}, and then varying the left-hand state. A reference Hugoniot locus, calculated from a numerical approximation of the mixed hyperbolic-parabolic system \eqref{eq:tlswv}, is shown. We also display the Hugoniot locus corresponding to the family of straight line segments. All the Hugoniot loci in the $h_1$-$(h_1u_1)$ plane and the $h_2$-$(h_2u_2)$ plane are shown in Figure \ref{fig:tl2}.

From Figure \ref{fig:tl2}, we observe that the Hugoniot locus calculated using straight line segments is clearly different from the one calculated from the underlying viscous two-layer shallow water equations \eqref{eq:tlswv}. On the other hand, all the three numerical schemes lead to Hugoniot loci that are very close to each other and to the reference Hugoniot locus. Minor differences are visible when we zoom in; see the bottom row of Figure \ref{fig:tl2}. We see that, among the three schemes, the ESPC scheme provides the best overall approximation, to the reference Hugoniot locus. However, both the ESPC-NV and the Roe schemes also provide a good approximation to the reference Hugoniot locus. The results show that (rather surprisingly) the numerical approximation of two-layer shallow water equations is not as sensitive to the viscous terms as the coupled Burgers system. The path-consistent Roe scheme performs adequately in approximating the correct solution. At the same time, the ESPC schemes provide a slightly more accurate approximation.

%-----------------------------------------------------------------------------------------------------------------

\subsection{Schemes with well-controlled dissipation (WCD)}

When approximating shocks to nonconservative hyperbolic systems, the above numerical experiments clearly indicate the superior performance of the {\sl schemes with controlled dissipation,} in contrast with solely imposing a formal path-consistent condition (directy imposed on, for instance, Godunov or Roe flux).  
However, as in the case of nonclassical shocks to conservation laws, schemes with controlled dissipation can fail to approximate strong (large amplitude) shocks. Evidence for this fact has been presented in the numerical experiments for the coupled Burgers' system (see Figure~\ref{fig:cp2}), as well as in Fjordholm and Mishra (2012) who considered the system of gas dynamics in Lagrangian  coordinates. 

Therefore, again following the extensive discussion in Section~3, it is now clear that we need {\sl schemes with well-controlled dissipation} (WCD) in order to approximate shocks with arbitrary large amplitude to nonconservative systems. Schemes with controlled dissipation provide a starting point but need to be further improved. The "full" equivalent equation, that is, a specification of the high-order terms $\mathcal{H}$, analogous to the equivalent equation \eqref{eq:ee} for conservative systems, needs to be taken into account so that the 
asymptotics for large shocks are taken into account in the design of the schemes.
 Furthermore, higher-order finite difference discretizations of the nonconservative term $A(u) u_x$ are required. Only then can we proceed analogously to the conservative case, and balance the leading-order terms of the equivalent equation to the high-order terms in order to design a suitable {\sl WCD condition.} The specification of this WCD condition and extensive experiments with the resulting schemes are the subject of the forthcoming paper (Beljadid, LeFloch, Mishra, and Pares, 2014). 

%==================================================================================

\section{Boundary layers in solutions to systems of conservation laws}
\label{sec:6}

\subsection{Preliminaries}

We now turn our attention to the {\sl initial and boundary value problem} for nonlinear hyperbolic systems of conservation laws \eqref{eq:109} with prescribed data:
\be
\label{eq:cl:bc}
\aligned
        u(0, x) &= u_0(x), \quad x \in \Omega = (X_l, \infty),
\\
 u(t, X_l)&= \overline u(t), \quad t \geq 0, 
\endaligned
\ee
for some given boundary point $X_l \in \RR$. This one-half boundary value problem can be readily generalized to include two boundaries
(that is, $\Omega = (X_l,X_r)$ for some $X_r> X_l$).
 The study of the initial and boundary value problem poses additional difficulties as compared to the study of the Cauchy problem: 
most importantly, the problem~\eqref{eq:109} and \eqref{eq:cl:bc} is ill-posed in the sense that, in general, it admits no solution unless the boundary condition $ u(t, X_l)= \overline u(t)$ is understood in the weaker sense
\bel{eq:DF}
 u(t, X_l) \in \Phi(\overline u(t)), \quad t \geq 0,  
\ee
as proposed by Dubois and LeFloch (1988). Here, $\Phi(\overline u(t))$ is a set containing $\overline u$, which {\sl depends} upon the way that the boundary value problem is handled, for instance via the Riemann problem (sharp boundary layers) or by the vanishing viscosity method (viscous boundary layers. 

Dubois and LeFloch (1988) studied sharp boundary layers for the Euler equations: by solving the so-called boundary Riemann problem and defining $\Phi(\overline u(t))$ as the set of all boundary values of Riemann problems with fixed left-hand state $\overline u(t)$ at $X_l$, these authors show that $\Phi(\overline u(t))$ can be decomposed into several strata (or submanifolds) in the state space 
whose local dimension increases with the strength of the shock layer. 

Alternatively, following (Benabdallah 1986, Dubois and LeFloch 1988, Gisclon 1996, Gisclon and Serre 1994, Joseph and LeFloch  1996,  Joseph and LeFloch 1999), one may consider the viscous approximations  (with $\eps>0$) 
\be
\label{eq:3560}
u^\eps_t + f(u^\eps)_x = \eps \, \bigl( B(u^\eps) \, u^\eps_x\bigr)_x, 
\qquad  x \in \Omega = (X_l, \infty), 
\ee
with given viscosity matrix $B=B(u)$, and supplement these equations with initial and
boundary conditions:  
\be
\label{e:vcl:bd}
\aligned
       u^{\epsilon}(0, x) &= u_0(x), \quad x \in \Omega = (X_l, \infty),
\\
 u^{\epsilon}(t, X_l) &= u_l (t), \quad t \geq 0.
\endaligned
\ee
We assume that this initial and boundary value problem~\eqref{eq:3560}-\eqref{e:vcl:bd} is locally well-posed (which is the case under mild conditions on the viscosity matrix $B$ and with sufficiently regularity on the data) and that for $\epsilon \to 0$ the solutions $u^\eps$ converge 
(in a suitable topology) to a limit $u=u(t,x)$. Due to a boundary layer phenomena, the limit $u$, in general, {\sl does not} satisfy 
the prescribed boundary condition $u_l (t)$ pointwisely. Instead, Dubois and LeFloch (1988)  showed that, if the system of conservation 
laws is endowed with an entropy-entropy flux pair $(U,F)$, then the solutions of the viscous approximation~\eqref{eq:3560} converge as $\epsilon \to 0$) 
to a solution of the initial and boundary 
problem~\eqref{eq:109} and \eqref{eq:cl:bc} in a sufficiently strong topology, 
then the following {\sl entropy boundary inequality} 
\be
    \label{eq:DLF}
F(\overline{u}(t)) - F(u_l(t)) - D_u U(u_l(t)) \cdot \big( f(\overline{u}(t)) -f(u_l(t))\big)  \leq 0.
\ee 
For the case of scalar conservatio laws, a version of this inequality was also derived earlier by Le Roux (1977) and Bardos, Le Roux, and N\'edelec (1979). 

As for other small-scale dependent models before, a major difficulty in the study of initial and boundary value problems was pointed out 
by Gisclon and Serre (1994) and Joseph and LeFloch (1996):  the limit of the viscous approximation~\eqref{eq:3560} \emph{depends} on the underlying viscosity mechanism. In other words, the limit of~\eqref{eq:356} in general changes if one changes the viscosity matrix $B$.

As an example, we consider the linearized shallow water equations \eqref{eq:lsws} with initial data \eqref{eq:lswinit} and boundary data \eqref{eq:ldir}. The system is a \emph{linear}, strictly hyperbolic, $2 \times 2$ system and is the simplest possible problem that can be considered in this context. We consider two different viscosity operators: an \emph{artificial} uniform (Laplacian) viscosity \eqref{eq:lswslap} and the \emph{physical} 
eddy viscosity \eqref{eq:lswsed}. The resulting limit solutions are shown in the left of Figure~\ref{fig:71}. As shown in the figure, there is a significant difference in solutions (near the boundary) corresponding to different viscosity 
operators.  

For an extended discussion of the initial boundary value problem for 
systems of conservation laws and its viscous approximation, we refer to Serre (2000, 2007) and the bibliography therein, including for the
 theory of discrete shock profiles. We stress that analytically establishing the convergence $\epsilon \to 0$ 
in the possibly characteristic regime and general diffusion matrices 
was addressed only rather recently in the successive works: Joseph and LeFloch (1999, 2002, 2006), 
Ancona and Bianchini (2006),  
Bianchini and Spinolo (2009), 
and Christoforou and Spinolo (2011).   
Many other results are also available in more specific cases, which we do not attempt to review here. 

%----------------------------------------------------------------------------------------------

\subsection{Standard finite difference (or finite volume) schemes}

As in the previous sections, it is common to use a finite difference or a finite volume scheme of the form \eqref{eq:fdm1}, with a suitable numerical flux $g_{i+1/2}$. Following Goodman (1982) and Dubois and LeFloch (1988) (see also the textbook by LeVeque 2003), the Dirichlet boundary conditions at $X=X_l$ are imposed by setting in the \emph{ghost} cell $[x_{-1/2},x_{1/2}]$:
\be
\label{eq:bcnum}
u^n_0 = u_l(t^n).
\ee

\begin{figure}[htbp]
\centering
\subfigure[Viscous profile]{\includegraphics[width=0.48\linewidth]{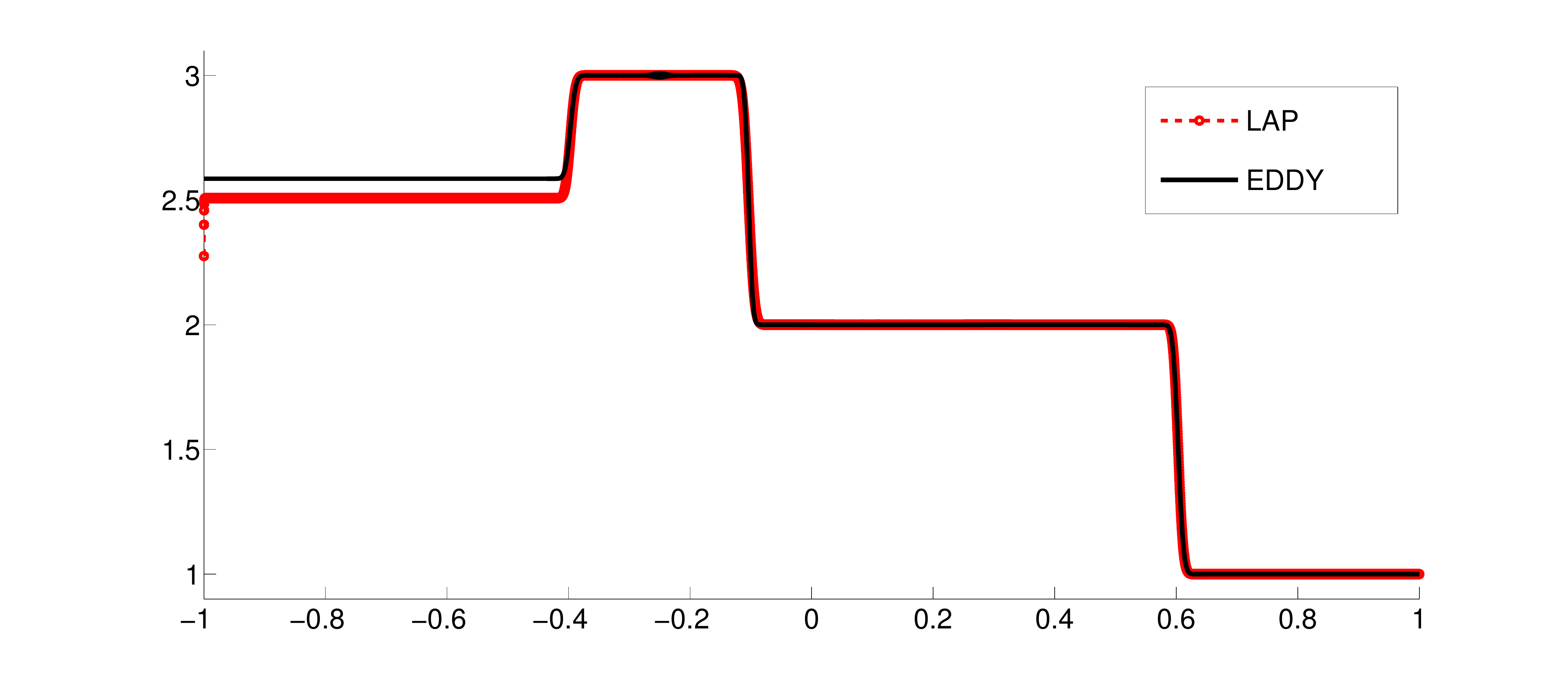}}
\subfigure[Roe scheme]{\includegraphics[width=0.48\linewidth]{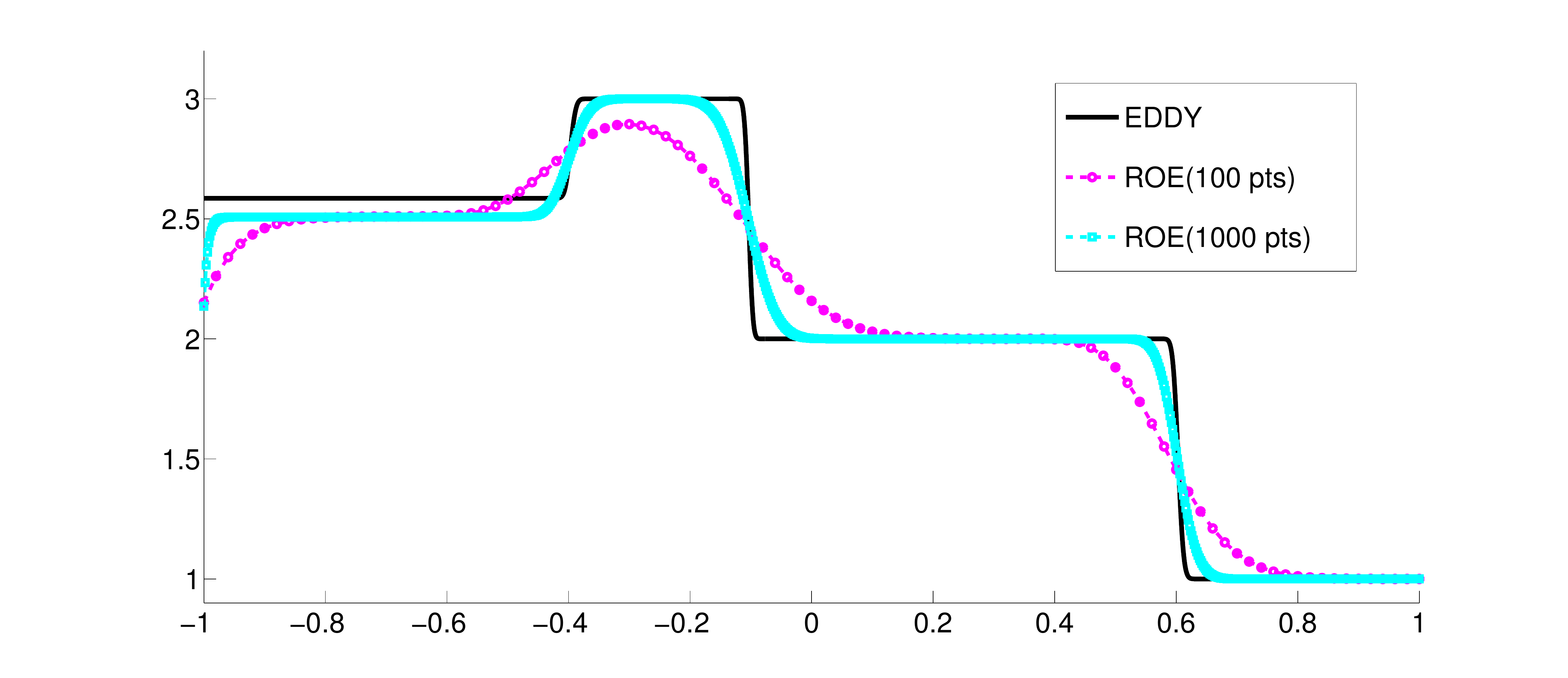}}
\caption{Linearized shallow water equations. 
\newline 
Left: Viscous profile at $t=0.25$ with viscosity \eqref{eq:lswslap} and 
eddy viscosity \eqref{eq:lswsed}. 
\newline 
Right: Roe (Godunov) scheme with the same data.}
\label{fig:71}
\end{figure}

One might expect that using a flux based on the Riemann solver at the boundary should suffice to approximate the correct solution of the 
initial boundary value problem. 
However, standard numerical schemes may not converge to the physical viscosity  solution of the initial boundary value problem for a system of conservation laws. We illustrate this by again considering the linearized shallow water equations \eqref{eq:lsws} with initial data 
\eqref{eq:lswinit} and boundary data \eqref{eq:ldir} at time $t=0.25$. The results with a standard Roe (Godunov) scheme for this linear system are presented in Figure~\ref{fig:71} (right). The figure clearly shows that the Roe scheme converges to a solution that is different from the physical-viscosity  solution of the system, realized as a limit of the eddy viscosity approximation \eqref{eq:lswsed}. In fact, the solution converges to the limit of the \emph{artificial} uniform viscosity approximation \eqref{eq:lswslap}.

As in the previous examples of numerical approximation of small-scale dependent shock waves, this failure to approximate the correct solution can be explained in terms of the equivalent equation \eqref{eq:ee} of the scheme \eqref{eq:fdm1}. As explained before, the numerical viscosity of the scheme may not match with the underlying physical viscosity $B$ in \eqref{eq:356}. As the solutions of the boundary value problems depend explicitly on the underlying small-scale mechanism, it is not surprising that the scheme fails to approximate the physically relevant solution as shown in Figure~\ref{fig:71}. 

%-----------------------------------------------------------------------------------------------------------------

\subsection{Schemes with controlled dissipation}

As in the previous examples, the equivalent equation suggests an approach to design numerical schemes which will approximate the physically relevant solution. As in the previous sections, this approach consists of the following stages.

\subsubsection{Design of entropy conservative scheme}

As outlined before, the first step to choose a numerical flux $g^{\ast}_{i+1/2}$ such that the resulting finite difference scheme \eqref{eq:fdm1} satisfies a discrete version of the entropy identity \eqref{eq:efdm2}. The construction of such schemes follows Tadmor (1987) and already has been outlined. We follow this recipe for our construction and we use explicit formulas for the entropy conservative flux.

\subsubsection{Correct numerical diffusion operator}

As before, we need add numerical diffusion to stabilize the entropy conservative scheme. To this end, we choose the following flux:
\be
\label{eq:esf71}
g_{i+1/2} = g^{\ast}_{i+1/2} - \frac{1}{2} c_{\max} B(\widehat{u}_{i+1/2})[[u]]_{i+1/2}.
\ee
Here, $B$ is the underlying (small-scale) physical viscosity in \eqref{eq:3560}. 
The resulting numerical scheme \eqref{eq:fdm1} with numerical flux \eqref{eq:esf71}, has been termed as the \emph{CND} scheme ('correct numerical diffusion' scheme) and was shown by Mishra and Spinolo (2011)  to be i) entropy stable and ii) its equivalent equation matched the underlying parabolic regularization \eqref{eq:356} of the conservation law \eqref{eq:109}. Consequently, this scheme can be termed as a scheme with controlled dissipation in our terminology. 

%--------------------------------------------------------------------------------------------------------------------

\subsection{Numerical experiments}

\subsubsection{Linearized shallow water equations}

We consider the linearized shallow water equations of fluid flow (LeVeque 2003):
\be
\label{eq:lsws}
\begin{aligned}
h_t + \widetilde{u} h_x + \widetilde{h} u_x &= 0, \\
u_t + g h_x + \widetilde{u} u_x &= 0,
\end{aligned}
\ee 
where $h$ represents the height and $u$ the water velocity. The constant $g$ stands for the acceleration due to gravity and  $\widetilde{h},\widetilde{u}$ are the (constant) height and velocity states around which the shallow water
equations are linearized.

The physical-viscosity  viscosity mechanism for the shallow water system is the \emph{eddy viscosity}. Adding eddy viscosity to the linearized shallow water system results in the following mixed hyperbolic-parabolic system:
\be
\label{eq:lswsed}
\begin{aligned}
h_t + \widetilde{u} h_x + \widetilde{h} u_x &= 0,\\
u_t + g h_x + \widetilde{u} u_x &= \epsilon u_{xx}.
\end{aligned}
\ee
On the other hand, for the sake of comparison, we can also add an \emph{artificial} viscosity to the linearized shallow waters by including a Laplacian regularization, that is, 
\be
\label{eq:lswslap}
\begin{aligned}
h_t + \widetilde{u} h_x + \widetilde{h} u_x &= \epsilon h_{xx}, \\
u_t + g h_x + \widetilde{u} u_x &= \epsilon u_{xx}.
\end{aligned}
\ee

For the rest of this section, we specify the parameters
\be
\label{eq:lsw:wtl}
\widetilde{h} = 2,\quad \widetilde{u} = 1, \quad g = 1
\ee
and the initial data
\be
\label{eq:lswinit}
(h,u) (0, x) = \begin{cases}
                     U^- = (3,1),  \qquad & x < 0, \\
                     U^+= (1,1), & x > 0,
                     \end{cases}
\ee
together with the Dirichlet boundary data 
\be
\label{eq:ldir}
   (h, u) (-1, t) = U_l (t) =(2,  1), \qquad  \, t > 0.
\ee
As the linearized shallow water equations \eqref{eq:lsws} are a linear constant-coefficient system of equations, one can explicit solve the above initial and boundary value problem (see Mishra and Spinolo 2011) for limits of the eddy viscosity as well as the uniform viscosity. These exact solutions are used as reference solutions.  The boundary condition then holds only in the weak sense of Dubois and LeFloch \eqref{eq:DF} with suitably defined boundary layer sets.  

\begin{figure}[htbp]
\centering
\subfigure[Height ($h$)]{\includegraphics[width=9cm]{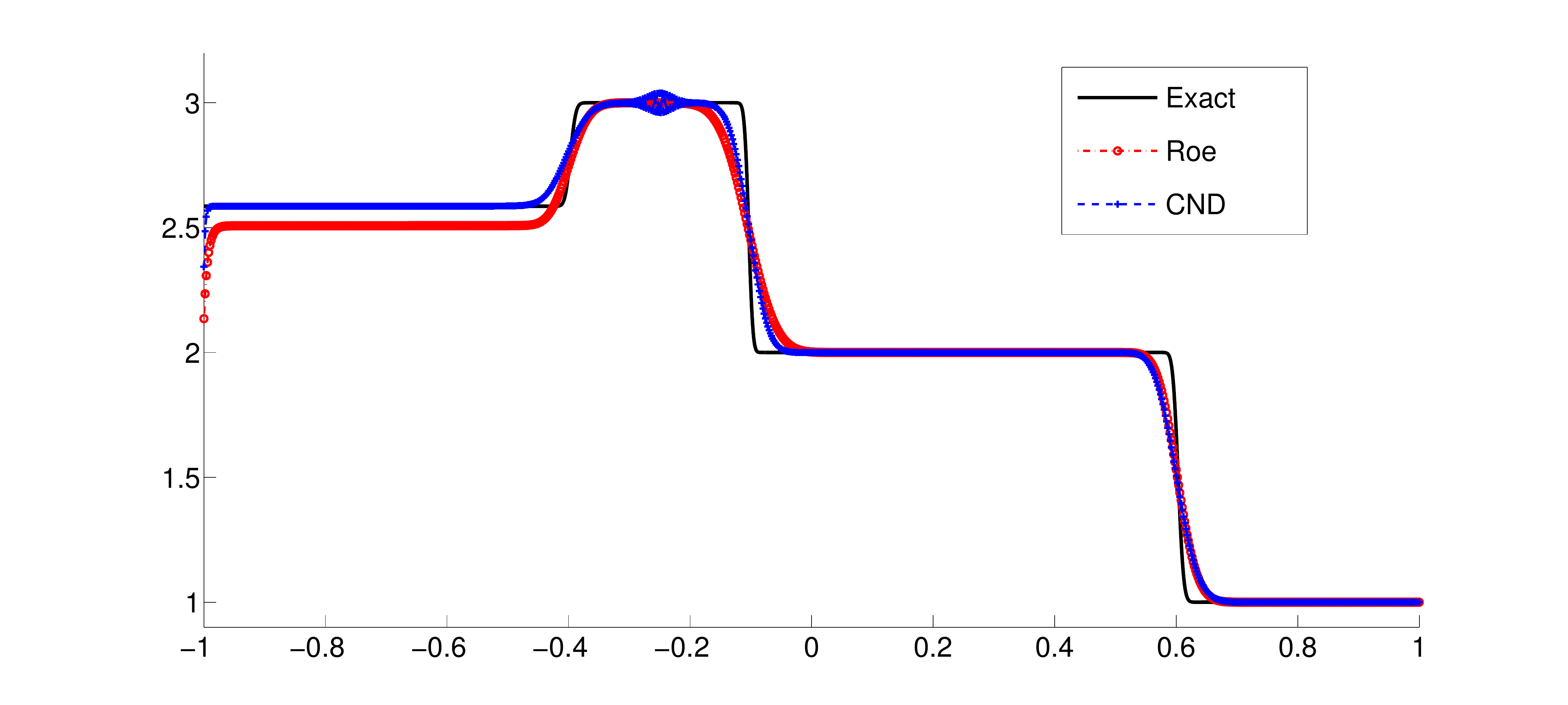}} \\
\subfigure[Velocity ($u$)]{\includegraphics[width=9cm]{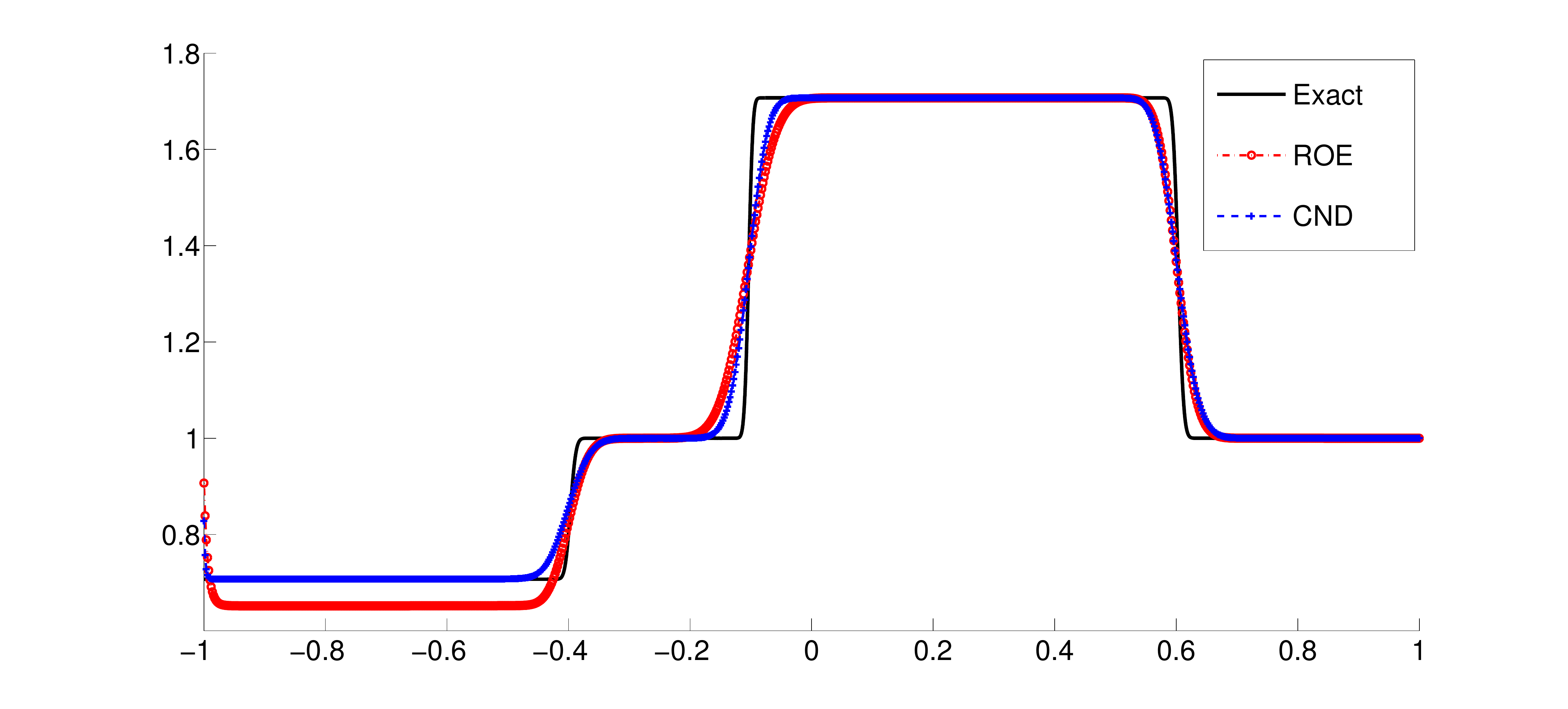}}
\caption{Solutions of the linearized shallow water equations \eqref{eq:lsws} at time $t=0.25$ with initial data \eqref{eq:lswinit} and boundary data \eqref{eq:ldir} computed with the Roe and CND schemes with $1000$ mesh points. The exact solution computed is provided for comparison.}
\label{fig:72}
\end{figure}

The numerical solutions computed with the standard Roe scheme and the CND scheme  at time $t=0.25$ are shown in Figure~\ref{fig:72}. As we are interested in computing the physical-viscosity  solutions of the linearized shallow water equations, obtained as a limit of the eddy viscosity \eqref{eq:lswsed}. Observe that the full Dirichlet boundary conditions can be imposed at the boundary, even for the case of eddy viscosity as $\tilde{u} > 0$. Both the numerical solutions are computed with a $1000$ mesh points.

\begin{figure}[htbp]
\centering
\subfigure[Height ($h$)]{\includegraphics[width=9cm]{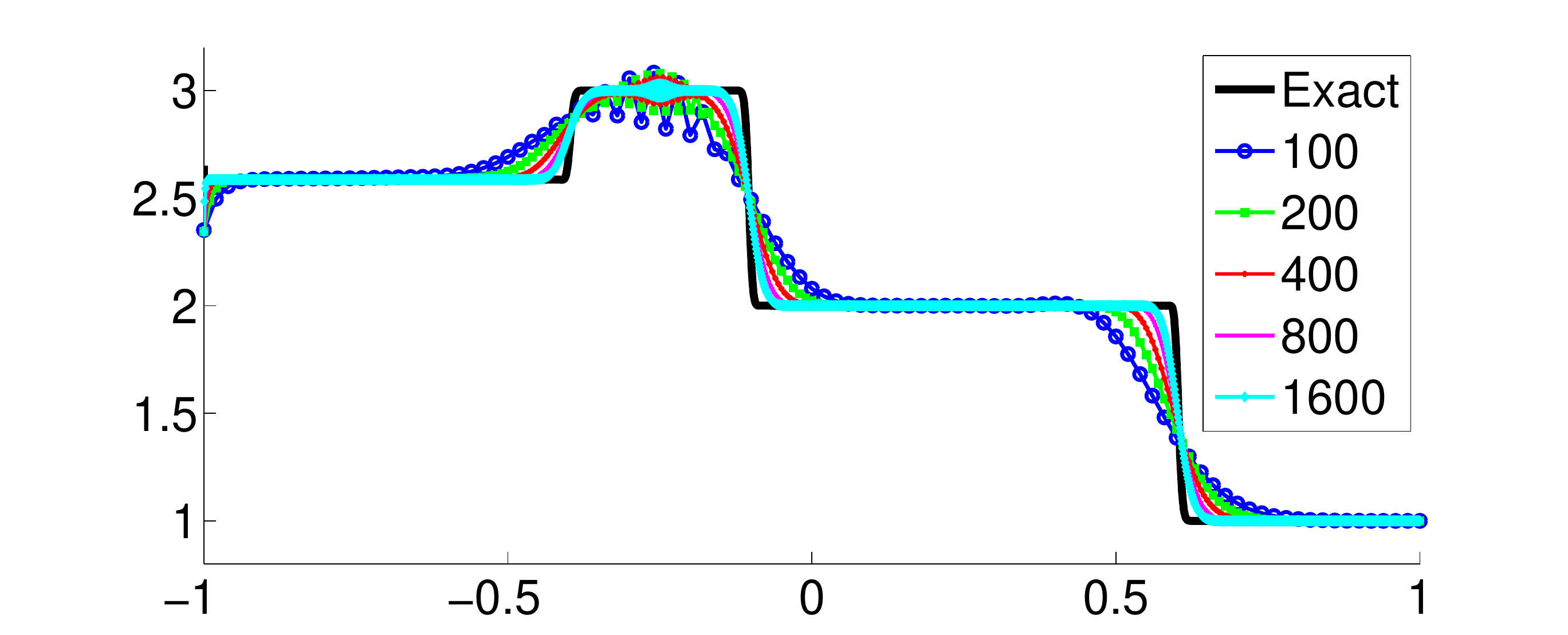}} \\
\subfigure[Velocity ($u$)]{\includegraphics[width=9cm]{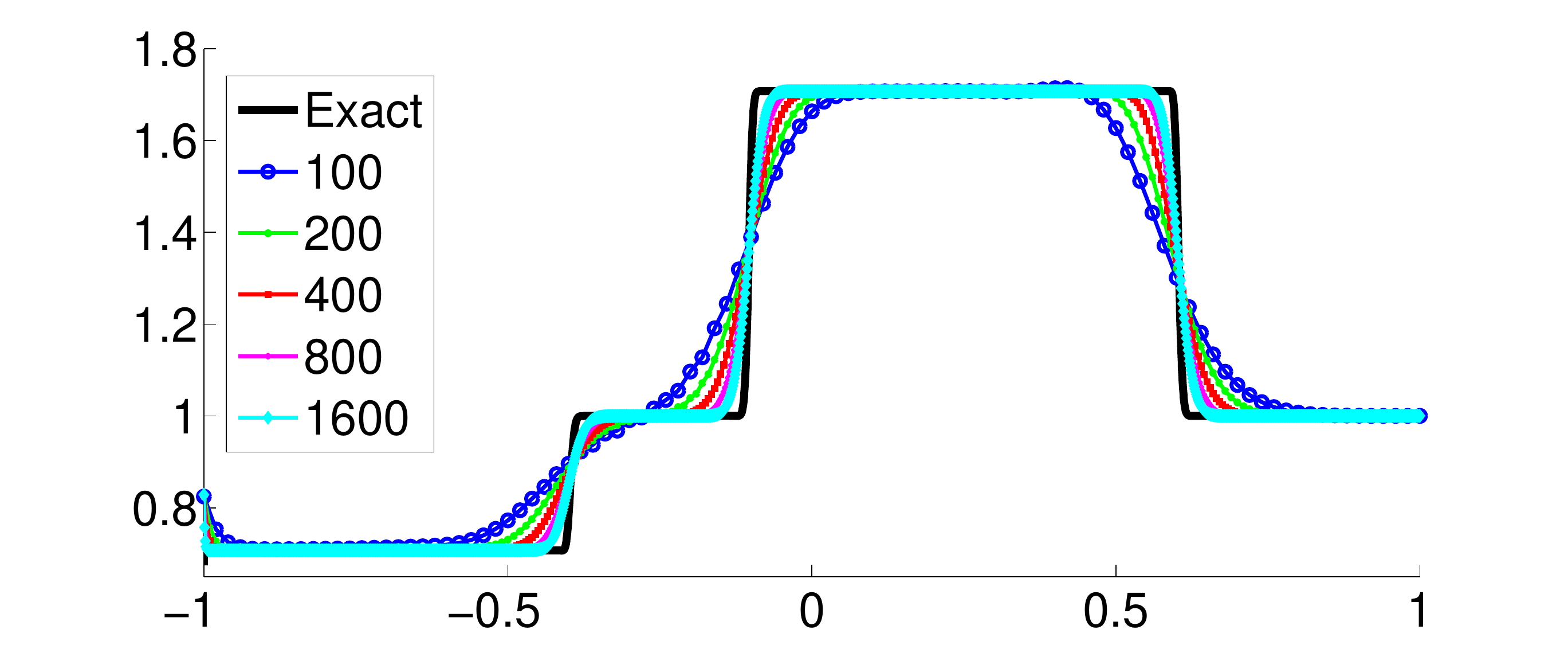}}
\caption{Solutions of the linearized shallow water equations \eqref{eq:lsws} at time $t=0.25$ with initial data \eqref{eq:lswinit} and boundary data \eqref{eq:ldir} computed with the CND scheme at different mesh resolutions. The exact solution is provided for comparison.}
\label{fig:720}
\end{figure}

The results in Figure~\ref{fig:72} clearly show that the Roe scheme does not converge to the desired solution, realized as the limit of the physical-viscosity  \eqref{eq:lswsed}. On the other hand, the solutions computed with the CND scheme approximate the physical-viscosity  solution quite well. There are some small amplitude oscillations in the height with the CND scheme. This is a consequence of the singularity of the viscosity matrix $B$ in this case. As there is no numerical viscosity in the scheme approximating the height, this results in small amplitude oscillations. These small amplitude oscillations might lead to small amplitude oscillations in the velocity. However, these oscillations are damped considerably due to the numerical viscosity used to approximate the velocity. At this resolution, it is not possible to observe these oscillations in the velocity. A mesh refinement study, the results of which are plotted in Figure~\ref{fig:720}, further establishes that the approximate solutions generated by the CND scheme does converge to the physically relevant solution (the limit of the eddy viscosity approximation) as the mesh is refined. Furthermore, the amplitude of the height oscillations reduces as the mesh is refined. 

\subsubsection{Nonlinear Euler equations}

The Euler equations of gas dynamics (in one space dimension) read 
\be
\label{eq:euler}
\begin{aligned}
\rho_t + (\rho u)_x &= 0, \\
(\rho u)_t + (\rho u^2 + p)_x &= 0, \\
E_t   + ((E+p)u)_x &= 0,
\end{aligned}
\ee
where $\rho$ denotes the fluid density and $u$ the fluid velocity. The total energy $E$ and the pressure $p$ are related by the ideal gas equation of state
\bel{eq:eos}
E = \frac{p}{\gamma -1 } + \frac{1}{2} \rho u^2,
\ee
the so-called adiabatic constant $\gamma>1$ being a constant specific of the gas. The system is hyperbolic with eigenvalues
\be
\label{eq:euleig}
\lambda_1 = u - c, \quad \lambda_2 = u, \quad \lambda_3 = u +c,
\ee
where $c = \sqrt{\frac{\gamma p}{\rho}}$ is the sound speed.
Furthermore, the equations are augmented with the entropy inequality
 \be
 \label{eq:eulent}
 \left(\frac{-\rho s}{\gamma - 1} \right)_t + \left(\frac{-\rho u s}{\gamma - 1}\right)_x \leq 0,
 \ee
with thermodynamic entropy
$$
s = \log(p) - \gamma \log(\rho).
$$

The compressible Euler equations are derived by ignoring kinematic viscosity and heat conduction. Taking these small-scale effects into account results in the \emph{Navier-Stokes system} for compressible fluids: 
\be
\label{eq:cns}
\begin{aligned}
\rho_t + (\rho u)_x &= 0, \\
(\rho u)_t + (\rho u^2 + p)_x &= \nu u_{xx}, \\
E_t   + ((E+p)u)_x &= \nu \left(\frac{u^2}{2}\right)_{xx} + \kappa \theta_{xx}.
\end{aligned}
\ee
Here, $\theta$ represents the temperature 
$$
\theta = \frac{p}{(\gamma-1)\rho}, 
$$
while $\nu$ is the viscosity coefficient and $\kappa$ the coefficient of heat conduction.
For the sake of comparison, we can also add a uniform (Laplacian) diffusion to obtain the compressible Euler equations with \emph{artificial} viscosity:
\be
\label{eq:eulerlap}
\begin{aligned}
\rho_t + (\rho u)_x &= \epsilon \rho_{xx}, \\
(\rho u)_t + (\rho u^2 + p)_x &= \epsilon (\rho u)_{xx}, \\
E_t   + ((E+p)u)_x &= \epsilon E_{xx}.
\end{aligned}
\ee

\begin{figure}[htbp]
\centering
\subfigure[Density ($\rho$)]{\includegraphics[width=9cm]{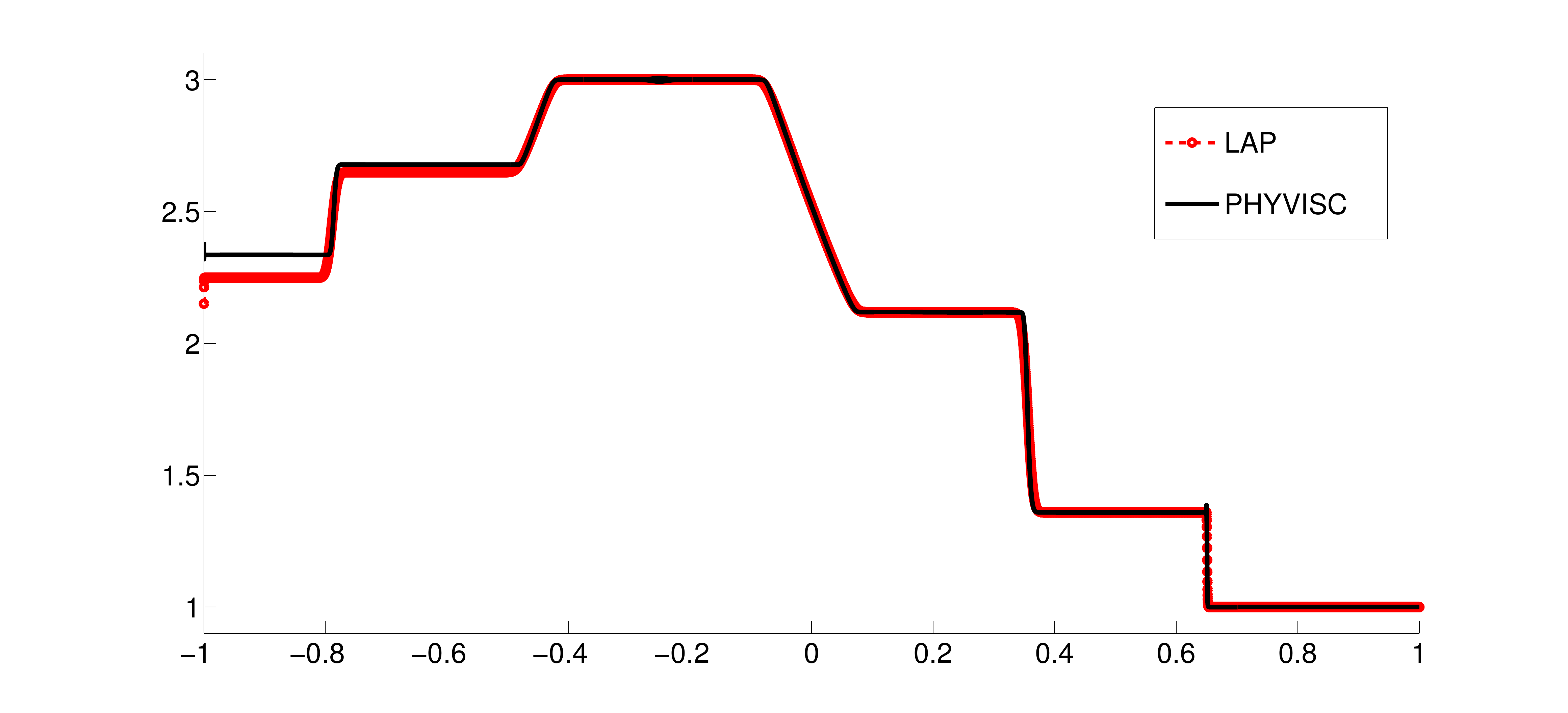}} \\
\subfigure[Velocity ($u$)]{\includegraphics[width=9cm]{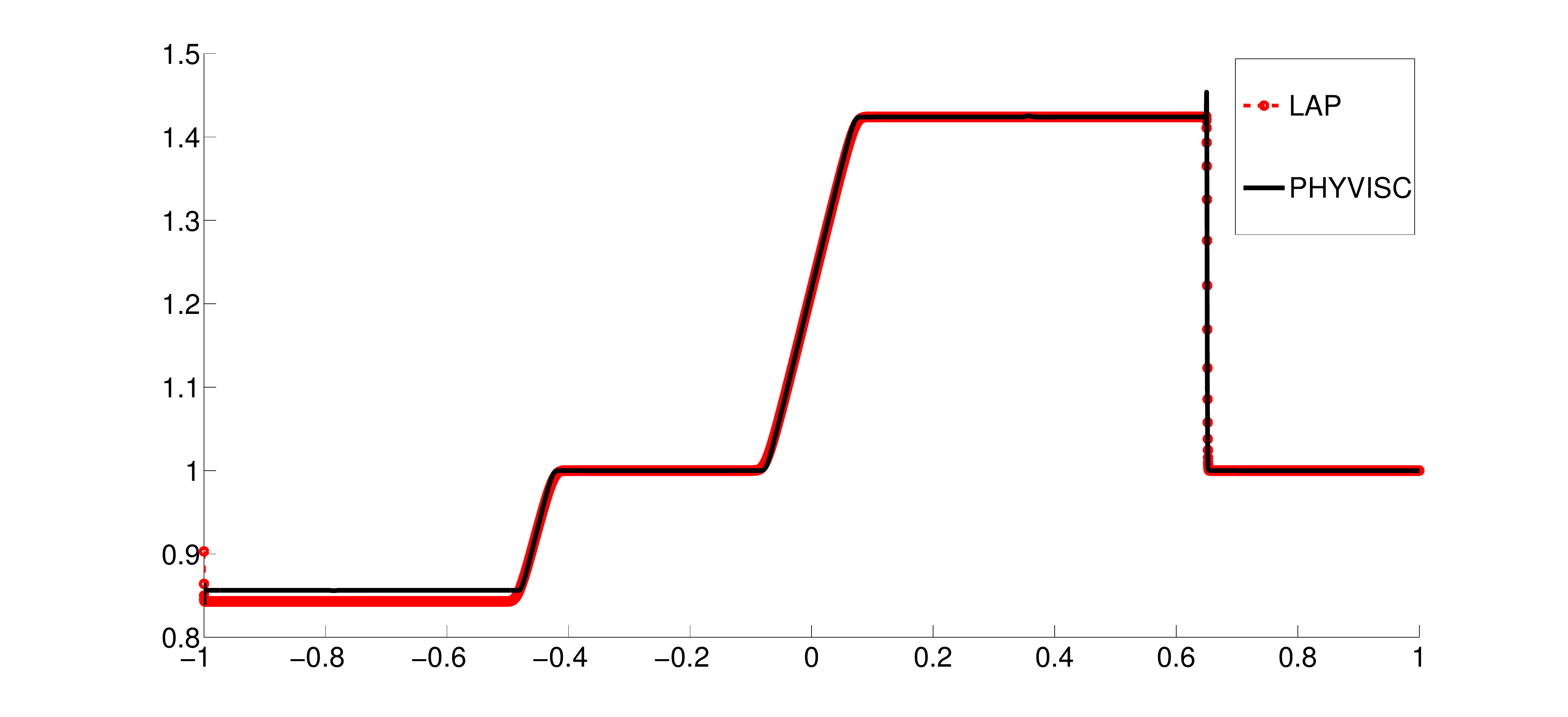}} \\
\subfigure[Pressure ($p$)]{\includegraphics[width=9cm]{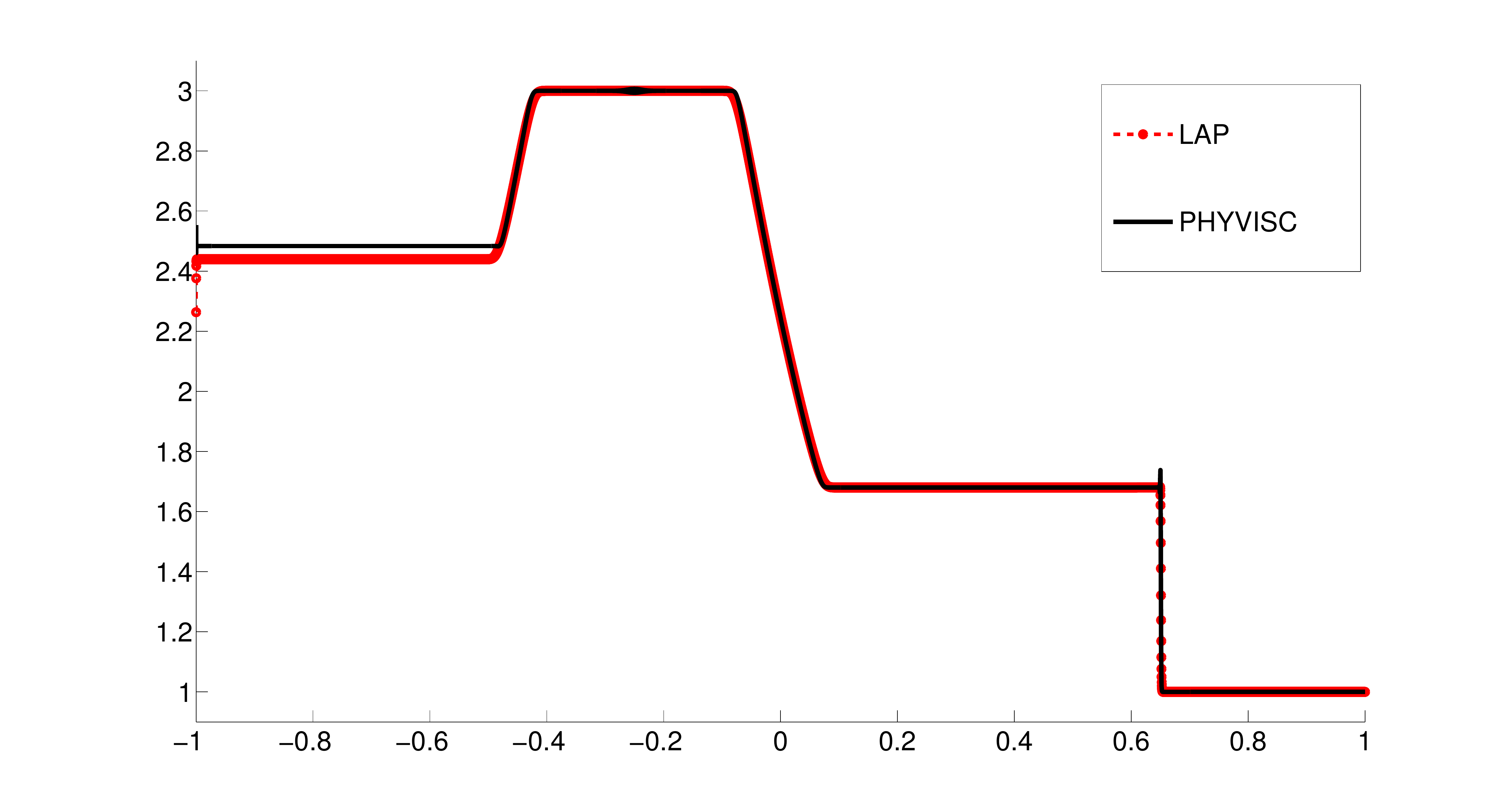}} 
\caption{Limit solutions at time $t=0.25$  of the compressible Euler equations \eqref{eq:euler} with initial data \eqref{eq:eulinit} and boundary data \eqref{eq:eulbd}. The limits of the physical viscosity i.e compressible Navier-Stokes equations \eqref{eq:cns} and the artificial (Laplacian) viscosity \eqref{eq:eulerlap} are compared.}
\label{fig:73}
\end{figure}

Although explicit solutions are not know (even for the boundary Riemann problem), we can still rely on a numerical approximation of the regularized equations \eqref{eq:eulerlap} and \eqref{eq:cns}, as done in Mishra and Spinolo (2011). To illustrate the difference in limit solutions for different regularizations, we consider both \eqref{eq:eulerlap} and \eqref{eq:cns} in the domain $[-1,1]$ with initial data
\be
\label{eq:eulinit}
(\rho_0,u_0,p_0) = \begin{cases}
                     (3.0,1.0,3.0), \qquad &\; x < 0, \\
                     (1.0,1.0,1.0), & \; x > 0.
                     \end{cases}
\ee
We impose absorbing boundary conditions at the right-hand boundary and Dirichlet boundary conditions at the left-hand boundary with boundary data
\be
\label{eq:eulbd}
\Big(\rho(-1,t),u(-1,t),p(-1,t) \Big)= (2.0,1.0,2.0), 
\ee
and, for simplicity, we set $\nu = \kappa = \epsilon$. The results for the finite difference scheme approximating the uniform viscosity  \eqref{eq:eulerlap} and the physical viscosity \eqref{eq:cns} at time $t=0.25$ are presented in Figure \ref{fig:73}. The figure shows that the there is a clear difference in the limit solutions of this problem, obtained from the compressible Navier-Stokes equations \eqref{eq:cns} and the Euler equations with artificial viscosity \eqref{eq:eulerlap}. The difference is more pronounced in the density variable near the left-hand boundary. Both the limit solutions were computed by setting $\epsilon = 10^{-5}$ and on a very fine mesh of $32000$ points. 

The above example also illustrates the limitations of using a mixed hyper\-bolic-parabolic system like the compressible Navier-Stokes equations \eqref{eq:cns}. In order to resolve the viscous scales, we need to choose $\Dx = \mathcal{O}\left(\frac{1}{\epsilon} \right)$, with $\epsilon$ being the viscosity parameter. As $\epsilon$ is very small in practice, the computational effort involved is prohibitively expensive. In the above example, we needed $32000$ points to handle $\epsilon = 10^{-5}$. Such ultra fine grids are not feasible, particularly in several space dimensions.

Hence, we will use a scheme with controlled dissipation such as the CND scheme. This scheme requires both entropy conservative fluxes as well as numerical diffusion operators. We use the entropy conservative fluxes for the Euler equations that were recently developed by Ismail and Roe (2007) and a discrete version of the Navier-Stokes viscosity \eqref{eq:cns} as in Mishra and Spinolo (2011). 

\begin{figure}[htbp]
\centering
\subfigure[Density ($\rho$)]{\includegraphics[width=9cm]{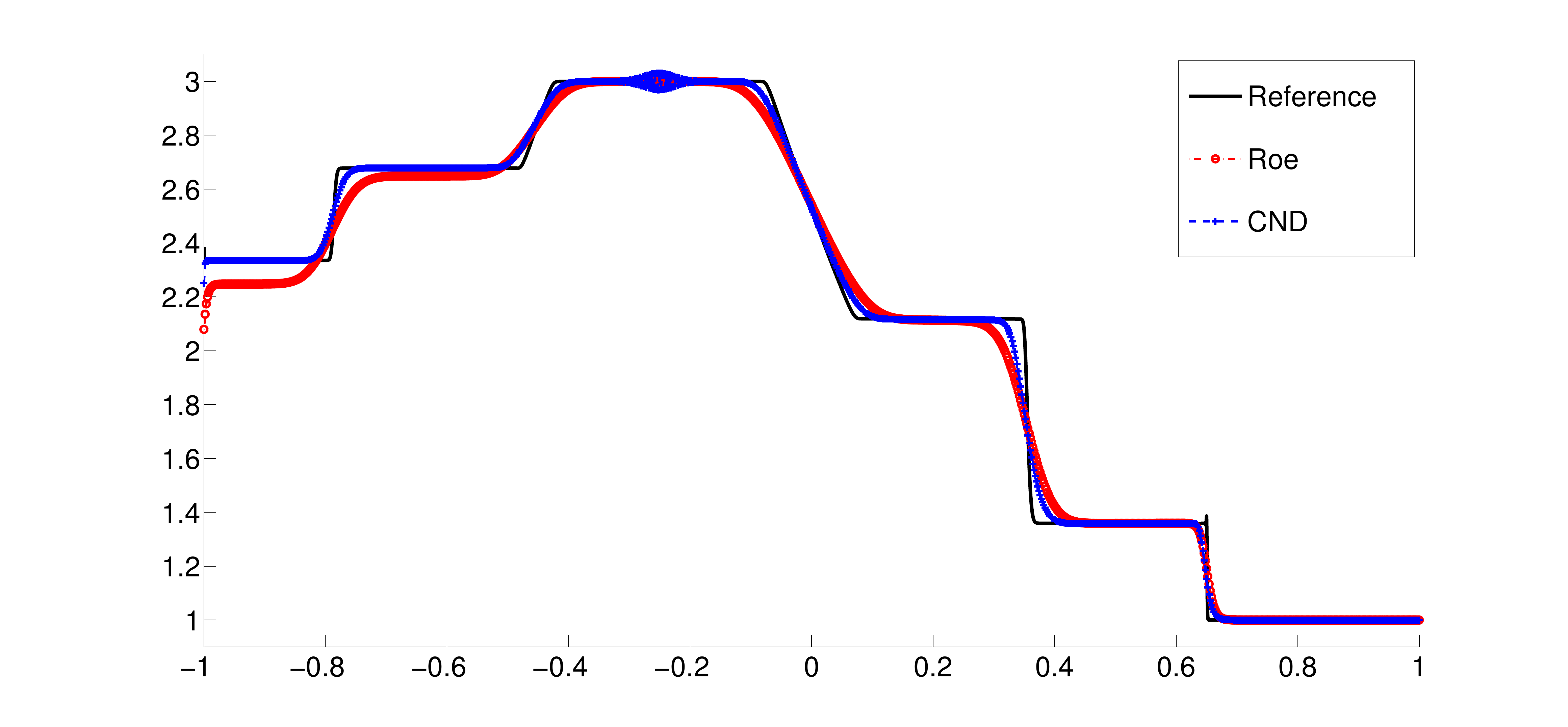}} \\
\subfigure[Velocity ($u$)]{\includegraphics[width=9cm]{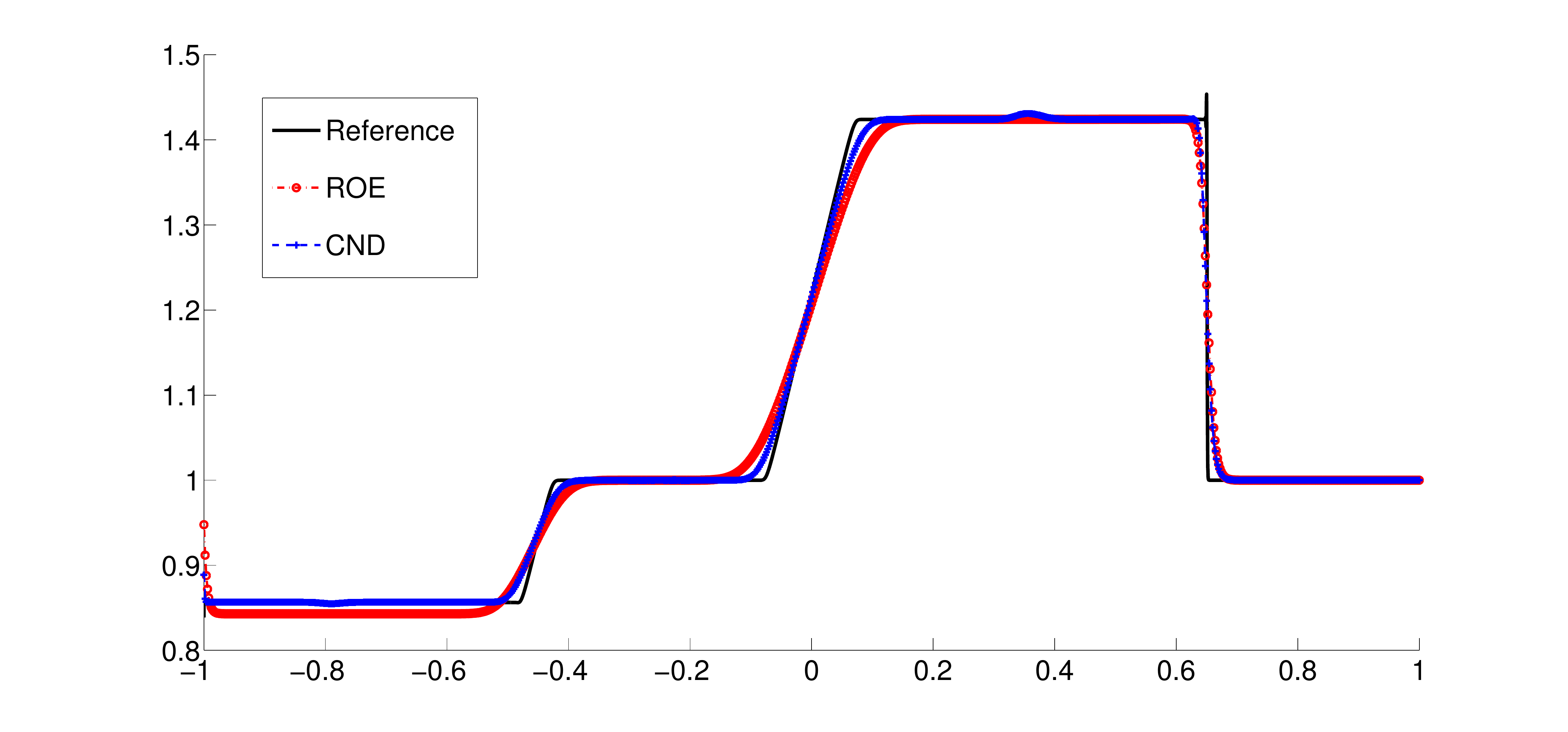}} \\
\subfigure[Pressure ($p$)]{\includegraphics[width=9cm]{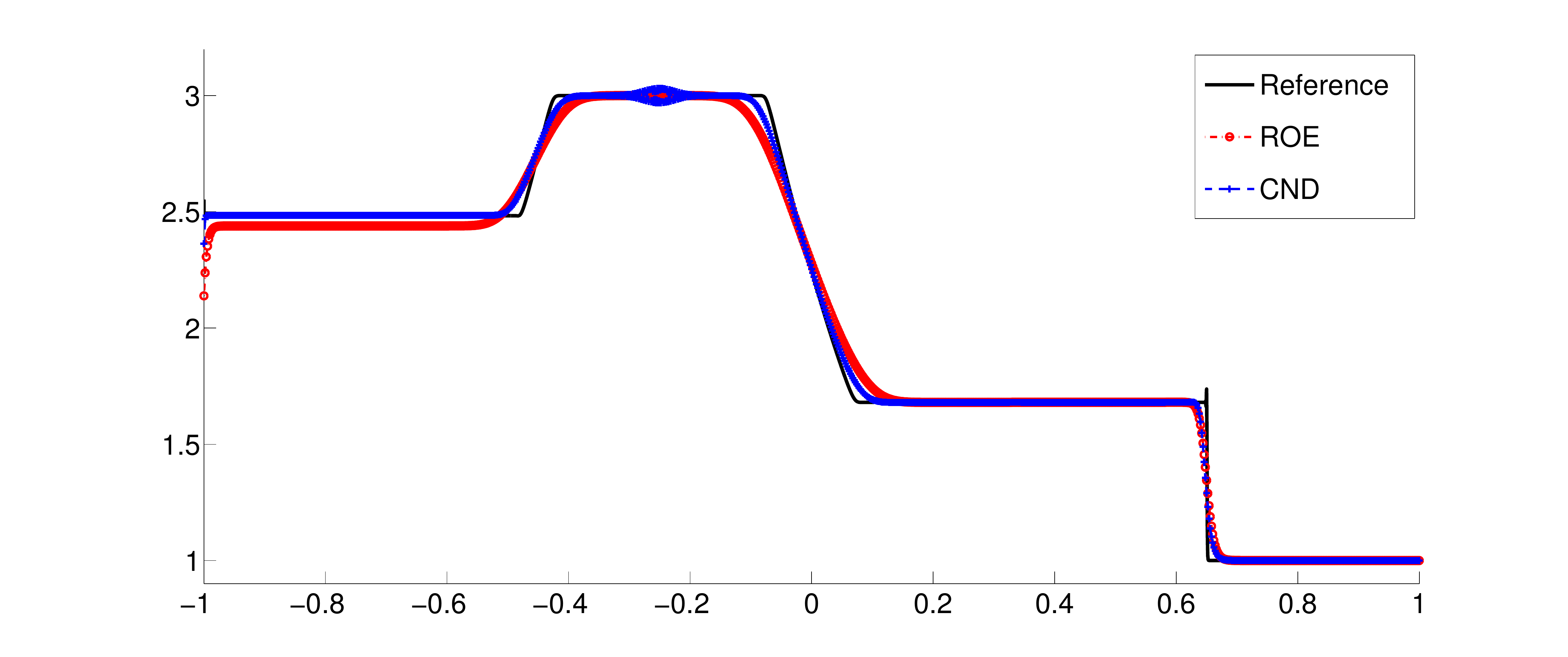}} 
\caption{Approximate solutions of the compressible Euler equations \eqref{eq:euler} with initial data \eqref{eq:eulinit} and boundary data \eqref{eq:eulbd} at time $t=0.25$. We compare the Roe and CND schemes on $1000$ mesh points with a reference solution of the compressible Navier-Stokes equations \eqref{eq:cns} with $\kappa = \nu = 10^{-5}$.}
\label{fig:74}
\end{figure}

We discretize the initial and boundary value problem for the compressible Euler equations \eqref{eq:euler} on the computational domain $[-1,1]$ with initial data \eqref{eq:eulinit} and Dirichlet data \eqref{eq:eulbd}. The results with the CND scheme and a standard Roe scheme at time $t=0.25$ are shown in Figure~\ref{fig:74}. We present approximate solutions, computed on a mesh of $1000$ points, for both schemes. Both the Roe and the CND schemes have converged at this resolution. As we are interested in approximating the physical-viscosity  solutions of the Euler equations, realized as a limit of the Navier-Stokes equations, we plot a reference solution computed on a mesh of $32000$ points of the compressible Navier-Stokes equations \eqref{eq:cns} with $\kappa = \nu = 10^{-5}$. The figure shows that the Roe scheme clearly converges to an incorrect solution near the left-hand boundary. This lack of convergence is most pronounced in the density variable. Similar results were also obtained with the standard Rusanov, HLL and HLLC solvers (see the book by LeVeque (2003) for a detailed description of these solvers).

On the other hand, the CND scheme converges to the physical-viscosity  solution. There are slight oscillations with the CND scheme as the numerical diffusion operator is singular. However, these oscillations do not impact on the convergence properties of this scheme. Although, the Roe scheme does not generate any spurious oscillations, yet it converges to an incorrect solution of the Euler equations. On the other hand, the CND scheme does converge to the physical-viscosity  solution (the Navier-Stokes limit) in spite of some spurious oscillations. It appears that the oscillations are the price that one has to pay in order to resolve the physical-viscosity  solution correctly.  Moreover, the CND scheme is slightly more accurate than the Roe scheme when both of them converge to the same solution (see near the interior contact).
\subsubsection{Schemes with well controlled dissipation}
The CND scheme of Mishra and Spinolo (2011) is a scheme with controlled dissipation in this case. It approximates shocks of moderate amplitude quite robustly as shown in the previous numerical experiments. However, there might be a deterioration of performance once the shock amplitude is very large. In this case, a scheme with well-controlled dissipation (WCD) can be readily designed using the construction involving the equivalent equation, as outlined in the previous section.

%==================================================================================

\section{Concluding remarks}
\label{sec:7} 

\subsection{Numerical methods with sharp interfaces} 
\label{sec:34-b} 

Another class of numerical methods for approximating small-scale dependent solutions to hyperbolic problems is based on front tracking schemes and random choice schemes. Although they may not be always suited for physical applications where no information is a priori known on the solutions and high-order accuracy is sought, these methods have some definite advantages. Most importantly, numerical shocks are not smeared and are represented as discontinuities and it is only the location of the shock and the left- and right-hand values that are approximated. 
Front-traching and random choice methods have been first used to establish the existence of small-scale dependent solutions and allows to achieve the following results: 
\bei 

\item Existence theory for the initial value problem for solutions with nonclassical  shocks by Amadori, Baiti, Piccoli, and LeFloch (1999), Baiti, LeFloch, and Piccoli (1999, 2000, 2001, 2004), LeFloch (2002), Laforest and LeFloch (2010, 2014). 

\item Existence theory for the initial value problem for nonconservative systems by LeFloch and Liu (1993), which was based on the nonconservative Riemann solver and the theory of nonconservative products develope by Dal Maso, LeFloch, and Murat (1990, 1995). 
 
\item Existence theory for  the initial value problem for hyperbolic systems of conservation laws by 
Amadori (1997), Amadori and Colombo (1997), 
Ancona and Marson (1999), 
Karlsen, Lie, and Risebro (1999),
Donadello and Marson (2007).  
% Amadori and Baiti. 

\eei 

In particular, this strategy was numerically implemented by Chalons and LeFloch (2003) (nonclassical solutions via the random choice scheme) 
and, in combination with a level set technique, by 
 Zhong, Hou, and LeFloch (1996), 
 Hou, Rosakis, and LeFloch (1999), 
and Merkle and Rohde (2006, 2007), who treated a nonlinear elasticity model with trilinear law in two spatial dimensions. 
As observed by Zhong, Hou and LeFloch, this model exhibits complex interface needles attached to the boundary, whose computation is 
 numerically very challenging since they are strongly small-scale dependent. 
In addition, methods combining differences and interface tracking were also developed  
which ensure that the interface is sharp and (almost) exactly propagated; see Boutin, Chalons, Lagouti\`ere, and LeFloch (2008).  

%-----------------------------------------------------------------------------------------------------------

\subsection{Convergence analysis} 
\label{sec:34-b} 

As described in the current review, the design and numerical implementation of robust and efficient methods for small-scale dependent schemes are 
well established by now. However, a complete theory, with convergence results, is available only for random choice and front tracking schemes; this theory 
encompasses nonclassical undercompressive shocks, nonconservative hyperbolic systems, and the boundary value problem. 
However, only partial results are available concerning the rigorous convergence analysis of the proposed schemes, 
and the interested reader may consult and built upon the references cited in the bibliography. 

%------------------------------------------------------------------------------------------------------------ 

\subsection{Perspectives}
\label{sec:72} 

In summary, our main guidelines in order to design efficient and robust numerical methods for small-scale dependent shock waves 
can be summarized as follows: 

\begin{itemize}

\item Standard finite difference or finite volume schemes fail to properly approximate weak solutions containing small-scale dependent shocks, and a {\sl well-controlled dissipation} requirement is necessary.

\item This lack of convergence to physically-relevant solutions arises with most nonlinear hyperbolic problems, including: 

\bei 

\item  non-genuinely nonlinear systems (when dispersive effects are present and may contribute to the dynamics of shocks), 

\item nonconservative hyperbolic systems (since products of discontinuous functions by measures are regularization-dependent), and 

\item boundary value problems (due to the formation of
possibly characteristic boundary layers). 

\eei

\item Numerous applications lead to such problems, especially: 

\bei 

\item the dynamics of viscous-capillary fluids, 

\item the dynamics of two-phase fluids such as liquid-vapor or liquid-solid or solid-solid mixtures.

\item the model of  magnetohydrodynamics with Hall effect incuded, which  leads to the most challenging system, on account of 
lack of genuine nonlinearity, lack of strict hyperbolicity, and presence of dispersive terms. 

\eei 

\item Although the classes of problems under consideration are quite distinct physical origin, a single source of difficulty was identified here, that is, the importance of properly computing the global effect of small-scale terms. Standard schemes for all these problems fail because the leading terms in the equivalent equation are different from the small-scale mechanisms of the underlying PDE.

\item Schemes with well-controlled dissipation have been developed in the past fifteen years, in order to 
accurately compute small-scale dependent shock waves. The key 
ingredient was to systematically design
 numerical diffusion operators that lead to the equivalent equation of the scheme matching the small-scale mechanisms 
of the underlying PDE at the leading order.

\item As the residual terms of the equivalent equation can become large with increasing shock strength, any scheme of a fixed order of accuracy will fail to converge to the correct solution for very large shocks, and this issue was addressed. 
 The proper notion of convergence for small-scale dependent shock waves is that the approximate solutions converge in $L^1$ as well as in terms of kinetic relations, as the order $p$ of the scheme is increased.

\item Computing kinetic functions, familes of paths, and admissible boundary sets 
is very useful in order to investigate the effects of the diffusion/dispersion ratio, regularization, order of accuracy of the schemes, 
 the efficiency of the schemes, as well as to make comparisons between several physical models.

\end{itemize}

%==================================================================================

\section*{Acknowledgements} 

The first author (PLF) was supported by a DFG-CNRS collaborative grant between France and Germany
entitled ``Micro-Macro Modeling and Simulation of Liquid-Vapor Flows'', by the Centre National de la Recherche Scientifique, and by the Agence Nationale de la Recherche. The second author (SM) was partially supported by ERC starting grant 306279 (SPARCCLE) from the European Research Council. 

%==================================================================================

\label{lastpage}
\end{document}